\documentclass[10pt,reqno]{article}

\usepackage[b5paper,top=1.2in,left=0.9in]{geometry}
\usepackage{verbatim}
\usepackage{amsthm}
\usepackage{amssymb}
\usepackage{amsmath}
\usepackage{graphicx}
\usepackage{appendix}
\usepackage{color}

\renewcommand{\d}{\mathrm{d}}
\newcommand{\D}{\mathrm{D}}
\renewcommand{\i}{\mathrm{i}}
\newcommand{\e}{\mathrm{e}}

\newtheorem{Thm}{Theorem}[section]
\newtheorem{Lem}[Thm]{Lemma}
\newtheorem{Prop}[Thm]{Proposition}
\newtheorem{Cor}[Thm]{Corollary}
\newtheorem{Rem}[Thm]{Remark}
\newtheorem{Def}[Thm]{Definition}
\newtheorem{Ex}[Thm]{Example}
\newtheorem{Fact}[Thm]{Fact}
\newtheorem{Nota}[Thm]{Notation}

\def\R{\mathbb{R}}
\def\Q{\mathbb{Q}}
\def\N{\mathbb{N}}
\def\C{\mathbb{C}}
\def\Z{\mathbb{Z}}

\def\to{\longrightarrow}

\def\cA{\mathcal{A}}
\def\cB{\mathcal{B}}
\def\cC{\mathcal{C}}
\def\cD{\mathcal{D}}
\def\cF{\mathcal{F}}
\def\cH{\mathcal{H}}
\def\cK{\mathcal{K}}

\def\cM{\mathcal{M}}
\def\cN{\mathcal{N}}
\def\cP{\mathcal{P}}

\def\cS{\mathcal{S}}
\def\cT{\mathcal{T}}
\def\cU{\mathcal{U}}

\def\cW{\mathcal{W}}
\def\a{\alpha}
\def\b{\beta}
\def\e{\epsilon}
\def\G{\Gamma}
\def\c{\gamma}
\def\D{\Delta}
\def\DD{\nabla}
\def\d{\delta}

\def\h{\theta}
\def\l{\lambda}
\def\L{\Lambda}
\def\S{\Sigma}
\def\s{\sigma}
\def\t{\tau}

\def\w{\omega}
\def\ze{\zeta}

\def\sl{\mathfrak{sl}}

\def\gl{\mathfrak{gl}}

\def\ox{\otimes}
\def\o+{\oplus}
\def\bo+{\bigoplus}
\def\x{\times}
\def\p[#1,#2]{\phi_{#1,#2}}
\def\til[#1]{\widetilde{#1}}
\def\what[#1]{\widehat{#1}}

\def\bC{\textbf{C}}

\def\bJ{\textbf{J}}

\def\bW{\textbf{W}}

\def\z[#1]{z_{#1}}

\def\oo{\infty}
\def\=>{\Longrightarrow}
\def\inj{\hookrightarrow}

\def\<{\langle}
\def\>{\rangle}
\def\corr{\longleftrightarrow}
\def\^{\wedge}
\def\+{\dagger}

\def\sub{\subset}
\def\inv{^{-1}}
\def\dis{\displaystyle}
\def\over[#1]{\overline{#1}}
\def\vec[#1]{\overrightarrow{#1}}
\def\mat[#1, #2]{\left[\begin{array}{ccccc}#1\end{array}\left|\begin{array}{c}#2\end{array}\right.\right]}
\def\xto[#1]{\xrightarrow{#1}}
\def\dd[#1,#2]{\frac{d#1}{d#2}}
\def\del[#1,#2]{\frac{\partial #1}{\partial #2}}
\def\Facts[#1]{\begin{Fact}\mbox{}\begin{itemize}#1\end{itemize}\end{Fact}}
\def\Notation[#1]{\begin{Nota}\mbox{}\begin{itemize}#1\end{itemize}\end{Nota}}
\def\Eqn[#1]{\begin{eqnarray*}#1\end{eqnarray*}}
\def\tab{\;\;\;\;\;\;}

\newcommand{\veca}[2][cccccccccccccccccccccccccccccccccccccccccc]{\left(\begin{array}{#1}#2 \\ \end{array} \right)}

\newcommand{\Eq}[1]{\begin{align}#1\end{align}}
\begin{document}

\title{Representation of the Quantum Plane, its Quantum Double and Harmonic Analysis on $GL_q^+(2,\R)$}

\author{  Ivan C.H. Ip\footnote{
          Kavli IPMU, The University of Tokyo,
          5-1-5 Kashiwanoha, Kashiwa,
	277-8583, Chiba, Japan
	\newline
          Email: ivan.ip@ipmu.jp}
          }

\date{\today}

\numberwithin{equation}{section}

\maketitle

\begin{abstract}
We give complete detail of the description of the GNS representation of the quantum plane $\cA$ and its dual $\hat{\cA}$ as a von Neumann algebra. In particular we obtain a rather surprising result that the multiplicative unitary $W$ is manageable in this quantum semigroup context. We study the quantum double group construction introduced by Woronowicz, and using Baaj and Vaes' construction of the multiplicative unitary $\bW_m$, we give the GNS description of the quantum double $\cD(\cA)$ which is equivalent to $GL_q^+(2,\R)$. Furthermore we study the fundamental corepresentation $T^{\l,t}$ and its matrix coefficients, and show that it can be expressed by the $b$-hypergeometric function. We also study the regular corepresentation and representation induced by $\bW_m$, and prove that the space of $L^2$ functions on the quantum double decomposes into the continuous series representation of $U_{q\til[q]}(\gl(2,\R))$ with the quantum dilogarithm $|S_b(Q+2i\a)|^2$ as the Plancherel measure. Finally we describe certain representation theoretic meaning of integral transforms involving the quantum dilogarithm function.
\end{abstract}

\newpage
\tableofcontents

\section{Introduction}\label{sec:intro}
The quantum plane $\cA_q$ as a Hopf algebra is generated by self-adjoint elements $A,B$ satisfying $AB=q^2BA$ and the coproduct
\Eq{\D(A)=A\ox A,\tab\D(B)=B\ox A+1\ox B} for $q=e^{\pi ib^2}$ with $|q|=1$. The coproduct reflects the fact that it is the quantum version of the classical $ax+b$ group, which is the group of affine transformations on the real line $\R$. Sometimes known as the quantum $ax+b$ group, it has been studied for example in \cite{PuW, WZ}. The main problem arose from the fact that $A$ and $B$ are realized as unbounded operators, and that $\D(B)$ is in general not self-adjoint. This poses quite some problems in the well-definedness of the algebra on the $C^*$-algebraic level. 

A class of well-behaved ``integrable" representations is studied in \cite{FK, Sch} where instead we impose the condition that both $A,B$ are \emph{positive} self-adjoint. In this case the operators can be realized as the``canonical representation":
\Eq{\label{Hirr}A = e^{2\pi bs},\tab B=e^{2\pi bp},}
where they act on $\cH_{irr}:=L^2(\R)$ as unbounded operators. This simplifies a lot of functional analytic problems, because now all the operators considered are positive essentially self-adjoint, with the help of certain transformations that can be carried out by the quantum dilogarithm function $g_b(x)$. On the other hand, with positivity, we can define using functional calculus a wide class of functions on $\cA_q$, and in this way we can define the $C^*$-algebra $\cA:=C^\oo(\cA_q)$ of ``functions vanishing at infinity" for $\cA_q$, which is expressed using an integral transform
\Eq{f:= \iint f(s,t)A^{ib\inv s}B^{ib\inv t} dsdt,} where $f(s,t)$ satisfies certain analytic properties. It turns out that this definition encodes the modular double counterpart as well. By definition, the other half $\cA_{\til[q]}$, first introduced by Faddeev \cite{Fa}, is generated by $\til[A],\til[B]$ satisfying $\til[A]\til[B]=\til[q]^2\til[B]\til[A]$ with the same coproduct, where $\til[q]=e^{\pi ib^{-2}}$, and they are related to the quantum plane by $\til[A]=A^{1/b^2}, \til[B]=B^{1/b^2}$. This $b\longleftrightarrow b\inv$ duality is manifest in the definition of $\cA$ and subsequently present in all later calculations involved.

The first main result of this paper is the derivation of the Haar functional on $\cA$ (cf. Theorem \ref{Haar}). It can be expressed simply by
\Eq{h\left( \iint f(s,t)A^{ib\inv s}B^{ib\inv t} dsdt\right)=f(0,iQ)} where $Q:=b+b\inv$. 
This comes as a surprise as there is no classical Haar measure on the $ax+b$ semigroup due to the lack of inverse. With $\cA$ a $C^*$-algebra equipped with a left invariant Haar weight, we can carry out the so-called Gelfand-Naimark-Segal (GNS)-construction, which essentially represents $\cA$ naturally on a Hilbert space $\cH$ induced by its own multiplication. In this paper we put $\cA$ in the context of the theory of locally compact quantum group in the von Neumann setting \cite{KV1,KV2}, in which the modular theory (Tomita-Takesaki's Theory) for von Neumann algebra is studied, and the main ingredient, the multiplicative unitary $W$, can be defined as a unitary operator on $\cH\ox\cH$. 

The multiplicative unitary $W$ is a unitary operator on the Hilbert space $\cH\ox\cH$ that satisfies the pentagon equation
\Eq{W_{23}W_{12}=W_{12}W_{13}W_{23}.}
For every quantum group a multiplicative unitary can be constructed using the coproduct, however, not every multiplicative unitary is related to a quantum group. In \cite{W} Woronowicz introduced the notion of \emph{manageability} that describes a class of well-behaved multiplicative unitary. The second main result of this paper is that this multiplicative unitary $W$ obtained for the quantum plane above is in fact manageable (cf. Theorem \ref{manageable}). This is rather striking because in \cite{WZ} it is mentioned that manageability is the property that distinguishes quantum groups from quantum semigroups, however as noted above, we have been restricting ourselves to positive operators. The main reason is the fact that the GNS construction provides us with a ``bigger" Hilbert space, so that there is more freedom of choice for the operators to satisfy the manageability condition. As a by-product of this discrepancy, it turns out that we obtain a new transformation rule for the quantum dilogarithm function (cf. Proposition \ref{32}) that is not available in the literature.

The motivation for the study of the quantum plane comes from the quantum double group construction introduced in \cite{PW} for compact quantum groups, which is the dual version of the Drinfeld double construction. In \cite{Pu} it is shown that the quantum double construction of the quantum $az+b$ group gives rise to $GL_q(2,\C)$ for certain root of unity $q$. It turns out that this construction can be carried over to locally compact quantum groups. In order to carry out a similar recipe for the quantum plane for general $q$ with $|q|=1$, it is necessary to define the dual space $\hat{\cA}$ in the same setting as $\cA$ on the $C^*$-algebraic level. The third main result of the paper is the derivation of the non-degenerate pairing between $\cA$ and $\hat{\cA}$ on the $C^*$-algebraic level, which remarkably involves the quantum dilogarithm function in place of the $q$-factorial. Following \cite{KV1} we can then describe the GNS-representation for $\hat{\cA}$ acting on the same GNS space $\cH$ for $\cA$, its modular theory and its multiplicative unitary as before.

After defining the dual space, we can apply the quantum double group construction and obtain a new algebra $\cD(\cA)$ generated by $\veca{\z[11]&\z[12]\\\z[21]&\z[22]}$ which is precisely the $GL_q^+(2,\R)$ quantum (semi)group as Hopf algebra, where the generators are again restricted to positive self-adjoint elements. We note also that the relations for $GL_q^+(2,\R)$ is a special case for the two parameter deformations $GL_{p,q}(2,\C)$ observed in \cite{DMMZ,SWZ} where $p:=1,q:=q^2$. On the other hand, the relations involved unmistakably resemble the quantum Minkowski spacetime relations defined in \cite{FJ}, but this time in the non-compact setting where $|q|=1$ and the variables $\z[ij]$ are positive self-adjoint. Therefore the quantum plane can be seen as a building block towards this ``split quantum Minkowski spacetime" $\cM_q$. Hence from the properties of $\cD(\cA)$, the quantum Minkowski spacetime can be easily extended to the $C^*$-algebraic and von Neumann algebraic level, and hence $L^2(GL_q^+(2,\R))$ is well-defined. Motivated from the representation theory of classical $SL^+(2,\R)$, we define the matrix coefficients $T_{s,\a}^{\l,t}(z)$ for the fundamental corepresentation of $\cD(\cA)$ (Definition \ref{fundT}). For $t=\l$ it is explicitly given by
\Eq{T_{s,\a}^\l(z)=\frac{1}{2\pi b}\left(\int_0^\oo (x\z[11]+\z[21])^{b\inv(l-i\a)}(x\z[12]+\z[22])^{b\inv(l+i\a)} x^{b\inv(-l+is)} \frac{dx}{x}\right)N^{\frac{Q}{2b}},} which can be seen to be a generalization of the matrix coefficient in the compact case \cite{FL}. On the other hand, a closed form expression (cf. Corollary \ref{matrixFb}) using the $b$-hypergeometric function $F_b$ is found, which clearly gives a quantum analogue of the classical formula to the matrix coefficients for $SL(2,\R)$.

 To complete the description of $\cD(\cA)$ in the von Neumann setting, it is necessary to compute its multiplicative unitary $\bW_m$. A derivation is obtained from the construction given in \cite{BV,K} for general quantum groups where $\bW_m$ is obtained as a product of 6 $W$'s of the base quantum group and its dual. Restricting to the present simpler setting and using the properties of $W$ we simplify the expression to just 4 $W$'s. Explicitly they are given by
\Eq{\bW=W_{13}V_{32}''\hat{W}_{24}V_{32}^*}
(cf. Proposition \ref{Wmsimple}).
This is encouraging since it has been shown in certain context that the $R$ matrix (satisfying the hexagon or Yang Baxter equation) can be expressed in 4 $R$'s \cite[Ex 7.3.3]{Ma}, or 4 $W$'s \cite{Ka1} from the smaller group using instead the Heisenberg Double construction.

With the new multiplicative unitary constructed, it is straightforward, despite tedious, to construct both the left and right regular corepresentations for $\cD(\cA)$. Furthermore, I. Frenkel \cite{F} has constructed a non-degenerate pairing between $\cM_q$ with $U_q(\gl(2,\R))$, hence from the corepresentation we can derive the fundamental representation $\cP_{\l,t}$ from $T^{\l,t}$, the left and right regular representation from $\bW_m$ of $U_q(\gl(2,\R))$ as well as its modular double on $L^2(GL_q^+(2,\R))$ by the pairing. The final main result (cf. Theorem \ref{main}) of this paper is the quantum ``Peter Weyl theorem"
\Eq{L^2(GL_q^+(2,\R))\simeq \int_\R^\o+ \int_{\R_+}^\o+ \cP_{\l,s/2}\ox \cP_{\l,-s/2} |S_b(Q+2i\l)|^2 d\l ds}
as a representation of $U_{q\til[q]}(\gl(2,\R))_L\ox U_{q\til[q]}(\gl(2,\R))_R$, which is a generalization of the statement for $U_q(\sl(2,\R))$ announced in \cite{PT1}, however the details of the proof are never published. The remarkable fact is that the quantum dilogarithm appears as the Plancherel measure for $L^2(GL_q^+(2,\R))$, which comes from the spectral analysis of the Casimir operators. In short, it states that the operator
\Eq{\bC=e^{2\pi bx}+e^{-2\pi bx}+e^{-2\pi bp}}
can be diagonalized as multiplication operator
\Eq{\bC=e^{2\pi b\l}+e^{-2\pi b\l},}
acting on the Hilbert space $L^2(\R_+, |S_b(Q+2i\l)|^2d\l)$, using certain eigenfunction transform (cf. Theorem \ref{Casimir}).
Furthermore, we note that only the fundamental series $\cP_{\l}$ appears in the decomposition of the regular representation. Called the ``self-dual" principal series representation, they are known \cite{PT1} to be closed under tensor product, which is rather interesting since the same does not hold in the classical group setting \cite{Puk}.

As a corollary, we have expressed the multiplicative unitary canonically as a direct integral of the fundamental representations (cf. Corollary \ref{Wdecompose}):
\Eq{\bW_m = \int_\R\int_{\R_+}^{\o+} T^{\l,t/2} d\mu(\l) dt,}
which generalizes the canonical definition for the compact quantum group given in \cite{PW}. An interesting problem will be to investigate its classical limit. As we know, the functions on the full group $SL(2,\R)$ contains both the continuous and discrete series representation, hence the above results may give an insight to the decomposition of the functions on the \emph{positive} semigroup and its decomposition, and distinguish the principal continuous series as the fundamental component. Furthermore, we also know that classically $SL(2,\R)\simeq SU(1,1)$, but the behaviors change drastically in the quantum level. Therefore another interesting problem will be a comparison with several known harmonic analysis on the quantum $SU(1,1)$ group in the operator algebraic setting, where the discrete series and the so-called \emph{strange series} also play a role \cite{GKK,KK}.

Another aspect of this paper is the study of several properties of the quantum dilogarithm function $G_b$. The quantum dilogarithm function played a prominent role in this quantum theory. This function and its many variants are being studied \cite{Gon, KaN, Ru, Vo} and applied to vast amount of different areas, for example the construction of the '$ax+b$' quantum group by Woronowicz et al. \cite{PW, WZ}, the harmonic analysis of the non-compact quantum group $U_q(\sl(2,\R))$ and its modular double \cite{BT,PT1,PT2}, the $q$-deformed Toda chains \cite{KLS}, current algebra and Virasoro algebra on the lattice \cite{FV}, and hyperbolic knot invariants \cite{Ka2}. Recently attempts have also been made to cluster algebra \cite{FG, KaN} and quantization of the Teichm\"{u}ller space \cite{CF, FK, Ka3}. One of the important properties of this function is its invariance under the duality $b\leftrightarrow b\inv$ that helps encoding the details of the modular double in $\cA$, and also relates, for example, to the self-duality of Liouville theory \cite{PT1}. The classical limit and several relations of $G_b$ are studied in the previous paper \cite{Ip}. In this paper, we give a proof of the important relation (cf. Corollary \ref{delta})
\Eq{\lim_{\e\to 0}\int_C\frac{G_b(Q-\e)G_b(iz)}{G_b(Q+iz-\e)} f(z)dz = f(0)}
with suitable contour $C$, which is important to make sense rigorously of certain calculations involving the Haar functional, as well as the non-degenerate pairing needed to obtain the regular representations. This statement and its variants is often assumed, for example in \cite{Vo}, but is never proved. 

As noted above, there are various version of the quantum dilogarithm. Their relations and visualization of their graphs are presented in \cite{Ip2}. We choose this particular definition $G_b, g_b$ by Teschner \cite{BT} for various reasons. First of all, it gives a closed form expression for the $q$-Binomial Theorem (Lemma \ref{qbi}), so that it gives precisely the quantum analogue of the classical Gamma function $\G(x)$, see \cite{Ip}. Furthermore, various transformations such as the Tau-Beta Theorem (Lemma \ref{tau}), the ``45-relation"
 (Lemma \ref{45}) and the ``69-relation" (Proposition \ref{69}) can be written in closed form without any extra constants and exponentials. Finally the multiplicative unitary constructed above can be expressed in closed form by $g_b$ in a very simple way. Another popular version $G(a_+,a_-;z)$ is given by Ruijsenaars \cite{Ru} and their relation is given by (cf. Lemma \ref{GbG})
\Eq{G_b(z)=G(b,b\inv;iz-\frac{iQ}{2})e^{\frac{\pi i}{2}z(z-Q)}.}
We have chosen the pair $(a_+,a_-)$ to be $(b,b\inv)$ which demonstrates the self-duality of the Liouville theory in the physical literature. We believe similar treatments in our paper can be made for a general pair with $q=e^{\pi i\frac{a_+}{a_-}}$. 

It is well-known that by studying the matrix coefficients of various classical matrix groups, we recover certain functional relations between different special functions, for example the hypergeometric functions ${}_mF_n$. Various examples can be found e.g. in \cite{Vi}. In this quantum setting, it is not surprising that all these various calculations involving the quantum dilogarithm will provide us certain representation theoretic meaning of their relations by integral transforms. By applying the classical limit for $G_b$ obtained in \cite{Ip} (cf. Theorem \ref{classicallimit}), we recover Barnes' first and second lemma for the Gamma function. Furthermore, a relation involving $G_b$ that is not commonly used is observed (cf. Proposition \ref{32}):
\Eq{\int_C G_b(\a+i\t)G_b(\b-i\t)G_b(\c-i\t)e^{-2\pi i(\b-i\t)(\c-i\t)}d\t = G_b(\a+\c)G_b(\a+\b)}
These will be briefly discussed in the last section.

As a side note, the frequent use of the term ``modular" in this paper requires some clarification. The modular double introduced by Faddeev refers in the general case to the quantum groups related by the transformation of the modular group for the complex parameter 
\Eq{q=e^{\pi i\t} \mapsto \til[q]=e^{\pi i\frac{a\t+b}{c\t+d}},\tab \veca{a&b\\c&d}\in SL(2,\Z),}
while the modular theory of Tomita-Takesaki refers in the classical case to the modularity of the Haar measure, for example the modular function $\D$ that relates the left translates of the right Haar measure
\Eq{\mu(g\inv A)=\D(g)\mu(A),} where $g\in G$ and $A$ is a Borel set in $G$.

The present paper is organized as follows. In Section 2 we recall several technical results and motivations that is needed in later sections, including the Mellin transform, the GNS-construction and the significance of the multiplicative unitary. We also collect several results concerning the calculations involving $q$-commuting variables. In Section 3 we recall the definition and properties of the quantum dilogarithm function $G_b$, and describe its integral transformation formula. In Section 4 we define the quantum plane $\cA$ in the Hopf algebra and $C^*$-algebraic level, the modular double, and describe completely its GNS construction and the multiplicative unitary. Section 5 deals with the dual space $\hat{\cA}$, we derive a non-degenerate pairing and obtain its GNS construction on the same space as $\cA$. In Section 6 we develop two useful transformation that shed light to the action of the multiplicative unitary, as well as the action of the quantum plane on a Hilbert space. In Section 7 we carry out the quantum double group construction, study its fundamental corepresentation, and express explicitly the new multiplicative unitary, and the corepresentation induced by it. In Section 8 we look at the dual picture and obtain a pairing with $U_q(\gl(2,\R))$, derive the regular representations, and prove the main theorem on the decomposition of $L^2(GL_q^+(2,\R))$ into principal continuous series of $U_q(\gl(2,\R))$. Finally the last section discuss certain integral transformations of $G_b$ arising from representation theoretic calculations.

\textbf{Acknowledgments.} I would like to thank my advisor Professor Igor Frenkel for suggesting the current project and providing useful insights to the problems. I would also like to thank L. Faddeev for pointing out several useful references, and to J. Teschner and Hyun Kyu Kim for helpful discussions.

\section{Preliminaries}\label{sec:Pre}
In this section, we recall several technical results that will be needed in the description of the quantum plane. We will first remind the definitions and properties of the Mellin transform which serves as a motivation to define the quantum plane algebra in Section \ref{sec:QPlane}. Then we describe the details of weights and multiplier algebra of a $C^*$-algebra, and its GNS construction and modular theory in the von Neumann setting. Next we explain the significance of the multiplicative unitary that is important in the study of locally compact quantum groups. Finally we collect several technical results concerning the calculations involving $q$-commuting variables.


\subsection{Mellin transform}\label{sec:Pre:Mellin}
In this subsection, let us recall the Mellin transform of a function and its properties. 

\begin{Thm}Let $f(x)$ be a continuous function on the half line $0<x<\oo$. Then its Mellin transform is defined by
\Eq{\phi(s):=(\cM f)(s)=\int_0^\oo x^{s-1} f(x)dx,} whenever the
integration is absolutely convergent for $a<Re(s)<b$. By the Mellin
inversion theorem, $f(x)$ is recovered from $\phi(s)$ by
\Eq{f(x):=(\cM\inv\phi)(x)=\frac{1}{2\pi}\int_{c-i\oo}^{c+i\oo}x^{-s}\phi(s)ds,}
where $c\in\R$ is any values in between $a$ and $b$.
\end{Thm}

We have the following analyticity theorem \cite{PK}:
\begin{Thm}\label{MTpole}
(Strip of analyticity) If $f(x)$ is a locally integrable function on $(0,\oo)$ such that it has decay property:
\Eq{\label{growth}f(x)=\left\{\begin{array}{cc}O(x^{-a-\e})&x\to 0^+\\O(x^{-b+\e})&x\to +\oo \end{array}\right.} 
for every $\e>0$ and some $a<b$, then the Mellin transform defines an analytic function $(\cM f)(s)$ in the strip
$$a<Re(s)<b.$$

(Analytic continuation) Assume $f(x)$ behaves algebraically for $x\to 0^+$, i.e. 
\Eq{f(x)\sim \sum_{k=0}^\oo A_k x^{a_k},} 
where $Re(a_k)$ increases monotonically to $\oo$ as $k\to\oo$. Then the Mellin transform $(\cM f)(s)$ can be
analytically continued into $Re(s)\leq a=-Re(a_0)$ as a meromorphic function with simple poles at the points $s=-a_k$
with residue $A_k$.

A similar analytic property holds for the continuation to the right half plane.

(Growth) If $f(x)$ is a holomorphic function of the complex variable $x$ in the sector $-\a<\arg x<\b$ where
$0<\a,\b\leq \pi$, and satisfies the growth property \eqref{growth} uniformly in any sector interior to the above
sector, then $(\cM f)(s)$ has exponential decay in $a<Re(s)<b$ with
\Eq{(\cM f)(s)=\left\{\begin{array}{cc}O(e^{-(\b-\e)t})&t\to+\oo\\O(e^{(\a-\e)t})&t\to -\oo \end{array}\right.} 
for any $\e>0$ uniformly in any strip interior to $a<Re(s)<b$, where $s=\s+it$.

(Parseval's Formula) 
\Eq{\int_0^\oo f(x)g(x)x^{z-1}dx =\frac{1}{2\pi i}\int_{c-i\oo}^{c+i\oo} (\cM f)(s)(\cM g)(z-s)ds,} 
where $Re(s)=c$ lies in the common strip for $\cM f$ and $\cM g$. In
particular we have 
\Eq{\int_0^\oo |f(x)|^2 dx =\frac{1}{2\pi}\int_{-\oo}^\oo |(\cM f)(\frac{1}{2}+it)|^2dt.}
\end{Thm}


\subsection{Weight and multiplier algebra of a $C^*$-algebra}\label{sec:pre:C*}
In this subsection, we recall the definition of weights on a $C^*$-algebra, and the language of multiplier algebra. Most of the notions are adopted from \cite{KV1,T}. Let $\cA$ be a $C^*$-algebra and $\cA^+$ its positive self-adjoint elements.
\begin{Def} A weight on a $C^*$-algebra $\cA$ is a function $\phi: \cA^+\to [0,\oo]$ such that
\begin{eqnarray}
\phi(x+y)&=&\phi(x)+\phi(y),\\
\phi(rx)&=&r\phi(x),
\end{eqnarray}
for $x,y\in \cA^+, r\in \R_{>0}$.
\end{Def}
\begin{Def} Given a weight $\phi$ on $\cA$, we define
\begin{eqnarray}
\cM_{\phi}^+&=&\{a\in \cA^+: \phi(a)<\oo\},\\
\cN_{\phi}&=&\{a\in \cA: \phi(a^*a)<\oo\},\\
\cM_{\phi}&=&\mbox{span}\{y^*x: x,y\in \cN_{\phi}\}.
\end{eqnarray}
Then it is known that $\cM_{\phi}^+=\cM_{\phi}\cap \cA^+$, and that $\phi$ extends uniquely to a map $\cM_{\phi}\to \C$.
A weight is called faithful iff $\phi(a)=0\=> a=0$ for every $a\in \cA^+$.
\end{Def}

Next we recall a useful notion of a multiplier algebra. Let $\cB(\cH)$ be the algebra of bounded linear operators on a Hilbert space $\cH$.

\begin{Def}If $\cA\sub \cB(\cH)$ as operators, then the multiplier algebra $M(\cA)$ of $\cA$ is the $C^*$-algebra of operators
\Eq{M(\cA)=\{b\in \cB(\cH): b\cA \sub \cA, \cA b\sub \cA\}.}
In particular, $\cA$ is an ideal of $M(\cA)$.
\end{Def}
\begin{Ex} Important examples include
\Eq{M(\cK(\cH))=\cB(\cH),}
where $\cK(\cH)$ are compact operators on $\cH$, and
\Eq{M(C_0(X))=C_b(X),}
where $X$ is a locally compact Hausdorff space, $C_0(X)$ is the algebra of $\C$-valued functions on $X$ vanishing at infinity equipped with the sup-norm, and $C_b(X)$ is the $C^*$-algebra of all bounded continuous functions on $X$.
\end{Ex}

\begin{Prop}
Let $\cA$ and $\cB$ be $C^*$-algebras. A homomorphism $\phi:\cA\to M(\cB)$ is called non-degenerate if the linear span of $\phi(\cA)\cB$ and of $\cB\phi(\cA)$ are both equal to $\cB$. Then $\phi$ extends uniquely to a homomorphism $M(\cA)\to M(\cB)$. In particular, by taking $\cB=\C$, every weight $\w$ on $\cA$ has a unique extension to $M(\cA)$.
\end{Prop}

Using the notion of a multiplier algebra, the concept of a multiplier Hopf algebra is introduced in \cite{VD} (see also \cite{T}). In particular, the coproduct $\D$ of $\cA$ will be a non-degenerate homomorphism $\D:\cA\to M(\cA\ox \cA)$. The coassociativity is well-defined from the proposition above.

Given a multiplier Hopf * algebra $\cA$ with coproduct $\D$, we can define left and right invariance of a functional.
\begin{Def}\label{leftinv} A linear functional $h$ on $\cA$ is called a left invariant Haar functional if it satisfies
\Eq{(1\ox h)(\D x) = h(x)\cdot 1_{M(\cA)},}
where $1_{M(\cA)}$ is the unital element in $M(\cA)$. 

Similarly, a right invariant Haar functional satisfies
\Eq{(h\ox 1)(\D x)=h(x)\cdot 1_{M(\cA)}.}
\end{Def}

\subsection{GNS representation and Tomita-Takesaki's theory}\label{sec:pre:GNS}
Let us recall the main objects in the study of GNS representation of a $C^*$-algebra (see e.g. \cite{KV1}).

\begin{Def} A Gelfand-Naimark-Segal (GNS) representation of a $C^*$-algebra $\cA$ with a weight $\phi$ is a triple
$$(\cH, \pi, \L),$$
where $\cH$ is a Hilbert space, $\L:\cA\to\cH$ is a linear map, and $\pi:\cA\to\cB(\cH)$ is a representation of $\cA$ on $\cH$ such that $\L(\cN)$ is dense in $\cH$, and
\begin{eqnarray}
\pi(a)\L(b)&=&\L(ab) \tab\forall a\in\cA,b\in \cN,\\
\<\L(a),\L(b)\>&=&\phi(b^*a)\tab\forall a,b\in\cN,
\end{eqnarray}
where $\cN=\{a\in\cA: \phi(a^*a)<\oo\}$.
\end{Def}

The Tomita-Takesaki's Theory \cite{Ta,T} provides a detailed description of the GNS-construction in the von Neumann algebraic setting:
\begin{Def} Giving a weight $\phi$ that determines the Hilbert space structure by the GNS construction, the operator
\Eq{\label{T}T: x\to x^*} is closable and has a polar decomposition as 
\Eq{T=J\DD^{\frac{1}{2}},}
where $J$ is called the modular conjugation, and $\DD$ is called the modular operator.
\end{Def}

A very important property of $J$ is the following
\begin{Thm}[Murray-von Neumann]\label{MvN} Let $M$ be the completion of $\cA\sub \cB(\cH)$ in the weak operator topology as a von Neumann algebra. Then considering $M\sub \cB(\cH)$, we have
\Eq{JM J=M',}
where $M'$ is the commutant of $M$.
\end{Thm}

\begin{Def} For $x\in\cA$, the operator
\Eq{\s_t^{\phi}(x) := \DD^{it}x\DD^{-it}} is called the modular automorphism group.

We have $\phi(a^*b) = \phi(b\s_0^\phi(a^*))$

\end{Def}

On the Hopf * algebra level, we have the following properties:

\begin{Prop}
The antipode $S$ has a polar decomposition \Eq{S=\t_{-i/2}\circ R,}
where $\t_{-i/2}$ denotes the analytic generator of $(\t_t)_{t\in\R}$, called the \emph{scaling group}, which is a group of automorphisms of $M$, and $R$, called the \emph{unitary antipode}, is an anti-automorphism of $M$. 

We have $R^2=1$ and $S^2=\t_{-i}$.

\end{Prop}

Finally some properties concerning the left and right invariant Haar functional:
\begin{Prop}
The left invariant Haar functional $\phi$ and right invariant Haar functional $\psi$ are related by $\psi=\phi \circ R$.

We have \Eq{\phi\circ\t_t=\nu^{-t}\phi,\tab \psi\circ\s_t^{\phi}=\nu^{-t}\psi,}where $\nu>0$ is called the scaling constant.

Furthermore there exists an element $\d\in \cA$ such that $$\s_t^\phi(\d)=\nu^t\d,$$ where $\d\in \cA$ is called the modular element.
\end{Prop}

\begin{Prop} \label{GNSright}The map $\L_R: \cA\to \cH$ defined by
\Eq{\L_R(a):=\L_L(a\d^{-\frac{1}{2}}),} gives the GNS representation of $\cA$ with the right Haar functional on the same space $\cH$, where $\d$ is the modular element defined above.
\end{Prop}

\subsection{The multiplicative unitary}\label{sec:pre:mu}
Multiplicative unitaries are fundamental to the theory of quantum groups in the setting of $C^*$-algebras and von Neumann algebras, and to generalization of Pontryagin duality. In particular, a multiplicative unitary encodes all the structure maps of a quantum group and its dual. A very good exposition is given in \cite{T}.

First let us define the leg notation.
\begin{Def}Let $\cH$ be a Hilbert space and $W\in \cB(\cH\ox\cH)$ be a bounded operator. Then we define
$W_{ij}\in\cB(\cH^{\ox k})$ by letting $W$ acts on the factors at the position $i,j$. In particular, the operators
 $W_{12},W_{23},W_{13}\in \cB(\cH\ox\cH\ox\cH)$ are given by the formulas
\Eq{T_{12}:=T\ox Id_H,\tab T_{23}:=Id_H\ox T,}
\Eq{T_{13}:=\S_{12}T_{23}\S_{12}=\S_{23}T_{12}\S_{23},}
where $\S\in\cB(\cH\ox\cH)$ denotes the flip $f\ox g\mapsto g\ox f$.
\end{Def}

\begin{Def} Let $\cH$ be a Hilbert space. A unitary operator $W\in \cB(\cH\ox \cH)$ is called a multiplicative unitary if it satisfies the pentagon equation
\Eq{W_{23}W_{12}=W_{12}W_{13}W_{23},}
where the leg notation is used.
\end{Def}

Given a GNS representation $(\cH,\pi, \L)$ of a locally compact quantum group $\cA$, we can define a unitary operator \Eq{\label{WW}W^*(\L(a)\ox \L(b)) = (\L\ox \L)(\D(b)(a\ox 1)).}

It is known that $W$ is a multiplicative unitary \cite[Thm 3.16, Thm 3.18]{KV1}, and the coproduct on $\cA$ defining it can be recovered from $W$:
\begin{Prop} Let $x\in\cA \inj \cB(\cH)$ as operator. Then 
\Eq{W^*(1\ox x)W = \D(x)} as operators on $\cH\ox \cH$.
\end{Prop}
\begin{proof} For $x, f,g\in \cA$, we have $x\cdot \L(f)=\L(xf)$, hence
\begin{eqnarray*}
W^*((1\ox x)\cdot (\L(f)\ox \L(g)))&=&W^* (\L(f)\ox \L(xg))\\
&=&(\L\ox \L)(\D(xg)(f\ox 1))\\
&=&(\L\ox \L)(\D(x)\D(g)(f\ox 1))\\
&=&\D(x)\cdot (\L\ox \L)(\D(g)(f\ox 1))\\
&=&\D(x)W^*(f\ox g).
\end{eqnarray*}
\end{proof}

As a motivation, let us consider an example involving classical group:
\begin{Ex}\cite[Ex 7.1.4]{T}\cite{W} Let $G$ be a locally compact group with right Haar measure $\l$. Then the operator 
\Eq{(W f)(x,y):=f(xy,y)}
is a multiplicative unitary in $\cB(L^2(G,\l)\ox L^2(G,\l))=\cB(L^2(G\x G,\l\x\l))$. The pentagon equation for $W$ is equivalent to the associativity of the multiplication in $G$. Indeed, for $f\in L^2(G\x G\x G, \l \x \l \x \l)\simeq L^2(G,\l)^{\ox 3}$ and $x,y,z\in G$, we have
\begin{eqnarray*}
(W_{23}W_{12}f)(x,y,z)&=&f(x(yz),yz,z),\\
(W_{12}W_{13}W_{23}f)(x,y,z)&=&f((xy)z,yz,z).
\end{eqnarray*}
\end{Ex}

An important property of $W$ is that it encodes the information of the dual quantum group $\hat{\cA}$. By definition, $\hat{\cA}$ is the closure of the linear span of \Eq{\{(\w\ox 1)W: \w\in \cB(\cH)^*\}\sub\cB(\cH),} and we actually have 
\Eq{W\in M(\cA\ox \hat{\cA}),} where $M$ stands for the multiplier algebra.

Similarly, the multiplicative unitary for the dual $\hat{\cA}$ is given by $\hat{W}:=W_{21}^*$, hence it is known from the Pontryagin duality that $\hat{W}\in M(\hat{\cA}\ox \cA)$ as well. From the pentagon equation, we then obtain:

\begin{Cor}\label{DWW}
As an element $W\in M(\cA\ox \hat{\cA})$, together with the pentagon equation, we have 
\Eq{\label{DW}(\D\ox 1)W=W_{13}W_{23},}
\Eq{(1\ox \hat{\D})W=W_{13}W_{12}.}
\end{Cor}

Finally let us define the notion of ``manageability" introduced by Woronowicz \cite{W, WZ} that describes a class of well-behaved multiplicative unitary. It is shown that any manageable multiplicative unitary gives rise to a quantum group on the $C^*$-algebraic level.
\begin{Def}\label{manageable}
A multiplicative unitary $W$ is manageable if there exists a positive self-adjoint operator $Q$ acting on $\cH$ and a unitary operator $\til[W]$ acting on $\over[\cH]\ox \cH$ such that $\ker(Q)=\{0\}$,
\Eq{(Q\ox Q)W = W(Q\ox Q),} 
and
\Eq{\<x\ox u, W(z\ox y)\>_{\cH\ox \cH} = (\over[z]\ox Qu, \til[W](\over[x]\ox Q\inv y)\>_{\over[\cH]\ox \cH}}
for any $x,z\in \cH, y\in D(Q\inv), u\in D(Q)$.
\end{Def}
Here $\over[\cH]$ denote the complex conjugate of $\cH$, so that the map $x\in \cH \to \over[x]\in \over[\cH]$ is an anti-unitary map. The inner product on $\over[\cH]$ is given by $\<\over[x],\over[y]\>_{\over[\cH]}:=\<y,x\>_\cH$.

\subsection{$q$-commuting operators}\label{sec:pre:qcom}
Throughout the paper, we will consider positive operators $A,B$ satisfying $AB=q^2BA$ for $q=e^{\pi ib^2}, |q|=1$. 

We will realize the operators using the canonical pair $A=e^{2\pi bx}, B=e^{2\pi bp}$ where $p=\frac{1}{2\pi i}\frac{d}{dx}$. Then it is well-known that 
\begin{Prop} Both $A,B$ are positive unbounded operators on $L^2(\R)$ and they are essentially self-adjoint. The domain for $A$ is given by 
\Eq{\cD_A=\{f(x)\in L^2(\R): e^{2\pi bx}f(x) \in L^2(\R)\}}
and the domain for $B$ is given by the Fourier transform of $\cD_A$.
\end{Prop}
Hence we can apply functional calculus and obtain various functions in $A$ and $B$. In particular, for any function defined on $x>0\in\R$ such that $|f(x)|=1$, $f(A)$ will be a unitary operator.

In this paper, we will be using a dense subspace $\cW\sub L^2(\R)$ introduced in \cite{Gon, Sch}.
\begin{Def}\label{WDef}The dense subspace $\cW\sub L^2(\R)$ is the linear span of functions of the form
\Eq{e^{-\a x^2+\b x}P(x),} where $\a,\b\in\C$ with $Re(\a)>0$, and $P(x)$ is a polynomial in $x$.
\end{Def}

Then it is known \cite{Sch} that $\cW$ forms a core for both $A$ and $B$ and it is also stable under Fourier transform. It is obvious that these functions have analytic continuation to the whole complex plane, and they have rapid decay along the real direction. All the operators in the remaining sections will first be defined on $\cW$ and extended by continuity to all of $L^2(\R)$ or its natural domain if the operator is unbounded.

We will also use the notion 
\Eq{\label{WxW}\cW\hat{\ox}\cW\sub L^2(\R\x \R, dsdt),}
where the extra exponent $e^{\c st},\c\in\C$ is allowed.

For convenience, we describe some computations involving $A$ and $B$, and similarly for the operators $\hat{A}$ and $\hat{B}$ with the opposite commutation relations.
\begin{Lem} \label{qqcom} For $AB=q^2BA$, $\hat{A}\hat{B}=q^{-2}\hat{B}\hat{A}$, both $q\inv BA\inv$ and $q\hat{B}\hat{A}\inv$ are positive self-adjoint operators. By the Baker-Campbell-Hausdorff formula
\Eq{e^{2\pi bx}e^{2\pi bp} = e^{2\pi b(x+p)}e^{(2\pi b)^2[x,p]/2}=q e^{2\pi b(x+p)},}
we deduce for example
\Eq{q\inv BA\inv = q\inv e^{2\pi bp}e^{-2\pi bx} = e^{2\pi b(p-x)}.}
Hence we can describe
\Eq{(q\inv BA\inv)^{ib\inv \t} = e^{\pi i\t^2}B^{ib\inv \t}A^{-ib\inv \t}=e^{-\pi i\t^2}A^{-ib\inv \t}B^{ib\inv \t},}
\Eq{(q \hat{B}\hat{A}\inv)^{ib\inv \t} =e^{\pi i\t^2}\hat{A}^{-ib\inv \t}\hat{B}^{ib\inv \t}= e^{-\pi i\t^2}\hat{B}^{ib\inv \t}\hat{A}^{-ib\inv \t}.}
Furthermore we have commutation relations of the form
\Eq{\label{comW1} (B\ox 1) e^{\frac{i}{2\pi b^2}\log A\inv\ox \log \hat{A}} = e^{\frac{i}{2\pi b^2}\log A\inv\ox \log \hat{A}} (B\ox \hat{A}\inv),}
\Eq{\label{comW2} (1\ox \hat{B}) e^{\frac{i}{2\pi b^2}\log A\inv\ox \log \hat{A}} = e^{\frac{i}{2\pi b^2}\log A\inv\ox \log \hat{A}} (A\ox \hat{B}),}
and similar variants.
\end{Lem}


\section{The quantum dilogarithm}\label{sec:QDilog}
We recall the definition of the quantum dilogarithm given in \cite{Ip} (see also \cite{BT}), an important special function that will be used throughout the paper.
\subsection{Definition and properties}\label{sec:QDilog:Def}
Throughout this section, we let $q=e^{\pi i b^2}$ where $b^2\in \R\setminus\Q$ and $0<b^2<1$, so that $|q|=1$. We also denote $Q=b+b\inv$.

Let $\w:=(w_1,w_2)\in\C^2$.

\begin{Def}The Double Zeta function is defined as
\Eq{\ze_2(s,z|\w):=\sum_{m_1,m_2\in\Z_{\geq0}}(z+m_1w_1+m_2w_2)^{-s}.}

The Double Gamma function is defined as
\Eq{\G_2(z|\w):=\exp\left(\frac{\partial}{\partial
s}\ze_2(s,z|\w)|_{s=0}\right).}

Let
\Eq{\G_b(x):=\G_2(x|b,b\inv).}

The Quantum Dilogarithm is defined as the function:
\Eq{S_b(x):=\frac{\G_b(x)}{\G_b(Q-x)}.}

The following form is often useful, and will be used throughout
this paper: \Eq{G_b(x):=e^{\frac{\pi i}{2}x(x-Q)}S_b(x).}
\end{Def}

Let us also relate $G_b$ to another well-known expression by Ruijsenaars \cite{Ru}
\begin{Lem}\label{GbG}
\Eq{G_b(z)=G(b,b\inv;iz-\frac{iQ}{2})e^{\frac{\pi i}{2}z(z-Q)},}
where
\Eq{G(a_+,a_-;z):=\exp\left(i\int_0^\oo \frac{dy}{y}\left(\frac{\sin 2yz}{2\sinh(a_+y)\sinh(a_-y)}-\frac{z}{a_+a_-y}\right)\right)}
with $|Im(z)|<(a_++a_-)/2$ and extends meromorphically to the whole complex plane.
\end{Lem}

The quantum dilogarithm satisfies the following properties:
\begin{Prop} Self-Duality:
\Eq{S_b(x)=S_{b\inv}(x),\tab G_b(x)=G_{b\inv}(x);}

Functional equations: \Eq{\label{funceq}S_b(x+b^{\pm1})=2\sin(\pi
b^{\pm1}x)S_b(x),\tab G_b(x+b^{\pm 1})=(1-e^{2\pi ib^{\pm 1}x})G_b(x);}

Reflection property:
\Eq{\label{reflection}S_b(x)S_b(Q-x)=1,\tab G_b(x)G_b(Q-x)=e^{\pi ix(x-Q)};}

Complex Conjugation: \Eq{\overline{S_b(x)}=\frac{1}{S_b(Q-\over[x])},\tab \overline{G_b(x)}=\frac{1}{G_b(Q-\bar{x})},}
in particular \Eq{\label{gb1}\left|S_b(\frac{Q}{2}+ix)\right|=\left|G_b(\frac{Q}{2}+ix)\right|=1 \mbox{ for $x\in\R$};}

Analyticity:

$S_b(x)$ and $G_b(x)$ are meromorphic functions with poles at
$x=-nb-mb\inv$ and zeros at $x=Q+nb+mb\inv$, for $n,m\in\Z_{\geq0}$;

\label{asymp} Asymptotic Properties:
\Eq{G_b(x)\sim\left\{\begin{array}{cc}\bar{\ze_b}&Im(x)\to+\oo\\\ze_b
e^{\pi ix(x-Q)}&Im(x)\to-\oo\end{array},\right.}
where
\Eq{\ze_b=e^{\frac{\pi i}{4}+\frac{\pi i}{12}(b^2+b^{-2})};}

\label{residue} Residues:
\Eq{\lim_{x\to 0} xG_b(x)=\frac{1}{2\pi},} or more generally,
\Eq{Res\frac{1}{G_b(Q+z)}=-\frac{1}{2\pi}\prod_{k=1}^n(1-q^{2k})\inv\prod_{l=1}^m(1-\widetilde{q}^{2l})\inv}
at $z=nb+mb\inv, n,m\in\Z_{\geq0}$ and $\widetilde{q}=e^{\pi i
b^{-2}}$.
\end{Prop}

From the asymptotic properties, we have the following useful corollary that is needed when we deal with interchanging of order of integrations:
\begin{Cor} \label{asymp} For $s,t\in\C, x\in\R$, the asymptotic behavior of the ratio is given by
\Eq{\left|\frac{G_b(s+ix)}{G_b(t+ix)}\right|\sim\left\{\begin{array}{cc}1&x\to+\oo\\e^{2\pi x Re(t-s)}&x\to-\oo\end{array}.\right.}
\end{Cor}

By analytic continuation in $b$, there exists a classical limit for $G_b(x)$ given in \cite{Ip}:
\begin{Thm}\label{classicallimit} For $q=e^{\pi ib^2}$, by letting $b^2\to i0^+$ we have
\Eq{\lim_{b\to 0} \frac{G_b(bx)}{\sqrt{-i}|b|(1-q^2)^{x-1}}=\G(x),}
where $\sqrt{-i}=e^{-\frac{\pi i}{4}}$ and $-\frac{\pi}{2}<\arg(1-q^2)<\frac{\pi}{2}$. The limit converges uniformly for any compact set in $\C$.
\end{Thm}

We will also need another important variant of the quantum dilogarithm function:
\Eq{g_b(x):=\frac{\bar{\ze_b}}{G_b(\frac{Q}{2}+\frac{1}{2\pi i b}\log x)},}
where $|g_b(x)|=1$ when $x\in\R_{>0}$ due to \eqref{gb1}.
\begin{Lem} Let $u, v$ be self-adjoint operators with $uv=q^2vu$, $q=e^{\pi i b^2}$. Then
\Eq{\label{qexp}g_b(u)g_b(v)=g_b(u+v),}
\Eq{\label{qpenta}g_b(v)g_b(u)=g_b(u)g_b(q\inv uv)g_b(v).}
\end{Lem}

\eqref{qexp} and \eqref{qpenta} are often referred to as the quantum exponential and the quantum pentagon relations, which follows from the other properties:
\Eq{\label{qsum1}g_b(u)^*vg_b(u) = q\inv uv+v,}
\Eq{\label{qsum2}g_b(v)ug_b(v)^*=u+q\inv uv.}

\subsection{Integral transformations}\label{sec:QDilog:Int}
Here we describe several properties of the quantum dilogarithm involving integrations.

\begin{Lem}\label{FT} \cite[(3.31), (3.32)]{BT} We have the following Fourier transformation formula:
\Eq{\int_{\R+i0} e^{2\pi i t r}\frac{e^{-\pi i t^2}}{G_b(Q+i
t)}dt=\frac{\bar{\ze_b}}{G_b(\frac{Q}{2}-ir)}=g_b(e^{2\pi b r}), }
\Eq{\int_{\R+i0} e^{2\pi i t r}\frac{e^{-\pi Qt}}{G_b(Q+i
t)}dt=\ze_b G_b(\frac{Q}{2}-ir)=\frac{1}{g_b(e^{2\pi br})}=g_b^*(e^{2\pi br}) ,} where the contour goes above the pole
at $t=0$.

Using the reflection properties, we also obtain \Eq{\int_{\R+i0} 
e^{2\pi i t r}e^{\pi Q t}G_b(-i
t)dt=\frac{\bar{\ze_b}}{G_b(\frac{Q}{2}-ir)},}
\Eq{\int_{\R+i0} 
e^{2\pi i t r}e^{\pi it^2}G_b(-it)dt=\ze_b G_b(\frac{Q}{2}-ir),}
where again the contour goes above the pole at $t=0$.
\end{Lem}

\begin{Lem} \label{qbi}\cite[B.4]{BT} $q$-Binomial Theorem:
For positive self-adjoint variables $u,v$ with $uv=q^2vu$, we have:
\Eq{(u+v)^{ib\inv t}=\int_{C}\veca{it\\i\t}_b u^{ib\inv (t-\t)}v^{ib\inv \t}d\t ,}
where the $q$-beta function (or $q$-binomial coefficient) is given by
\Eq{\veca{t\\\t}_b=\frac{G_b(-\t)G_b(\t-t)}{G_b(-t)},}
and $C$ is the contour along $\R$ that goes above the pole at $\t=0$
and below the pole at $\t=t$.
\end{Lem}

\begin{Lem}\label{tau} \cite[Lem 15]{PT2} We have the Tau-Beta Theorem:
\Eq{\int_C e^{-2\pi \t \b}
\frac{G_b(\a+i\t)}{G_b(Q+i\t)}d\t =\frac{G_b(\a)G_b(\b)}{G_b(\a+\b)},}
where the contour $C$ goes along $\R$ and goes above the poles of
$G_b(Q+i\t)$ and below those of $G_b(\a+i\t)$. By the asymptotic properties of $G_b$, the integral converges for $Re(\b)>0, Re(\a+\b)<Q$.
\end{Lem}

\begin{Lem}\label{45}Rewriting the integral transform in \cite{Vo} in terms of $G_b$, we obtain the 4-5 relation given by:
\Eq{\int_C d\t e^{-2\pi \c\t}\frac{G_b(\a+i\t)G_b(\b+i\t)}{G_b(\a+\b+\c+i\t)G_b(Q+i\t)}=\frac{G_b(\a)G_b(\b)G_b(\c)}{G_b(\a+\c)G_b(\b+\c)},}
where the contour $C$ goes along $\R$ and goes above the poles of the denominator, and below the poles of the numerator. By the asymptotic properties of $G_b$, the integral converges for $Re(\c)>0$.
\end{Lem}

An important corollary of Lemma \ref{tau} is the following:
\begin{Cor} \label{delta} Let $f(z)\in \cW$ (cf. Definition \ref{WDef}), then we have
\Eq{\lim_{\e\to 0}\int_C\frac{G_b(Q-\e)G_b(iz)}{G_b(Q+iz-\e)} f(z)dz = f(0),}
where the contour go above the poles of the denominator and below the poles of the numerator. Alternatively, we can shift the poles so that the expression is an integration over $\R$:
\Eq{\lim_{\e\to 0}\int_\R\frac{G_b(Q-2\e)G_b(iz+\e)}{G_b(Q+iz-\e)} f(z)dz = f(0).}
Formally as distribution,
\Eq{\frac{G_b(Q)G_b(ix)}{G_b(Q+ix)}=\d(x),\tab x\in\R.}
\end{Cor}

\begin{proof} By Lemma \ref{tau}, we have
\begin{eqnarray*}
&&\lim_{\e\to 0}\int_\R\frac{G_b(Q-2\e)G_b(iz+\e)}{G_b(Q+iz-\e)} f(z)dz\\
&=&\lim_{\e\to 0}\int_\R \int_C e^{-2\pi\t(iz+\e)}\frac{G_b(Q-2\e+i\t)}{G_b(Q+i\t)}f(z)d\t dz,\\
\end{eqnarray*}
where the contour $C$ goes above the pole at $\t=0$.

Now by Corollary \ref{asymp}, the integrand in $\t$ has asymptotics
$$\left|e^{-2\pi\t\e}\frac{G_b(Q-2\e+i\t)}{G_b(Q+i\t)}\right|\sim\left\{\begin{array}{cc}e^{-2\pi \t\e}& \t\to\oo\\ e^{2\pi \t\e}&\t\to-\oo\end{array}.\right.$$
Hence the integrand is absolutely convergent in both $\t$ and $z$. Hence we can interchange the integration order of $\t$ and $z$,
\begin{eqnarray*}
&=&\lim_{\e\to 0}\int_C \int_\R e^{-2\pi\t(iz+\e)}\frac{G_b(Q-2\e+i\t)}{G_b(Q+i\t)}f(z)dz d\t \\
&=&\lim_{\e\to 0}\int_C e^{-2\pi\t\e}\frac{G_b(Q-2\e+i\t)}{G_b(Q+i\t)}(\cF f)(\t)d\t ,\\
\end{eqnarray*}
where $\cF$ is the Fourier transform. Since $f\in \cW, \cF f\in \cW$ as well, hence the integrand is absolutely convergent independent of $\e$. Hence we can interchange the limit and finally obtain
\begin{eqnarray*}
&=&\int_\R \lim_{\e\to 0}e^{-2\pi\t(iz+\e)}\frac{G_b(Q-2\e+i\t)}{G_b(Q+i\t)}(\cF f)(\t)d\t  \\
&=&\int_\R(\cF f)(\t)d\t \\
&=&f(0),
\end{eqnarray*}
where the last line follows from the properties of Fourier transform.
\end{proof}
\begin{Rem}
This is just the analogue of the delta function as a hyperfunction. Since the integrand is $0$ unless $z$ is close to 0, it suffices to notice that near $z=0$, from the analytic properties for $G_b$  we have $G_b(z) \sim \frac{1}{G_b(Q-z)}\sim\frac{1}{2\pi z}$, so that
\begin{eqnarray*}
&&\lim_{\e\to0}\frac{G_b(Q-2\e)G_b(iz+\e)}{G_b(Q+iz-\e)}\\
&\sim&\lim_{\e\to0} \frac{2\e}{2\pi i(iz-\e)(iz+\e)}\\
&=&\lim_{\e\to0}\frac{1}{2\pi i}\left(\frac{1}{z-i\e}-\frac{1}{z+i\e}\right)\\
&=&\d(z).
\end{eqnarray*}

Intuitively, restricting to $x\in\R$, $\frac{G_b(Q)G_b(ix)}{G_b(Q+ix)}$ is always 0 unless $x=0$, in which case it is $\oo$. The properties of $G_b$ say that it gives the right normalization to be a delta function.
\end{Rem}

The other variants which we will also use include:
\begin{Cor} For $f(x)\in\cW$, we have
\begin{eqnarray}
\lim_{\e\to0}\int_\R \frac{G_b(\e-ix)G_b(\e+ix)}{G_b(2\e)}f(x)dx &=& f(0),\\
\lim_{\e\to0}\int_\R \frac{G_b(\e+ix-b)G_b(Q+b-2\e)}{G_b(Q+ix-\e)}f(x)dx &=& -q^2 f(0)+f(-ib),
\end{eqnarray}
or formally
\begin{eqnarray}
\frac{G_b(-ix)G_b(ix)}{G_b(0)}&=&\d(x),\\
\frac{G_b(ix-b)G_b(Q+b)}{G_b(Q+ix)}&=&-q^2\d(x)+\d(x+ib),
\end{eqnarray}
which follows from the functional properties of $G_b(x)$ and similar arguments as above.
\end{Cor}

A particular important case in the study of modular double (see the next section) is proved in \cite[Lem 3]{BT} and \cite{Vo} which can also be obtained from the above arguments:
\begin{Cor}\label{u+v} Let $u,v$ be positive self-adjoint operators satisfying $uv=q^2vu$. By taking $t\to -ib\inv$ in Lemma \ref{qbi} we obtain
\Eq{(u+v)^{1/b^2} = u^{1/b^2}+v^{1/b^2}.}
\end{Cor}

\section{The quantum plane}\label{sec:QPlane}
In this section, we will define the main object of study in this paper, the $C^*$-algebra $\cA$ of ``functions vanishing at infinity" on the quantum plane. Its Haar functional is established, and we obtain the GNS construction of $\cA$ on a Hilbert space $\cH\simeq L^2(\R\x\R)$. Finally we construct the multiplicative unitary $W$, and show that it is in fact manageable.

Note that in analogy to the classical semi-group with inverses absent, in our context we don't have a well-defined antipode $S$ due to positivity. Hence we may also call the object informally as a ``quantum semigroup". However, as we see in Definition \ref{unitaryantipode}, there is instead a well-defined unitary antipode $R$.

With the existence of the coproduct (Corollary \ref{coproductM}), a left invariant faithful KMS Haar functional (Theorem \ref{Haar}, Theorem \ref{GNSmap}, Corollary \ref{KMS}), a right invariant Haar functional (Proposition \ref{right}) and the density condition (Theorem \ref{density}), the $C^*$-algebra $\cA$ can be put in the context of ``locally compact quantum group" in the sense of \cite{KV1,KV2}. This allows us to apply the theory developed in those papers.

\subsection{Motivations from classical $ax+b$ group}\label{sec:QPlane:ax+b}
Recall that the classical $ax+b$ group is the group $G$ of affine transformations on the real line $\R$, where $a>0$ and $b\in\R$, and they can be represented by a matrix of the form (we use the transposed version):
\Eq{g(a,b)=\veca{a&0\\b&1},}
where the multiplication is given by
\Eq{\label{ax+bprod} g(a_1,b_1)g(a_2,b_2)=\veca{a_1a_2&0\\b_1a_2+b_2&1}.}
We can then talk about the space $C_\oo(G)$ of functions on $G$ vanishing at infinity, and a dense subspace of functions $\cW_G$ of the form
\Eq{\cW_G=\{g(\log a)f(b): f,g\in \cW\sub L^2(\R)\},}
where $f,g\in \cW$ are the rapidly decaying functions defined in Definition \ref{WDef}. Then $C_\oo(G)$ is the sup-norm closure of $\cW_G$.

In the study of quantum plane, it is necessary to consider the semigroup $G_+$ of positive elements. Restricting $f(a,b)\in\cW_G$ to $b>0$, we can use the Mellin transform and express the functions as
\Eq{f(a,b)=\int_{\R+i0}\int_\R  F(s,t)a^{is}b^{it}dsdt.}
Since the function is continuous at $b=0$, it has at most $O(1)$ growth as $b\to 0^+$. Hence from Theorem \ref{MTpole}, we conclude that $F(s,t)$ is entire analytic with respect to $s$, and holomorphic on Im$(t)>0$. Furthermore, it has rapid decay in $s,t$ along the real direction, and can be analytically continued to Im$(t)\leq 0$. Moreover, since $f(a,b)$ is analytic at $b=0$, the analytic structure of $f(a,b)$ on $b$ is given by $\sum_{k=0}^\oo A_k b^k$ for some constants $A_k$. Hence $F(s,t)$ has possible simple poles at $t=-in, n=0,1,2,...$.

Finally according to the Parseval's Formula for Mellin transform, there is an $L^2$ norm on functions of $G_+$ given by
\Eq{\|f(a,b)\|^2 = \int_\R \int_\R|F(s,t+\frac{1}{2}i)|^2 dt ds}
induced from the Haar measure $\frac{da}{a} db$. However, it is no longer left invariant due to the lack of inverses.

\subsection{The quantum plane algebra and its modular double}\label{sec:QPlane:alg}
Throughout the section, we let $q=e^{\pi ib^2}, \til[q]=e^{\pi ib^{-2}}$ where $b^2\in\R\setminus \Q, 0<b^2<1$. We have $|q|=|\til[q]|=1$.

The quantum plane $\cA_q$ is formally the Hopf algebra generated by positive self-adjoint operators $A, B$ satisfying
\Eq{AB=q^2BA,}
with coproduct
\begin{eqnarray}
\D(A)&=&A\ox A,\\
\label{copB}\D(B)&=&B\ox A+1\ox B.
\end{eqnarray}
Hence it is clear from the previous section that the quantum plane is just the quantum analogue of the classical $ax+b$ group. Their relationships are studied in detail in \cite{Ip}.
By the following proposition, the coproducts are positive essentially self-adjoint, hence they are well-defined.

\begin{Prop}\label{x+p} The operator $e^{2\pi bx}+e^{2\pi bp}$ is positive and essentially self-adjoint.
In particular the coproduct $\D(B)=B\ox A+1\ox B$ defined in \eqref{copB} is positive and essentially self-adjoint operator on $L^2(\R)\ox L^2(\R)$.
\end{Prop}
\begin{proof} 
It follows from \eqref{qsum2} by applying $u=A, v=q\inv BA\inv$, and notice that $g_b(v)$ is unitary. Hence this gives a unitary transformation that sends $e^{2\pi bx}$ to $e^{2\pi bx}+e^{2\pi bp}$. In particular the functional analytic property carries over to the new operator.
\end{proof}
Since we are working with positive operators, this important property allows us to avoid the analysis of the coproduct of $B$ that is studied extensively in \cite{PuW} and \cite{WZ}.

An interesting object in the study of quantum plane is the modular double element \cite{Fa,FK,KLS}. $\cA_{\til[q]}$ is formally the Hopf algebra generated by positive self-adjoint elements
\Eq{\til[A]=A^{1/b^2},\tab \til[B]=B^{1/b^2},}
such that they satisfy
\Eq{\til[A]\til[B]=\til[q]^2\til[B]\til[A].}

Then according to the formula above, using Corollary \ref{u+v} we conclude that they have the same bi-algebra structure:
\begin{eqnarray}
\D(\til[A])&=&\til[A]\ox\til[A],\\
\D(\til[B])&=&(B\ox A+1\ox B)^{1/b^2}\nonumber\\
&=&(B\ox A)^{1/b^2}+(1\ox B)^{1/b^2}\nonumber\\
&=&\til[B]\ox \til[A]+1\ox \til[B].
\end{eqnarray}

Motivated from the definition given in the previous section for the classical $ax+b$ group, we define the $C^*$-algebra $C_\oo(\cA_q)$ of ``functions on the quantum plane vanishing at infinity" as follows.

\begin{Def}\label{qplanedef} $C_\oo(\cA_q)$ is the space of all functions
\Eq{C_\oo(\cA_q):=\left\{\int_{\R+i0}\int_\R f(s,t)A^{ib\inv s}B^{ib\inv t}dsdt\right\}^{\mbox{norm closure}}} 
where $f$ is entire analytic in $s$, and meromorphic in $t$ with possible simple poles at
$$t=-ibn-i\frac{m}{b},\tab n,m=0,1,2,...$$
and for fixed $v>0$, the functions $f(s+iv,t)$ and $f(s,t+iv)$ are of rapid decay (faster than any exponential) in both $s$ and $t$.
To define the $C^*$ norm, we realize $A^{ib\inv s}f(x)=e^{2\pi is}f(x)$ and $B^{ib\inv t}f(x)=e^{2\pi ip}f(x)=f(x+1)$ as unitary operators on $L^2(\R)$ and take the operator norm.

For simplicity, we will denote the space simply by $\cA:= C_\oo(\cA_q)$.
\end{Def}

\begin{Rem}\label{MellinQ}
As in the classical case, by Mellin transform in $s$, we can consider the linear span of $$g(\log A)\int_{\R+i0}f_2(t)B^{ib\inv t} dt$$ instead, where $g(x)$ is entire analytic and of rapid decay in $x\in\R$. This form will be useful later when we deal with the pairing between $\cA$ and its dual space $\hat{\cA}$. We can also conclude from the definition that $A,A\inv, \til[A], \til[A]\inv, B$ and $\til[B]$ are all elements of the multiplier algebra $M(\cA)$.
\end{Rem}

This definition differs from the classical case in certain aspects. First of all, there is a $b\longleftrightarrow b\inv$ duality in this definition. Indeed, it is obvious that we have
\Eq{C_\oo (\cA_q,b) = C_\oo (\cA_{\til[q]},b\inv)}
because $A^{ib\inv s}B^{ib\inv t} = \til[A]^{ibs}\til[B]^{ibt}$ and the analyticity of $f(s,t)$ has the $b\corr b\inv$ duality.
Hence it encodes the information of both $\cA_q$ and $\cA_{\til[q]}$. For example, we see by Theorem \ref{MTpole} that as a function of $B$, at $B=0$ it admits a series representation
\Eq{F(B)\sim \sum_{m,n\geq 0}\a_{mn} B^{m+n/b^2} =\sum_{m,n\geq 0}a_{mn}B^{m}\til[B]^n.}
This choice of poles in the definition is needed in order for the coproduct $\D(\cA)$ to lie in $M(\cA\ox \cA)$, the multiplier algebra, because of the appearances of the quantum dilogarithm $G_b(-i\t)$, hence it is a ``minimal" choice for all the calculations to work.

\begin{Prop}\label{coprod} Let us denote simply by $f:=\int_{\R+i0}\int_\R f(s,t)A^{ib\inv s}B^{ib\inv t} dsdt$. Then the coproduct on $\cA_q$ can be naturally extended to $\cA$ by
\begin{eqnarray*}
\D(f)&=&\int_{\R+i0}\int_\R f(s,t) (A^{ib\inv s}\ox A^{ib\inv s}) (B\ox A+1\ox B)^{ib\inv t}dsdt\\
&=&\int_{\R+i0}\int_\R\int_C f(s,t) \frac{G_b(-i\t)G_b(i\t-it)}{G_b(-it)}A^{ib\inv s}B^{ib\inv (t-\t)}\ox A^{ib\inv (s+t-\t)}B^{ib\inv \t}d\t dsdt\\
&=&\int_{\R+i0}\int_\R\int_{\R+i0} f(s,t+\t) \frac{G_b(-i\t)G_b(-it)}{G_b(-i(t+\t))}A^{ib\inv s}B^{ib\inv t}\ox A^{ib\inv (s+t)}B^{ib\inv \t}d\t dsdt,
\end{eqnarray*}
or formally
\Eq{\label{coproductformal} \D(f)(s_1,t_1,s_2,t_2)=f(s_1,t_1+t_2)\frac{G_b(-it_1)G_b(-it_2)}{G_b(-i(t_1+t_2))}\d(s_2,s_1+t_1).}
\end{Prop}

\begin{Cor}\label{coproductM} The map $\D$ sends $\cA\to M(\cA\ox \cA)$ where $M(\cA)$ is the multiplier algebra for $\cA$, and it is coassociative in the sense of multiplier Hopf algebra. Hence $\D$ is indeed a coproduct on $\cA$.
\end{Cor}
\begin{proof}
We need to show that $\D(f)(g\ox h)$ lies in $\cA\ox \cA$ for $f,g,h\in\cA$. From the formula for $\D$, we see that the poles for $f$ are canceled by $G_b(-i(t+\t))$, and two new set of poles at the specified locations are introduced by the numerator $G_b(-i\t)G_b(-it)$. Hence the integrand of the coproduct is a generalized function with the specified analytic properties. Since $C_\oo(\cA_q)$ is closed under product, which comes from the Mellin transform of certain series representation, the product of $\D(f)$ and any element $g\ox h \in \cA\ox \cA$ will lie in $\cA\ox\cA$.

Coassociativity is then obvious from the construction.
\end{proof}

Finally let us describe the antipode $S$ and the unitary antipode $R$. On the Hopf algebra level, the antipode is given by $S(A)=A\inv, S(B)=-BA\inv$. Formally we extend it to $\cA$ by
$$S\left(\int_{\R+i0}\int_\R f(s,t)A^{ib\inv s}B^{ib\inv t}dsdt\right)=\left(\int_{\R+i0}\int_\R f(s,t)(e^{\pi i}BA\inv)^{ib\inv t}A^{-ib\inv s}dsdt\right),$$
and hence we have
\begin{Def}\label{unitaryantipode} We define the antipode to be
\Eq{S(f)=f(-s-t,t)e^{-\pi Qt}e^{\pi i(2st+t^2)}.}
Then $S^2(f)(s,t)=f(s,t)e^{-2\pi Qt}$, hence we conclude the scaling group to be
\Eq{\t_{-i}(A^{ib\inv s})=A^{ib\inv s},\tab \t_{-i}(B^{ib\inv t})=e^{-2\pi Qt}B^{ib\inv t},}
or more generally 
\Eq{\t_t(A)=A,\tab\t_t(B)=e^{-2\pi bQt}B.}
Then the unitary antipode $R=\t_{i/2}S$ is given by
\begin{eqnarray}
R(A)&=&A\inv,\\
R(B)&=&-e^{-\pi ibQ}BA\inv = e^{-\pi ib^2-\pi i}BA\inv = q\inv BA\inv.
\end{eqnarray}
\end{Def}
We note that $R$ sends positive operators to positive operators, which differs from the usual choice (see e.g. \cite[Thm 8.4.33]{T}). Furthermore, when we later realize $f(s,t)$ as element in $\cH:=L^2(\R^2)$ in the von Neumann picture, the antipode $S$ can be seen to be an unbounded operator on $\cH$ due to the factor $e^{-\pi Qt}$ coming from the negative sign, while $R$ is unitary. Hence $R$ is better suit in this positive semigroup context.

It is known that as polynomial algebras $\cA_q$ and $\cA_{\til[q]}$ are $C^*$-algebras of Type $II_1$ \cite{Fa,KLS}. However, by expressing $A=e^{2\pi bx}, B=e^{2\pi bp}$, the definition of $\cA$ states that it is nothing but the algebra of all bounded operators on $L^2(\R)$ via the Weyl formula \cite{Fa}, hence it is a $C^*$-algebra of Type I. Therefore $\cA$ is precisely the modular double $\cA_q\ox\cA_{\til[q]}$, and it arises as a natural framework to the study of this subject.

Finally we will also see in the GNS construction that a natural $L^2$-norm is given by (cf. Theorem \ref{GNSmap})
\Eq{\|f(s,t)\|^2 = \int_\R \int_\R |f(s,t+\frac{iQ}{2})|^2 dsdt,}
which is an analogue of the classical Parseval's formula.

\subsection{The Haar functional}\label{sec:QPlane:Haar}
Recall from Definition \ref{leftinv} that a linear functional on $\cA$ is called a left invariant Haar functional if it satisfies
\Eq{(1\ox h)(\D x) = h(x)\cdot 1_{M(\cA)},}
and a right invariant Haar functional if it satisfies
\Eq{(h\ox 1)(\D x)=h(x)\cdot 1_{M(\cA)},}
where $1_{M(\cA)}$ is the unital element in the multiplier algebra $M(\cA)$.

\begin{Thm}\label{Haar} There exists a left invariant Haar functional on $\cA$, given by
\Eq{h\left(\int_{\R+i0}\int_\R f(s,t)A^{ib\inv s}B^{ib\inv t}dsdt \right) = f(0,iQ).}
\end{Thm}
\begin{proof} From Proposition \ref{coprod}, the coproduct is given by
\begin{eqnarray*}
&&\D\left(\int_{\R+i0}\int_\R f(s,t)A^{ib\inv s}B^{ib\inv t}dsdt \right) \\
&=&\int_{\R+i0}\int_\R\int _C f(s,t)\frac{G_b(-i\t)G_b(i\t-it)}{G_b(-it)}A^{ib\inv s}B^{ib\inv (t-\t)}\ox A^{ib\inv (s+t-\t)}B^{ib\inv \t} d\t dsdt \\
&=&\int_\R\int_{\R-i0}\int _C f(s,t+\t-s)\frac{G_b(-i\t)G_b(is-it)}{G_b(is-it-i\t)}A^{ib\inv s}B^{ib\inv (t-s)}\ox A^{ib\inv t}B^{ib\inv \t} d\t dsdt. \\
\end{eqnarray*}
In the last step, we push the contour of $t$ back to the real line, and the contour of $\t$ to be the real line that goes above the pole at $\t=0$ only. This will force the contour of $s$ to move below the pole at $s=t$.

Now as a function of $\t$ the integrand is well-defined for $0<Im(\t)<Q$. Therefore to apply the Haar functional, we take the analytic continuation to $\t=iQ$ by $$\dis\t=\lim_{\e\to 0}iQ-i\e.$$ Using Corollary \ref{delta}, we have
\begin{eqnarray*}
&=&\lim_{\e\to 0}\int_{\R-i0} f(s,iQ-s-i\e)\frac{G_b(Q-\e)G_b(is)}{G_b(Q+is-\e)}A^{ib\inv s} B^{ib\inv (-s)}ds\\
&=&f(0,iQ)
\end{eqnarray*}
as desired.
\end{proof}

Using exactly the same technique, we can derive the counit for $\cA$:
\begin{Cor} The counit is given by
\Eq{\e(f) = \lim_{t\to 0} G_b(Q+it)\int_\R f(s,t)ds}
and satisfies
\Eq{(1\ox \e)\D = Id_\cA = (\e\ox 1)\D.}
\end{Cor}
\subsection{The GNS description}\label{sec:QPlane:GNS}
The GNS representation enables us to bring the quantum plane into the Hilbert space level, where the algebra is realized as operators on certain Hilbert space. This Hilbert space is nothing but the space of ``$L^2$ functions on the quantum plane", which gives the necessary background later to describe the $L^2$ functions on $GL_q^+(2,\R)$ using the quantum double construction. This will allow us to state the main theorem (Theorem \ref{main}) on the decomposition of $L^2(GL_q^+(2,\R))$ into the fundamental principal series representations of $U_{q\til[q]}(\gl(2,\R))$ with certain Plancherel measure.

Now equipped with the left invariant Haar functional $h$, we can describe completely the GNS representation of $\cA$ using $h$ as the weight.

\begin{Thm} \label{GNSmap}
Let $\cH=L^2(\R\x \R)$ be the completion of $L^2(\R) \ox \cW$ equipped with the inner product
\Eq{\<f,g\> = \int_\R\int_\R \overline{g(s,t+\frac{iQ}{2})}f(s,t+\frac{iQ}{2})dsdt,}
and the $L^2$ norm
\Eq{\|f\|^2 = \iint |f(s,t+\frac{iQ}{2})|^2dsdt.}
Then the GNS map is simply given by
\Eq{\L\left(\int_{\R+i0}\int_\R f(s,t)A^{ib\inv s}B^{ib\inv t}dsdt \right) = f(s,t).}
In particular, the Haar functional is a faithful weight.

The action $\pi$ is given by the multiplication of $\cA$ on itself:
\Eq{\pi(a)\L(b):=\L(ab),\tab a,b\in\cA.}
\end{Thm}

First of all, let us note that $\cA$ is closed under *:
\begin{Prop} The conjugation map $T$ is given by
\Eq{T(f)=\over[f](-s,-t)e^{2\pi ist},}
where we denote by $\over[f](z):=\overline{f(\over[z])}$.

We observe that $\over[f](-s,-t)e^{2\pi ist}$ still satisfies all the analytic properties required to be an element of $\cA$: it is entire analytic in $s$, meromorphic in $t$ with possible poles at $t=-inb-imb\inv$, etc. see Definition \ref{qplanedef}.
\end{Prop}
\begin{proof}
\begin{eqnarray*}
&&\left(\int_{\R+i0}\int_\R f(s,t)A^{ib\inv s}B^{ib\inv t}dsdt \right)^* \\
&=&\int_{\R+i0}\int_\R \over[f(s,t)] B^{-ib\inv \over[t]}A^{-ib\inv \over[s]}d\over[s]d\over[t]\\
&=&\int_{\R-i0}\int_\R \overline{f(\over[s],\over[t])} B^{-ib\inv t}A^{-ib\inv s}dsdt\\
&=&\int_{\R-i0}\int_\R  \overline{f(\over[s],\over[t])}e^{2\pi i s t}A^{-ib\inv s}B^{-ib\inv t}dsdt\\
&=&\int_{\R+i0}\int_\R \overline{f(-\over[s],-\over[t])}e^{2\pi i st}A^{ib\inv s}B^{ib\inv t}dsdt\\
&=&\int_{\R+i0}\int_\R  \over[f](-s,-t)e^{2\pi i st}A^{ib\inv s}B^{ib\inv t}dsdt.
\end{eqnarray*}
\end{proof}

Next we describe the product in terms of $f(s,t)$:
\begin{Prop} \label{producttwisted}The product of two elements in $\cA$ is given by:
\Eq{(f\cdot g)(s,t)=\int_C\int_\R f(s-s',t-t')g(s',t')e^{2\pi is'(t-t')} ds' dt'}
i.e. a twisted convolution. Here $C$ is the contour that, with possible poles, goes below $t'=t$ and above $t'=0$.
\end{Prop}
\begin{proof}
\begin{eqnarray*}
&&\left(\int_{\R+i0}\int_\R f(s,t)A^{ib\inv s}B^{ib\inv t}dsdt\right)\left(\int_{\R+i0}\int_\R g(s',t')A^{ib\inv s'}B^{ib\inv t'}ds'dt'\right)\\
&=&\int_{\R+i0}\int_\R\int_{\R+i0}\int_\R f(s,t)g(s',t')e^{2\pi is't}A^{ib\inv (s+s')}B^{ib\inv (t+t')} ds'dt' dsdt\\
&=&\int_{\R+i0}\int_\R\left(\int_C\int_\R f(s-s',t-t')g(s',t')e^{2\pi is'(t-t')}ds'dt'\right)A^{ib\inv s}B^{ib\inv t} dsdt,\\
\end{eqnarray*}
where in the last step, the shift in $t$ will push the contour of $t'$ to go below $t'=t$.
\end{proof}

\begin{proof}[Proof of Theorem \ref{GNSmap}] Using the formula for the Haar functional, and the definition of the GNS inner product, we have
\begin{eqnarray*}
\<f,g\>&=&h(g^*f) \\
&=&h(\over[g](-s,-t)e^{2\pi ist}\cdot f(s,t))\\
&=&h\left(\int_{\R+i0}\int_\R \over[g](s'-s,t'-t)e^{2\pi i(s-s')(t-t')}f(s',t')e^{2\pi is'(t-t')} ds'dt'\right)\\
&=&h\left(\int_{\R+i0}\int_\R \over[g](s'-s,t'-t)f(s',t')e^{2\pi is(t-t')} ds'dt'\right)\\
&=&\int_{\R+i0}\int_\R \over[g](s',t'-iQ)f(s',t') ds'dt'\\
&=&\int_{\R}\int_\R \over[g](s',t'-\frac{iQ}{2})f(s',t'+\frac{iQ}{2}) ds'dt'\\
&=&\int_{\R}\int_\R \over[g(s',t'+\frac{iQ}{2})]f(s',t'+\frac{iQ}{2}) ds'dt',
\end{eqnarray*}
where in the second to last line, we do a shift in the contour of $t$ by $\frac{iQ}{2}$ by holomorphicity.
\end{proof}

Finally we read off the representation $\pi$ of $\cA$ on the functions $f(s,t)\in\cH$:
\begin{Prop} The representation $\pi: \cA\inj \cB(\cH)$ is given by
\Eq{\pi(A)=e^{-2\pi bp_s},\tab \pi(B)=e^{2\pi bs}e^{-2\pi bp_t},}
where $p_s=\frac{1}{2\pi i}\del[,s]$, so that $e^{2\pi bp_s}\cdot f(s)=f(s-ib)$. Similarly for $p_t$.
\end{Prop}
\begin{proof}
$A$ acts on $f(s,t)$ by:
\begin{eqnarray*}
&&A^{ib\inv s'}\iint f(s,t) A^{ib\inv s}B^{ib\inv t}dsdt\\
&=&\iint f(s,t) A^{ib\inv (s+s')}B^{ib\inv t}dsdt\\
&=&\iint f(s-s',t) A^{ib\inv s}B^{ib\inv t}dsdt.\\
\end{eqnarray*}
Hence $A^{ib\inv s'}\cdot f(s,t) = f(s-s', t)$, or $A= e^{-2\pi bp_s} .$

Similarly, 
\begin{eqnarray*}
&&B^{ib\inv t'}\iint f(s,t) A^{ib\inv s}B^{ib\inv t}dsdt\\
&=&\iint f(s,t) e^{2\pi ist'}A^{ib\inv s}B^{ib\inv (t+t')}dsdt\\
&=&\iint f(s,t-t') e^{2\pi ist'}A^{ib\inv s}B^{ib\inv t}dsdt.\\
\end{eqnarray*}
Hence $B^{ib\inv t'}\cdot f(s,t) = e^{2\pi i st'}f(s,t-t')$, or $B=e^{2\pi bs}e^{-2\pi b p_t}$.
\end{proof}

We will take this Hilbert space $\cH$ as our canonical choice, when we introduce right measure and the dual space. 

\subsection{Modular maps and the right picture}\label{sec:QPlane:modular}
In this section we describe the remaining objects defined in Section \ref{sec:pre:GNS}, as well as the description of the right invariant picture for the quantum plane.

From now on, let us restrict all operators on the dense subspace $\cW$, and we will extend any bounded operators on $\cW$ to the whole space $\cH$.

\begin{Prop} The conjugation map $T$ (cf. \eqref{T}) is given by
\Eq{\label{Tconvolution}T(f)=\over[f](-s,-t)e^{2\pi ist},}
$T^*$ is given by
\Eq{T^*(f)=\over[f](-s,-t)e^{2\pi ist}e^{2\pi sQ}.}
Hence the modular operator $\DD=T^*T$ is given by
\Eq{\DD(f)=f(s,t)e^{2\pi sQ},}
and the modular conjugation $J=T\DD^{-\frac{1}{2}}=\DD^{\frac{1}{2}}T$ is given by
\Eq{J(f)=\over[f](-s,-t)e^{2\pi ist}e^{\pi sQ}.}
We have $J^2=J^*J= Id$.
\end{Prop}
\begin{proof} It suffices to establish the map $T^*$. Since $T$ is anti-linear, from the definition of the adjoint of anti-linear maps, we have:
\begin{eqnarray*}
\<Tf, g\>&=&\over[\<f,T^* g\>]\\
RHS &=& \iint T^*g(\over[s],\over[t]+\frac{iQ}{2})\overline{f(s,t+\frac{iQ}{2})}d\over[s]d\over[t]\\
&=& \iint T^*g(s,t+\frac{iQ}{2})\overline{f(\over[s],\over[t]+\frac{iQ}{2})}dsdt\\
LHS &=& \iint \overline{g(\over[s],\over[t]+\frac{iQ}{2})}\over[f](-s,-t-\frac{iQ}{2})e^{2\pi i s(t+\frac{iQ}{2})}dsdt\\
&=& \iint \overline{g}(s,t-\frac{iQ}{2})\overline{f(-\over[s],-\over[t]+\frac{iQ}{2})}e^{2\pi i s(t+\frac{iQ}{2})}dsdt\\
&=&\iint \overline{g}(-s,-t-\frac{iQ}{2})\overline{f(\over[s],\over[t]+\frac{iQ}{2})}e^{2\pi i (-s)(-t+\frac{iQ}{2})}dsdt.\\
\end{eqnarray*}
Hence
$T^*g=\over[g](-s,-t)e^{2\pi i st}e^{2\pi sQ}$.
\end{proof}

\begin{Cor} The modular group $\s_t(x)=\DD^{it}x\DD^{-it}$ is given by
\begin{eqnarray}
\s_t(A)&=&e^{2\pi bQt}A,\\
\s_t(B)&=&B.
\end{eqnarray}
or in other words
\Eq{\label{modulargroup}\s_{t'}(f)=f(s,t)e^{2\pi Qt' s}}
\end{Cor}

With this expression of the modular group, it is now easy to see the following:
\begin{Cor} \label{KMS} The Haar weight is a KMS weight, in the sense of \cite{KV1}. That is
\begin{eqnarray}
h\circ\s_t&=&h\tab\tab\tab\tab\tab\mbox{ for every }t\in\R,\\
h(a^*a)&=&h(\s_{\frac{i}{2}}(a)\s_{\frac{i}{2}}(a)^*) \tab\mbox{ for every }a\in \cA.
\end{eqnarray}
\end{Cor}
\begin{proof}
The first equation is obvious from \eqref{modulargroup}. For the second equation, it follows by combining \eqref{Tconvolution}, \eqref{modulargroup} and Proposition \ref{producttwisted}.
\end{proof}
\begin{Prop} \label{right} The scaling constant $\nu$ is given by
\Eq{ \nu=e^{-2\pi Q^2}>0,}
the modular element $\d$ is given by
\Eq{\d=A^{-\frac{Q}{b}},}
and the right Haar weight is given by
\Eq{\psi(f)=f(-iQ,iQ)e^{-\pi iQ^2}.}
\end{Prop}
\begin{proof}
\begin{eqnarray*}
h(\t_{t'}(f))&=&h(\iint f(s,t)e^{-2\pi bQt' ib\inv t}A^{ib\inv s}B^{ib\inv t}ds dt\\
&=&h(\iint f(s,t)e^{-2\pi iQt't}A^{ib\inv s}B^{ib\inv t}ds dt\\
&=&f(0,iQ)e^{-2\pi iQt'iQ}\\
&=&f(0,iQ)e^{2\pi Q^2t'},
\end{eqnarray*}
hence $\nu=e^{-2\pi Q^2}$.

Obviously \Eq{\s_t(\d)=e^{2\pi b Qt(-\frac{Q}{b})}\d=e^{-2\pi Q^2t}\d=\nu^t\d.}

Finally for the right Haar weight, recall that
\begin{eqnarray*}
R(A)&=&A\inv,\\
R(B)&=& q\inv BA\inv,
\end{eqnarray*}
hence
\begin{eqnarray*}
\psi(f)&=&h(R(f))\\
&=&h(\iint f(s,t)(q\inv BA\inv)^{ib\inv t}A^{-ib\inv s}ds dt)\\
&=&h(\iint f(s,t)q^{-ib\inv t}q^{-ib\inv t(ib\inv t-1)}q^{-2b^{-2}t(s+t)}A^{-ib\inv (s+t)}B^{ib\inv t}ds dt)\\
&=&h(\iint f(s,t)e^{-\pi i(2st+t^2)}A^{-ib\inv (s+t)}B^{ib\inv t}ds dt\\
&=&h(\iint f(-s-t,t)e^{\pi i(2st+t^2)}A^{ib\inv s}B^{ib\inv t}ds dt\\
&=&f(-iQ,iQ)e^{\pi i (iQ)^2}\\
&=&f(-iQ,iQ)e^{-\pi i Q^2}.
\end{eqnarray*}
\end{proof}

We have an isometry $\cH_R\to \cH_L$ between the Hilbert space associated to the right and left invariant Haar weights respectively, given in Proposition \ref{GNSright} by
\Eq{\L_R(a)=\L_L(a\d^{\frac{1}{2}}) = \L_L(aA^{-\frac{Q}{2b}}).}
Explicitly, we can rewrite all the maps in the right picture. We list here for convenience.
\begin{Prop} We have
\Eq{\L_R\left(\int_{\R+i0}\int_\R f(s,t)A^{ib\inv s}B^{ib\inv t}dsdt\right)=f(s-\frac{iQ}{2},t)e^{-\pi Qt},}
with inner product given by
\Eq{\<f,g\>_R = \int_\R\int_\R \over[g(s-\frac{iQ}{2},t+\frac{iQ}{2})]f(s-\frac{iQ}{2},t+\frac{iQ}{2})e^{-2\pi Qt}dsdt.}
\begin{eqnarray}
T_Rf(s,t)&=&\over[f](-s,-t)e^{2\pi ist}e^{-\pi Qt},\\
T_R^*f(s,t)&=&\over[f](-s,-t)e^{2\pi ist+\pi Qt+2\pi Q s-\pi iQ^2},\\
\DD_Rf(s,t)&=&f(s,t)e^{2\pi sQ+2\pi tQ-\pi iQ^2},\\
J_Rf(s,t)&=&\over[f](-s,-t)e^{2\pi ist}e^{\pi sQ}e^{-\frac{\pi iQ^2}{2}}.
\end{eqnarray}
In particular, we recover the relation \cite{KV2}
\Eq{J_R=\nu^{i/4}J.}
\end{Prop}

\subsection{The multiplicative unitary}\label{sec:QPlane:multi}
In this section we will describe the multiplicative unitary $W$ explicitly. Recall (cf. \eqref{WW}) that it is defined as a unitary operator on $\cH\ox \cH$ by
\Eq{W^*(\L(a)\ox \L(b))=(\L\ox\L)(\D(b)(a\ox1)).}

We will use extensively the following variants of Lemma \ref{FT}:
\begin{Cor}\label{gbint} For $AB=q^2BA, \hat{A}\hat{B}=q^{-2}\hat{B}\hat{A}$, using Lemma \ref{qqcom} we have 
\begin{eqnarray}
g_b(B\ox q\hat{B}\hat{A}\inv) &=& \int_{\R+i0} \frac{B^{ib\inv\t}\ox\hat{A}^{-ib\inv \t}\hat{B}^{ib\inv \t}}{G_b(Q+i\t)}d\t,\\
\label{gb2}g_b(q\inv BA\inv\ox \hat{B}) &=& \int_{\R+i0} \frac{B^{ib\inv\t}A^{-ib\inv \t}\ox\hat{B}^{ib\inv \t}}{G_b(Q+i\t)}d\t,\\
\label{gb3}g_b^*(B\ox q\hat{B}\hat{A}\inv) &=& \int_{\R+i0} B^{ib\inv\t}\ox\hat{B}^{ib\inv \t}\hat{A}^{-ib\inv \t}G_b(-i\t)d\t,\\
g_b^*(q\inv BA\inv\ox \hat{B}) &=& \int_{\R+i0} A^{-ib\inv \t}B^{ib\inv\t}\ox\hat{B}^{ib\inv \t}G_b(-i\t)d\t.
\end{eqnarray}
Note that the arguments inside $g_b$ are all essentially self-adjoint, hence the expression is well-defined.
\end{Cor}
\begin{Prop} \label{W}The multiplicative unitary $W\in \cB(\cH\ox \cH)$ is given by
\begin{eqnarray}
W&=&\label{W1}e^{\frac{i}{2\pi b^2}\log A\inv\ox \log \hat{A}}g_b(B \ox q\hat{B}\hat{A}\inv)\\
&=&\label{W2}g_b(q\inv BA\inv \ox \hat{B})e^{\frac{i}{2\pi b^2}\log A\inv\ox \log \hat{A}},\\
W^*&=&e^{\frac{i}{2\pi b^2}\log A\ox \log \hat{A}}g_b^*(q\inv BA\inv \ox \hat{B})\\
&=&g_b^*(B \ox q\hat{B}\hat{A}\inv)e^{\frac{i}{2\pi b^2}\log A\ox \log \hat{A}},
\end{eqnarray}
where 
\Eq{A=e^{-2\pi bp_s},\tab B=e^{2\pi bs}e^{-2\pi bp_t}} as before, and
\Eq{\hat{A}=e^{2\pi bs},\tab \hat{B}=G_b(-it) \circ e^{2\pi bp_t-2\pi bp_s}\circ G_b(-it)\inv.}
\end{Prop}
\begin{Rem} Note that we actually have $W\in M(\cA\ox \cB(\cH))$. Moreover, in the next section, we will see that the hat operators are precisely the dual space elements, hence indeed $W\in M(\cA\ox \hat{\cA})$. Furthermore, the multiplicative unitary will be invariant under the exchange $b\to b\inv$ and $\cA_q\to \cA_{\til[q]}$.
\end{Rem}
\begin{proof} First of all notice that $\hat{A}$ and $\hat{B}$ are positive operators on the subspace $\cW$ (with the measure shifted in $t$ by $\frac{iQ}{2}$) since $G_b(\frac{Q}{2}-it)$ is unitary. Furthermore we have $\hat{A}\hat{B}=q^{-2}\hat{B}\hat{A}$. Moreover, both $q\inv BA\inv$ and $q\hat{B}\hat{A}\inv$ are positive operators, hence the expressions are well-defined. 

The formula for $W$ follows from $W^*$ by conjugation, and using the relations \eqref{comW1}, \eqref{comW2}.
For simplicity, we compute $W^*$ formally using the definition:
\begin{eqnarray*}
&&W^*(f\ox g) \\
&=& \D(g)(f\ox 1)\\
&=&\left(g(s_1,t_1+t_2)\frac{G_b(-it_1)G_b(-it_2)}{G_b(-i(t_1+t_2))}\d(s_2,s_1+t_1)\right)\cdot f(s_1,t_1)\\
&=&\int_C\int_\R g(s_1-s_1',t_1+t_2-t_1')\frac{G_b(it_1'-it_1)G_b(-it_2)}{G_b(-i(t_1+t_2-t_1'))}\\
&&\tab\d(s_2,s_1+t_1-t_1'-s_1') f(s_1',t_1')e^{2\pi is_1'(t-t_1')} ds_1'dt_1',\\
\end{eqnarray*}
where $C$ goes below $t_1'=t_1$.
Replacing $s_1'= s_1+t_1-s_2-t_1'$, we get
\begin{eqnarray*}
&=&\int_C f(s_1+t_1-s_2-t_1',t_1')g(t_1'-t_1+s_2,t_1+t_2-t_1')\\
&&\tab e^{2\pi i(s_1+t_1-s_2-t_1')(t-t_1')}\frac{G_b(it_1'-it_1)G_b(-it_2)}{G_b(-i(t_1+t_2-t_1'))} dt_1'.\\
\end{eqnarray*}
Renaming $t_1'=\t$ and do a shift $\t\mapsto t_1-\t$ we get
\begin{eqnarray*}
W^*(f\ox g)&=&\int_{\R+i0} f(s_1-s_2+\t,t_1-\t)g(s_2-\t,t_2+\t)e^{2\pi i(s_1-s_2+\t)\t}\frac{G_b(-i\t)G_b(-it_2)}{G_b(-i\t-it_2)} d\t,\\
\end{eqnarray*}
where the contour of $\t$ now goes above $\t=0$.

On the other hand, using Corollary \ref{gbint}, the formula we desired is
\begin{eqnarray*}
&&e^{\frac{i}{2\pi b^2}\log A\ox \log \hat{A}}g_b^*(q\inv BA\inv \ox \hat{B}) (f\ox g)\\
&=&e^{\frac{i}{2\pi b^2}\log A\ox \log \hat{A}}\left(\int_{\R+i0} A^{-ib\inv \t}B^{ib\inv \t}\ox \hat{B}^{ib\inv \t}G_b(-i\t)d\t\right) (f\ox g)\\
&=&e^{-2\pi ip_{s_1} s_2}\left(\int_{\R+i0} f(s_1+\t,t_1-\t)g(s_2-\t,t_2+\t)e^{2\pi i(s_1+\t)\t}\frac{G_b(-it_2)G_b(-i\t)}{G_b(-it_2-i\t)}d\t\right)\\
&=&\int_{\R+i0} f(s_1-s_2+\t,t_1-\t)g(s_2-\t,t_2+\t)e^{2\pi i(s_1-s_2+\t)\t}\frac{G_b(-i\t)G_b(-it_2)}{G_b(-i\t-it_2)} d\t.\\
\end{eqnarray*}
For completeness, let us note that the action of $W$ is given by
\Eq{\label{Waction}W(f\ox g)=\int_{\R+i0}\frac{G_b(-it_2)e^{2\pi i\t s_1}}{G_b(-it_2-i\t)G_b(Q+i\t)}f(s_1+s_2,t_1-\t)g(s_2-\t,t_2+\t)d\t.}
\end{proof}

In \cite{WZ}, it is commented that manageability of the multiplicative unitary is the property that distinguish quantum groups from quantum semigroups. However, we obtain the following:
\begin{Thm} $W$ defined above is manageable.
\end{Thm}
\begin{proof}
We need to look for the positive operator $Q$ and the unitary $\til[W]$ in Definition \ref{manageable}. It is known (e.g. see \cite{T, KV1}) that $Q$ can be expressed in terms of the scaling constant and the scaling group:
\Eq{Q=P^{1/2}, \tab P^{it}\L(x)= \nu^{1/2}\L(\t_t(a)).}
From the expressions derived in the previous sections, we obtain
\Eq{Q\cdot f(s,t) = e^{-\pi Q(t-\frac{iQ}{2})}f(s,t).}
Note that in fact $Q$ is positive, under the inner product of $\cH$ (with a shift in $t\mapsto t+\frac{iQ}{2}$).

Then we see that $Q\ox Q$ commutes with $A\ox1, 1\ox\hat{A}$ and $B\ox \hat{B}$. Therefore from the formula \eqref{W1} we conclude that $Q\ox Q$ commutes with $W$.

Next we will express $\til[W]$ from definition, and show that it is a unitary operator. First of all, let us note that restricting to holomorphic functions, when a shift in contour is involved, the inner product actually reads
\Eq{\<f,g\>_{\cH}=\int_\R\int_\R f(s,t+\frac{iQ}{2})\overline{g(\over[s],\over[t]+\frac{iQ}{2})}dsdt,}
see the derivation in the proof of Theorem \ref{GNSmap}.

By definition we have
$$\<f'\ox g', W(f\ox g)\>_{\cH\ox\cH}=\< \over[f]\ox Qg', \til[W](\over[f]'\ox Q\inv g)\>_{\over[\cH]\ox \cH},$$
where $\over[f]:= \over[f](z) = \overline{f(\over[z])}\in \over[\cH]$, and the inner product is given by
\begin{eqnarray}
\<\over[f],\over[g]\>_{\over[\cH]}&:=&\<g, f\>_{\cH}\nonumber\\
&=& \int_\R\int_\R \over[f](s,t-\frac{iQ}{2})g(s,t+\frac{iQ}{2})ds dt\nonumber\\
&=&\int_\R\int_\R\over[f](s,t-\frac{iQ}{2})\overline{\over[g](\over[s],\over[t]-\frac{iQ}{2})} ds dt.
\end{eqnarray}
Now using the action of $W$ given by \eqref{Waction}, we have
\begin{eqnarray*}
LHS&=&\iint\iint ds_1dt_1ds_2dt_2 f'(s_1,t_1+\frac{iQ}{2})g'(s_2,t_2+\frac{iQ}{2})\cdot\\
&&\overline{\int_{\R+i0} \frac{G_b(\frac{Q}{2}-i\over[t_2])e^{2\pi i\t \over[s_1]}}{G_b(\frac{Q}{2}-i\over[t_2]-i\t)G_b(Q+i\t)}f(\over[s_1]+\over[s_2],\over[t_1]-\t+\frac{iQ}{2})g(\over[s_2]-\t,\over[t_2]+\t+\frac{iQ}{2})d\t}\\
&=&\iint\iint \int_{\R+i0}d\t ds_1dt_1ds_2dt_2 f'(s_1,t_1+\frac{iQ}{2})g'(s_2,t_2+\frac{iQ}{2})e^{-2\pi i\t s_1}\cdot\\
&&\frac{G_b(\frac{Q}{2}-it_2-i\t)G_b(i\t)}{G_b(\frac{Q}{2}-it_2)}\over[f](s_1+s_2,t_1-\t-\frac{iQ}{2})\over[g](s_2-\t,t_2+\t-\frac{iQ}{2})\\
&=&\iint\iint \int_{\R+i0}d\t ds_1dt_1ds_2dt_2 \over[f](s_1,t_1-\frac{iQ}{2})g'(s_2,t_2+\frac{iQ}{2})e^{-2\pi i\t (s_1-s_2)}\cdot\\
&&\frac{G_b(\frac{Q}{2}-it_2-i\t)G_b(i\t)}{G_b(\frac{Q}{2}-it_2)}f'(s_1-s_2,t_1+\t+\frac{iQ}{2})\over[g](s_2-\t,t_2+\t-\frac{iQ}{2})\\
RHS&=&\iint\iint ds_1dt_1ds_2dt_2\over[f](s_1,t_1-\frac{iQ}{2}g'(s_2,t_2+\frac{iQ}{2})e^{-\pi Qt_2}\\
&&\overline{\til[W](\over[f]'\ox Q\inv g)(\over[s_1],\over[t_1]-\frac{iQ}{2},\over[s_2],\over[t_2]+\frac{iQ}{2})}.\\
\end{eqnarray*}
Comparing, we get
\begin{eqnarray*}
\overline{\til[W](\over[f]'\ox Q\inv g)(\over[s_1],\over[t_1]-\frac{iQ}{2},\over[s_2],\over[t_2]+\frac{iQ}{2})}&=&\int_{\R+i0}e^{\pi Qt_2}e^{-2\pi i\t(s_1-s_2)}\frac{G_b(\frac{Q}{2}-it_2-i\t)G_b(i\t)}{G_b(\frac{Q}{2}-it_2)}\\
&&f'(s_1-s_2,t_1+\t+\frac{iQ}{2})\over[g](s_2-\t,t_2+\t-\frac{iQ}{2})d\t\\
\til[W](\over[f]'\ox Q\inv g)(\over[s_1],\over[t_1]-\frac{iQ}{2},\over[s_2],\over[t_2]+\frac{iQ}{2})&=&\int_{\R+i0}e^{\pi Q\over[t_2]}e^{2\pi i\t(\over[s_1]-\over[s_2])}\frac{G_b(\frac{Q}{2}-i\over[t_2])}{G_b(\frac{Q}{2}-i\over[t_2]-i\t)G_b(Q+i\t)}\\
&&\over[f'](\over[s_1]-\over[s_2],\over[t_1]+\t-\frac{iQ}{2})g(\over[s_2]-\t,\over[t_2]+\t+\frac{iQ}{2})d\t\\
\end{eqnarray*}
so that
$$\til[W](\over[f]'\ox Q\inv g)(s_1,t_1,s_2,t_2)=\int_{\R+i0}\frac{e^{\pi Q(t_2-\frac{iQ}{2})}e^{2\pi i\t(s_1-s_2)}G_b(-it_2)}{G_b(-it_2-i\t)G_b(Q+i\t)}\over[f]'(s_1-s_2,t_1+\t)g(s_2-\t,t_2+\t)d\t.$$
Now renaming $G=Q\inv g$, i.e. $G(s_2,t_2)=e^{\pi Q(t_2-\frac{iQ}{2})}g(s_2,t_2)$, we have
$$g(s_2-\t,t_2+\t)=G(s_2-\t,t_2+\t)e^{-\pi Q(t_2+\t-\frac{iQ}{2})},$$
hence the action of $\til[W]$ becomes
\begin{eqnarray*}
\til[W](\over[f]\ox g)&=&\int_{\R+i0}\frac{e^{-\pi Q\t}e^{2\pi i\t(s_1-s_2)}G_b(-it_2)}{G_b(-it_2-i\t)G_b(Q+i\t)}\over[f](s_1-s_2,t_1+\t)g(s_2-\t,t_2+\t)d\t\\
&=&\int_{\R+i0}\over[f](s_1-s_2,t_1+\t)g(s_2-\t,t_2+\t)e^{2\pi i\t(s_1-s_2+t_2+\t-\frac{iQ}{2})}\\
&&\frac{G_b(Q+it_2+i\t)G_b(-i\t)}{G_b(Q+it_2)}d\t.\\
\end{eqnarray*}
Hence we can conclude that
\begin{eqnarray*}
\til[W]&=&e^{-2\pi is_2 p_{s_1}}\int_{\R+i0} (B^\#)^{ib\inv \t} \ox (\hat{B}^\#)^{ib\inv \t}(\hat{A}^{\#})^{ib\inv \t}G_b(-i\t)d\t\\
&=&e^{\frac{i}{2\pi b^2}\log A\inv \ox \log\hat{A}}g_b^*(B^\#\ox q\hat{B}^\#\hat{A}^\#),
\end{eqnarray*}
where 
$$B^\#=e^{2\pi b (s+p_{t})},\tab \hat{A}^\#=e^{2\pi b(t-\frac{iQ}{2})},$$
$$\hat{B}^\#=G_b(Q+it)\inv e^{2\pi b(p_{t}-p_{s})}\circ G_b(Q+it),$$
are all positive operators with respect to the measure in $\over[\cH]\ox\cH$, and we used an analogue of \eqref{gb3} with $\hat{A}^\#\hat{B}^\#=q^2\hat{B}^\#\hat{A}^\#$. Hence $\til[W]$ is a unitary operator.
\end{proof}
\begin{Rem} This is a rather curious and striking result since our quantum plane has been restricted to $B>0$ so that all the operators involved are positive essentially self-adjoint and nicely defined. In fact, the key difference is that in \cite{WZ}, the multiplicative unitary $W$ lies in $\cA\ox \cA$, while our $W$ lies in $\cA\ox \hat{\cA}$, and the Hilbert space $\cH\simeq L^2(\R\x\R)$ is ``twice" larger than the canonical space $L^2(\R)$ for the action of $\cA$ only. For example, we see that the operator $\hat{A}^\#=e^{2\pi b(t-\frac{iQ}{2})}$ does not lie in the action spanned by $\cA$ itself, which only contains actions on the variable $s$ and the shifting operator $e^{2\pi bp_t}$. Therefore it is worth studying a deeper meaning of manageability of $W$ in this wider context of quantum semigroups. The difference in the definition of $W$ mentioned above in fact produces a new transformation formula for $G_b$, see Section \ref{sec:meaning} for further detail.
\end{Rem}

\section{The dual space}\label{sec:Dual} We have encountered the dual space elements in the expression for the multiplicative unitary $W$. Following the recipe in \cite{KV1}, in order to describe the GNS representation for the dual space $\hat{\cA}$, which is self-dual for quantum plane, we need to establish a non-degenerate pairing between $\cA$ and $\hat{\cA}$. 

First of all let us describe its own GNS representation $(\hat{\cH},\hat{\pi},\hat{\L})$. After establishing the pairing, we can then relate them to the canonical space $\cH$.

\subsection{Definitions}\label{sec:Dual:def}
First let us describe the dual space on the Hopf algebra level.
\begin{Def} The algebraic dual space $\cA_q^*\simeq \cA_q$ is generated by self-adjoint operators $X,Y$ with $XY=q^2YX$ and the same coproduct.
\end{Def}
Following \cite{FK}, occasionally we will use the notation $\cB_q$ to denote $\cA_q^*$ when we want to think about the dual space as the quantum \emph{algebra} counterpart of the quantum \emph{group} $\cA_q$. This correspondence is discussed in \cite{Ip}. Similarly we use the notation $\cB_{q\til[q]}$ to denote the modular double $\cB_q\ox \cB_{\til[q]}$.

\begin{Def} We define $\hat{\cA}_q=\cA_q^{*op}$ to be the algebra generated by $X,Y$ with $XY=q^2YX$ and the opposite coproduct. Alternatively, by defining $\hat{A}=X\inv, \hat{B}=~q\inv YX\inv$, we define $\hat{\cA}_q$ to be generated by positive elements $\hat{A},\hat{B}$ with $\hat{A}\hat{B}=q^{-2}\hat{B}\hat{A}$ and the same coproduct as $A$ and $B$:
\begin{eqnarray}
\D(\hat{A})&=&\hat{A}\ox\hat{A},\\
\D(\hat{B})&=&\hat{B}\ox\hat{A}+1\ox \hat{B}.
\end{eqnarray}
\end{Def}

\begin{Def} Similar to the quantum plane, we define $\hat{\cA}:=\cC_\oo(\hat{\cA}_q)$ to be the closure of the linear span of elements of the form
\Eq{\hat{f}=\int_{\R+i0}\int_\R \hat{f}(s,t)\hat{B}^{ib\inv t}\hat{A}^{ib\inv s} dsdt,}
where $\hat{f}(s,t)$ has the same analytic properties as those in $\cA$.
\end{Def}

Since the spaces are self-dual, we immediately have the following
\begin{Prop} The left and right Haar functionals are given as before:
\begin{eqnarray}
\hat{h}_L(\hat{f}) &=& \hat{f}(0,iQ),\\
\hat{h}_R(\hat{f}) &=& \hat{f}(-iQ,iQ)e^{-\pi iQ^2}.
\end{eqnarray}
The GNS representation on $\hat{\cH}$ associated with $\hat{h}_L$ is simply given by
\Eq{\hat{\L}_L\left(\int_{\R+i0}\int_\R \hat{f}(s,t)\hat{B}^{ib\inv t}\hat{A}^{ib\inv s} dsdt\right)=f(s,t),}
and the inner product is given by (notice the extra exponent):
\Eq{\<\hat{f},\hat{g}\>:=\hat{h}_L(\hat{g}^*\hat{f})=\int_\R\int_\R \overline{\hat{g}(s,t+\frac{iQ}{2})}\hat{f}(s,t+\frac{iQ}{2})e^{2\pi Qs}dsdt.}
\end{Prop}

Following the same method as in the previous section, we also establish the following maps:
\begin{Prop} The product of two elements is given by
\Eq{(\hat{f}\cdot \hat{g})(s,t) = \iint \hat{f}(s-s',t-t')\hat{g}(s',t')e^{2\pi i(s-s')t} ds'dt'.}
The action of $\hat{A}$ and $\hat{B}$ is given by:
\Eq{\hat{\pi}(\hat{A})=e^{2\pi bt}e^{-2\pi bp_s},\tab \hat{\pi}(\hat{B})=e^{-2\pi bp_t}.}
The antipode is given by
\Eq{\hat{S}(\hat{f})=\hat{f}(-s-t,t)e^{\pi Qt}e^{\pi i t^2+2\pi is t}.}
The scaling group  is given by
\Eq{\hat{\t}_t(\hat{A})=\hat{A}\tab \hat{\t}_t(\hat{B})=e^{2\pi bQt}\hat{B}.}
The unitary antipode  is given by
\Eq{\hat{R}(\hat{A})=\hat{A}\inv\tab \hat{R}(\hat{B})=q\hat{B}\hat{A}\inv ,}
or explicitly by
\Eq{\hat{R}(\hat{f})=\hat{f}(-s-t,t)e^{\pi it^2+2\pi ist}.}
The scaling constant  is given by
\Eq{\hat{\nu}=e^{2\pi Q^2}=\nu\inv.}
Note that by general theory, for the dual space we have instead $\hat{S}=\hat{R}\hat{\t}_{i/2}$.
\end{Prop}

\begin{Prop} Acting on the space $\hat{\cH}$, we have
\begin{eqnarray}
\hat{T}:\hat{f}(s,t)&\mapsto& \over[\hat{f}](-s,-t)e^{2\pi ist},\\
\hat{T}^*:\hat{f}(s,t)&\mapsto& \over[\hat{f}](-s,-t)e^{2\pi ist}e^{-2\pi sQ},\\
\hat{\nabla}:\hat{f}(s,t)&\mapsto& \hat{f}(s,t)e^{-2\pi sQ},\\
\hat{J}:\hat{f}(s,t)&\mapsto&\over[\hat{f}](-s,-t)e^{2\pi ist}e^{-\pi sQ}.
\end{eqnarray}
\end{Prop}

\subsection{The non-degenerate pairing}\label{sec:Dual:pairing}
Recall that a Hopf pairing between two Hopf algebras $\cA$ and $\cA^*$ is a non-degenerate pairing such that
\Eq{\label{pairingdef}\<a, xy\> = \<\D(a),x\ox y\>,\tab \<ab,x\>=\< a\ox b, \D(x)\>.}
for $a,b\in\cA, x,y\in\cA^*$. (In this paper we will not consider counit and antipode in the pairing.)

In the current setting, we have the pairing between $\cA_q$ and $\cA_q^*$ from the compact case:
$$\<A,X\>=q^{-2},\tab \<B,X\>=0,$$
$$\<A,Y\>=0,\tab \<B,Y\>=c\in\C,$$
where $c$ is any complex number.

The pairings between $A,B$ and $\hat{A}=X\inv,\hat{B}=q\inv YX\inv$ are given by
$$\<A,\hat{A}\>=q^{2},\tab \<B,\hat{A}\>=0,$$
$$\<A,\hat{B}\>=0,\tab \<B,\hat{B}\>=cq:=c'.$$
We will choose $c'=1$.
From this it is extended to
$$\< B^mA^n,\hat{A}\>=q^{2n}\d_{m0},\tab\<B^mA^n ,\hat{B}\>=\d_{m1},$$
and then to
$$\<B^mA^n ,\hat{B}^{m'}\hat{A}^{n'}\>=q^{2n'(n+m)}[m]_q! \d_{mm'},$$
where $[m]_q!:=[m]_q[m-1]_q...[1]_q$ is the $q$-factorial, and $[m]_q=\frac{q^{m}-q^{-m}}{q-q\inv}$ is the $q$-number.

From the analogy between the quantum dilogarithm $G_b$ and the $\G_q$ function established in \cite{Ip}, we prove the following theorem:

\begin{Thm} The Hopf pairing between $\cA$ and $\cA^*\simeq \hat{\cA}^{cop}$ is given by
\begin{eqnarray}
&&\label{pairing1}\<\iint g(s,t) B^{ib\inv t}A^{ib\inv s}dsdt, \iint \hat{f}(s,t)\hat{B}^{ib\inv t}\hat{A}^{ib\inv s}dsdt\>\nonumber\\
&=&\iiint g(s,t)\hat{f}(s',t)G_b(Q+it)e^{-2\pi is'(s+t)} ds'dsdt\\
&&\<\iint g(s,t) A^{ib\inv s}B^{ib\inv t}dsdt, \iint \hat{f}(s,t)\hat{B}^{ib\inv t}\hat{A}^{ib\inv s}dsdt\>\nonumber\\
&=&\label{pairing}\iiint g(s,t)\hat{f}(s',t)G_b(Q+it)e^{-2\pi is'(s+t)}e^{-2\pi ist}ds'dsdt,
\end{eqnarray}
or using the Mellin transformed picture (see Remark \ref{MellinQ}):
\begin{eqnarray}
&&\<\int f(t) B^{ib\inv t}g(\log A)dt, \iint \hat{f}(s,t)\hat{B}^{ib\inv t}\hat{A}^{ib\inv s}dsdt\>\nonumber\\
&=&\label{pairing2}\iint g(-2\pi bs')\hat{f}(s',t)f(t)G_b(Q+it)e^{-2\pi its'}ds'dt\\
&&\<\int g(\log A)f(t) B^{ib\inv t}dt, \iint \hat{f}(s,t)\hat{B}^{ib\inv t}\hat{A}^{ib\inv s}dsdt\>\nonumber\\
&=&\iint g(-2\pi b(s'+t))\hat{f}(s',t)f(t)G_b(Q+it)e^{-2\pi its'}ds'dt.
\end{eqnarray}
Here the Hopf pairing between two $C^*$-algebra are defined by naturally extending the pairing in the defining relations \eqref{pairingdef} to a pairing between the multiplier algebras $M(\cA\ox\cA)$ and $M(\cA^*\ox\cA^*)$.
\end{Thm}

\begin{proof} It suffices to show that 
\Eq{\< f, \hat{g}\hat{h}\> = \<\D(f),\hat{g}\ox\hat{h}\>.}
The other relations are similar by duality.
We will prove the first form \eqref{pairing}. For simplicity we omit all the integrations after the pairing. Every variable is to be integrated.
\begin{eqnarray*} 
LHS&=&\<f(s,t), \iint \hat{g}(s-s',t-t')\hat{h}(s',t')e^{2\pi i(s-s')t'} ds'dt'\>\\
&=& f(s,t)\hat{g}(s''-s',t-t')\hat{h}(s',t')e^{2\pi i(s''-s')t'}G_b(Q+it)e^{-2\pi is'' (s+t)}e^{-2\pi ist}\\
&=& f(s,t)\hat{g}(s''-s',t-t')\hat{h}(s',t')G_b(Q+it)e^{2\pi i(s''-s')t'}e^{-2\pi is'' (s+t)}e^{-2\pi ist},\\
RHS&=&\<f(s_1,t_1+t_2)\frac{G_b(-it_1)G_b(-it_2)}{G_b(-it_1-it_2))}\d(s_2,s_1+t_1),\hat{g}(s_1,t_1)\hat{h}(s_2,t_2)\>\\
&=&f(s_1,t_1+t_2)\frac{G_b(-it_1)G_b(-it_2)}{G_b(-it_1-it_2)}G_b(Q+it_1)G_b(Q+it_2) \d(s_2,s_1+t_1)\\
&&\hat{g}(s_1',t_1)\hat{h}(s_2',t_2)e^{-2\pi is_1'(s_1+t_1)}e^{-2\pi is_2'(s_2+t_2)}e^{-2\pi is_1t_1}e^{-2\pi is_2t_2}\\
&=&f(s_1,t_1+t_2)\hat{g}(s_1',t_1)\hat{h}(s_2',t_2)G_b(Q+it_1+it_2)e^{-2\pi is_2't_2}e^{-2\pi is_1(t_1+t_2)}e^{-2\pi i(s_1+t_1)(s_1'+s_2')},\\
\end{eqnarray*}
where we used the reflection properties \eqref{reflection} of the $G_b$ function: 
\Eq{\frac{G_b(-is)G_b(-it)}{G_b(-is-it)}=\frac{G_b(Q+is+it)e^{2\pi ist}}{G_b(Q+is)G_b(Q+it)}.}
Now shifting $t_1\mapsto t_1-t_2,s_1'\mapsto s_1'-s_2'$, we obtain
\begin{eqnarray*}
&=&f(s_1,t_1)\hat{g}(s_1'-s_2',t_1-t_2)\hat{h}(s_2',t_2)G_b(Q+it_1)e^{-2\pi is_2't_2}e^{-2\pi is_1t_1}e^{-2\pi is_1'(s_1+t_1-t_2)}\\
&=&f(s_1,t_1)\hat{g}(s_1'-s_2',t_1-t_2)\hat{h}(s_2',t_2)G_b(Q+it_1)e^{2\pi i(s_1'-s_2')t_2}e^{-2\pi is_1'(s_1+t_1)}e^{-2\pi is_1t_1}.\\
\end{eqnarray*}
Hence we see that the expression is equal on renaming $s_1=s,t_1=t,s_1'=s'', s_2'=s'$ and $t_2=t'$.
\end{proof}
With the pairing established, we can prove the density condition needed in order for $\cA$ to be a locally compact quantum group.

\begin{Thm}\label{density} We have
\begin{eqnarray}
\cA&=&span\{(\w\ox1)\D(a)|\w\in\cA^*,a\in\cA\}^{closure}\\
&=&span\{(1\ox\w)\D(a)|\w\in\cA^*,a\in\cA\}^{closure}.
\end{eqnarray}
\end{Thm}
\begin{proof}
We will prove that $(\w\ox 1)\D(a)$ is dense in $1\ox\cA$, while the other statement is similar. Let us write
\begin{eqnarray*}
\w&=&\iint \hat{f}(s,t)\hat{B}^{ib\inv t}\hat{A}^{ib\inv s}dsdt\in\cA^*\\
a&=&\iint g(s,t) A^{ib\inv s}B^{ib\inv t}dsdt\in\cA.
\end{eqnarray*}
Using Proposition \ref{coprod}, the required pairing is (here we understood that it is an element in $1\ox\cA$ with coordinates $s_2,t_2$, and that all integrations converge absolutely with the appropriate contours):
\begin{flalign*}
&\<g(s_1,t_1+t_2)\frac{G_b(-it_1)G_b(-it_2)}{G_b(-i(t_1+t_2))}\d(s_2,s_1+t_1), \hat{f}(s,t)\>&\\
&=\iiint \hat{f}(s',t_1)g(s_1,t_1+t_2)\frac{G_b(-it_1)G_b(-it_2)}{G_b(-i(t_1+t_2))}G_b(Q+it_1)\d(s_2,s_1+t_1)e^{-2\pi is'(s_1+t_1)}ds'ds_1dt_1&\\
&=\iint \hat{f}(s',t_1)g(s_2-t_1,t_1+t_2)\frac{G_b(-it_2)}{G_b(-i(t_1+t_2))}e^{-2\pi is's_2}e^{\pi it_1^2-\pi t_1Q}ds'dt_1.&
\end{flalign*}
Now if we choose $\hat{f}(s',t_1)=\hat{f_1}(s')\hat{f_2}(t_1)$, the integration over $s'$ is just the Fourier transform $\hat{\cF}_1(s_2)$ of $\hat{f_1}(s')$, and we get
\begin{flalign*}
&=\int \hat{\cF}_1(s_2)\hat{f_2}(t_1)g(s_2-t_1,t_1+t_2)\frac{G_b(-it_2)}{G_b(-i(t_1+t_2))}e^{\pi it_1^2-\pi t_1Q}dt_1.&
\end{flalign*}
Finally we can choose any nice approximation of identity for $\hat{f_2}(t_1)\to\d(t_1)$, then in the limit we obtain simply
\begin{flalign*}
&\to\hat{\cF}_1(s_2)g(s_2,t_2),&
\end{flalign*}
which is then obviously all of $1\ox\cA$ for different choices of $\hat{\cF}$ and $g$.
\end{proof}
\begin{Cor} $\cA$ is a locally compact quantum group in the sense of \cite{KV1,KV2}.
\end{Cor}
\subsection{The GNS description}\label{sec:Dual:GNS}
Given our multiplicative unitary 
\Eq{W=e^{\frac{i}{2\pi b^2}\log A\inv\ox \log \hat{A}}g_b(B \ox q\hat{B}\hat{A}\inv)\in M(\cA \ox \cB(\cH)),} 
the space $\hat{\cA}\sub \cB(\cH)$ is originally defined to be
\Eq{ \hat{\cA}=\{ (\w\ox 1)W: \w\in \cA^*\}^{\mbox{norm closure}}.}
In order to find the GNS representation $\hat{\L}$ of $\hat{\cA}$ on the original space $\cH_L$, we introduce the contraction map $\l$ \cite[(8.2)]{KV1} which relates $\cA^*$ to the definition of $\hat{\cA}$ above using the multiplicative unitary, and the $\xi$ map \cite[Notation 8.4]{KV1}, which relates $\cA^*$ to the original GNS representation space $\cH_L$ of $\cA$ by the Riesz's theorem for Hilbert spaces. Then $\hat{\L}$ is just the composition of the two maps, which is naturally compatible with the previous GNS construction for $\cA$.
\begin{Def}
For $\w\in \cA^*, x\in\cA$, we define
\begin{eqnarray}
\l&:&\cA^*\to \hat{\cA}\nonumber\\
\l(\w)&=&(\w\ox 1)W.
\end{eqnarray}
and
\begin{eqnarray}
\xi&:&\cA^*\to \cH_L\nonumber\\
\w(x^*)&=&\<\xi(\w),\L(x)\>_L
\end{eqnarray}
Then the GNS map $\hat{\L}$ is given by
\begin{eqnarray}
\hat{\L}&:&\hat{\cA}\to \cH_L\nonumber\\
\hat{\L}(\l(\w))&=&\xi(\w).
\end{eqnarray}
\end{Def}

\begin{Prop}\label{lamda} $\l$ is the identity map from $\cA^*$ to $\hat{\cA}$ as elements in $\cB(\cH)$.
\end{Prop}
\begin{proof} Using the second form \eqref{W2} for $W$, and the integral form \eqref{gb2} for $g_b$, we get
$$W=\left(\int_{\R+i0} \frac{B^{ib\inv\t}A^{-ib\inv\t}\ox \hat{B}^{ib\inv\t}}{G_b(Q+i\t)}d\t\right)e^{\frac{i}{2\pi b^2}\log A\inv\ox\log \hat{A}}.$$
By the pairing \eqref{pairing} we have
\begin{eqnarray*}
\l(\w)&=&\left(\iint \hat{f}(s,t)\hat{B}^{ib\inv t}\hat{A}^{ib\inv s}dsdt\ox 1\right) (W)\\
&=&  \iint \hat{f}(s,t) \frac{(q^{2ib\inv s})^{-ib\inv t}}{G_b(Q+it)} G_b(Q+it)e^{-2\pi ist}\hat{B}^{ib\inv t}  \hat{A}^{\frac{i}{2\pi b^2} (2\pi bs)}dsdt\\
&=& \iint \hat{f}(s,t)\hat{B}^{ib\inv t} \hat{A}^{ib\inv s}dsdt.
\end{eqnarray*}
\end{proof}

\begin{Prop}
The map $\xi:\cA^*\sub\hat{\cH}\to \cH$ is given by:
\Eq{\xi:\hat{f}(s,t)\mapsto F(s,t):=\int_\R \hat{f}(s',iQ-t)G_b(-it)e^{2\pi i s'(s+t-iQ)}ds',}
it naturally extends to an invertible map $\xi:\hat{\cH}\to\cH$ with the inverse given by
\Eq{\xi\inv:F(s,t)\mapsto \hat{f}(s,t):=\int_\R \frac{F(s',iQ-t)}{G_b(Q+it)}e^{2\pi is(t-s')}ds'.}
\end{Prop}
\begin{proof}
Let $\w=\iint \hat{f}(s,t)\hat{B}^{ib\inv t}\hat{A}^{ib\inv s}dsdt, x= \iint g(s,t)A^{ib\inv s}B^{ib\inv t}dsdt$. Then
\begin{eqnarray*}
\w(x^*)&=&\<\xi(\w),\L(x)\>_L\\
LHS&=& \iiint \over[g](-s,-t)e^{2\pi ist}F(s',t)G_b(Q+it)e^{-2\pi is's-2\pi is't-2\pi ist}ds'dsdt\\
&=& \iiint \over[g](-s,-t)F(s',t)G_b(Q+it)e^{-2\pi is's-2\pi is't}ds'dsdt\\
RHS&=&\iint \overline{g(\over[s],\over[t]+\frac{iQ}{2})}\xi(\w)(s,t+\frac{iQ}{2})dsdt\\
&=&\iint \over[g](s,t-\frac{iQ}{2})\xi(\w)(s,t+\frac{iQ}{2})dsdt\\
&=&\iint \over[g](-s,-t)\xi(\w)(-s,-t+iQ)dsdt.
\end{eqnarray*}
Hence
\begin{eqnarray*}
\xi(\w)(s,t)&=&\int F(s',iQ-t)G_b(Q+i(iQ-t))e^{2\pi is's-2\pi is'(iQ-t)}ds'\\
&=&\int F(s',iQ-t)G_b(-it)e^{2\pi i s'(s+t-iQ)}ds'.
\end{eqnarray*}
We observe that $\xi$ can be realized as
$$\xi=\cF_s\circ P_t\circ G_b(Q+it)e^{-2\pi ist},$$
where $\cF_s$ is Fourier transform on $s$, $P_t f(t):=f(iQ-t)$ and $G_b(Q+it)e^{-2\pi ist}$ is just multiplication by this function. Hence we have
$$\xi\inv=G_b(Q+it)\inv e^{2\pi ist}\circ P_t\inv\circ \cF_s\inv,$$
or the formula desired.
\end{proof}

\begin{Prop} $\xi$ is an isometry, $\xi:\hat{\cH}\to \cH$.
\end{Prop}
\begin{proof} We have
\begin{eqnarray*}
&&\<\xi(\hat{f}),\xi(\hat{g})\>_L\\
&=&\int_\R\int_\R \overline{\hat{g}(s'',iQ-(t+\frac{iQ}{2}))G_b(-i(t+\frac{iQ}{2}))e^{2\pi is''(s+(t+\frac{iQ}{2})-iQ)}}\\
&& \hat{f}(s',iQ-(t+\frac{iQ}{2}))G_b(-i(t+\frac{iQ}{2}))e^{2\pi is'(s+(t+\frac{iQ}{2})-iQ)}ds'ds''dsdt\\
&=&\int_\R\int_\R \overline{\hat{g}(s'',-t+\frac{iQ}{2})G_b(\frac{Q}{2}-it)}e^{-2\pi is''(s+t+\frac{iQ}{2})}\\
&& \hat{f}(s',-t+\frac{iQ}{2})G_b(\frac{Q}{2}-it)e^{2\pi is'(s+t-\frac{iQ}{2})}ds'ds''dsdt\\
&=&\int_\R\int_\R \overline{\hat{g}(s'',t+\frac{iQ}{2})}\hat{f}(s',t+\frac{iQ}{2})e^{-2\pi is''(s-t+\frac{iQ}{2})}e^{2\pi is'(s-t-\frac{iQ}{2})}ds'ds''dsdt\\
&=&\int_\R\int_\R \overline{\hat{g}(s'',t+\frac{iQ}{2})}\hat{f}(s',t+\frac{iQ}{2})e^{-2\pi is''(-t+\frac{iQ}{2})}e^{2\pi is'(-t-\frac{iQ}{2})}e^{2\pi is(-s''+s')}ds'ds''dsdt\\
&=&\int_\R\int_\R \overline{\hat{g}(s,t+\frac{iQ}{2})}\hat{f}(s,t+\frac{iQ}{2})e^{-2\pi is(-t+\frac{iQ}{2})}e^{2\pi is(-t-\frac{iQ}{2})}dsdt\\
&=&\int_\R\int_\R\overline{\hat{g}(s,t+\frac{iQ}{2})}\hat{f}(s,t+\frac{iQ}{2})e^{2\pi Qs}dsdt\\
&=&\int_\R\int_\R \overline{\hat{g}(s,t+\frac{iQ}{2})}\hat{f}(s,t+\frac{iQ}{2})e^{2\pi Qs}dsdt\\
&=&\<\hat{f},\hat{g}\>_{\hat{L}}.
\end{eqnarray*}
\end{proof}

Hence we conclude that
\begin{Cor} $\hat{\L}(\w):=\xi(\w)$ gives the GNS representation for $\hat{\cA}$ on $\cH$ as desired.
\end{Cor}

\subsection{The multiplicative unitary and modular maps}\label{sec:Dual:multi}
Under the transformation $\xi$, we can now express all the operators defined earlier on $\hat{\cH}$ to $\cH$. We have
\begin{Prop}
The action $\pi: \hat{\cA}\to \cH$ is given by
\begin{eqnarray}
\pi(\hat{A}) &=& e^{-2\pi bs},\\
\pi(\hat{B}) &=& G_b(-it) \circ e^{2\pi b(p_t-p_s)} \circ G_b(-it)\inv.
\end{eqnarray}
\end{Prop}

Hence in particular $W$ is indeed a genuine element in $M(\cA\ox \hat{\cA})$.

We also note that by general theory, the multiplicative unitary for $\hat{\cA}$ is given by
\Eq{\hat{W}=\S W^* \S,}
where $\S$ is the permutation operator on the tensor product. The coproduct induces by $\hat{W}$ is precisely the one we defined earlier, that $\hat{A}$ and $\hat{B}$ transformed as how $A$ and $B$ do.

\begin{Prop} The modular maps act on $\cH$ by:
\begin{eqnarray}
\hat{T}(f) &=& \over[f](s+t-iQ,-t)e^{-\pi Qt}e^{-\pi it^2}\\
\hat{T^*}(f)&=&\over[f](s+t,-t)e^{-\pi Qt}e^{-\pi it^2}\\
\hat{\DD}(f)&=&f(s+iQ,t)\\
\hat{J}(f)&=&\over[f](s+t-\frac{iQ}{2},-t)e^{-\pi Qt}e^{-\pi it^2}.
\end{eqnarray}
\end{Prop}

From these actions, we can verify all the well-known properties between these maps on $\cH$ (see \cite{KV2} Prop 2.1, 2.11,2.12):
\begin{Prop} \label{RJ} For $x\in \cA, y\in\hat{\cA}$, we have the properties:
\begin{eqnarray}
\hat{J}\L_R(x)&=&\L(R(x)^*),\\
\hat{\DD}^{-\frac{1}{2}}\L(x) &=& \L_R(\t_{i/2}(x)),\\
\hat{T}^*\L(x)&=&\L(S\inv(x)^*),\\
\label{nu}\hat{J}J&=&\nu^{i/4}J\hat{J},\\
\label{R1}R(x)&=&\hat{J}x^*\hat{J},\\
\label{R2}\hat{R}(y)&=&Jy^*J.
\end{eqnarray}
Furthermore, if we define on $\cH$ the operator $G$ (\cite{KV1} Prop 3.22):
\Eq{G \L((h_R\ox 1)(\D(x^*)(y\ox 1)))= \L((h_R\ox 1)(\D(y^*)(x\ox1)))}
and its polar decomposition $G=IN^{1/2}$, then we have
$\hat{T}^*=G, \hat{J}=I$ and $\hat{\DD}=N\inv$.
\end{Prop}


\section{Transformation to new Hilbert spaces}\label{sec:Trans}
To prepare for the construction of the quantum double, we found it useful to introduce a transformation of the Hilbert space $\cH$ as well as $\hat{\cH}$, so that the action of $\cA$ and $\hat{\cA}$ acts as the canonical Weyl algebra. Furthermore, this transformation will bring the inner product to the canonical form for $L^2(\R)$.
\subsection{Unitary transformations}\label{sec:Trans:unit}
In this section, we list the transformations that will be used. They are all unitary transformations on $L^2(\R\x\R)$ equipped with the standard Lebesgue measure. For an operator $P$, its action under any transformation $\cT$ is given by
\Eq{P\mapsto \cT\circ P \circ \cT\inv.}

\begin{Prop} The following lists the transformations and their effects on the multiplication and differential operators $s$ and $p_s$.
$$f(s,t)\longmapsto:$$
\begin{eqnarray}
f(-s,t)&:& s\to -s,\tab p_s\to -p_s,\\
f(s,-t)&:& t\to -t,\tab p_t\to -p_t,\\
f(t,s)&:& s\corr t,\tab p_s\corr p_t,\\
f(s,t)e^{\pm2\pi ist}&:& p_s\to p_s\mp t,\tab p_t \to p_t\mp s,\\
f(s,t)e^{\pm \pi is^2}&:& p_s\to p_s\mp s,\\
f(s,t)e^{\pm \pi it^2}&:& p_t\to p_t\mp t,\\
f(s\pm t,t)&:& s\to s\pm t,\tab p_t\to p_t\mp p_s,\\
f(s,t\pm s)&:& t\to t\pm s,\tab p_s\to p_s\mp p_t,\\
\cF(f)=\int f(\t)e^{-2\pi i\t x}d\t &:&x\to -p_x,\tab p_x\to x,\\
\cF\inv(f)=\int f(\t)e^{2\pi i\t x}d\t &:& x\to p_x,\tab p_x\to -x,
\end{eqnarray}
where $\cF$ is the Fourier transform with
$$\cF^2 = -Id.$$
All the transformations above preserve the dense subspace $\cW\hat{\ox}\cW$ (cf. \eqref{WxW}).
\end{Prop}
\subsection{The representation space}\label{sec:Trans:rep}
We introduce the transformations such that $\cA$ acts on $\cH$ canonically, and $\hat{\cA}$ acts on $\hat{\cH}$ canonically. We will use capital letter to denote the transformed function. All the new Hilbert spaces are $L^2(\R\x\R)$.

\begin{Def} We define $\cT: \cH \to \cH_{rep}$ by
\Eq{\cT: f(s,t)\mapsto \int_\R f(\a, t-s+\frac{iQ}{2})e^{2\pi i\a s} d\a,}
\Eq{\cT\inv: F(s,t)\mapsto \int_\R F(\a, t+\a-\frac{iQ}{2})e^{-2\pi i\a s}d\a,}
or simply $\cT = (t\mapsto t-s)\circ \cF_s\inv\circ (t\mapsto t+\frac{iQ}{2})$.

Similarly we define $\hat{\cT}:\hat{\cH}\to\hat{\cH}_{rep}$ by
\Eq{\hat{\cT}: \hat{f}(s,t)\mapsto \int \hat{f}(\a, t-s+\frac{iQ}{2})e^{-2\pi i\a t}e^{\pi Q\a}d\a,}
\Eq{\hat{\cT}\inv: \hat{F}(s,t)\mapsto \int \hat{F}(\a,t+\a-\frac{iQ}{2})e^{2\pi is(\a+t)}d\a,}
or simply $\hat{\cT}=(t\mapsto t-s)\circ \cF_s \circ (t\mapsto t+\frac{iQ}{2})\circ e^{-2\pi is t}$.
\end{Def}

\begin{Prop} Under the transformations, both the spaces $\cH_{rep}$ and $\hat{\cH}_{rep}$ carry the usual inner product with respect to the standard Lebesgue measure:
\Eq{\<F(s,t), G(s,t)\> = \int_\R \int_\R \over[G(s,t)]F(s,t) dsdt.}
\end{Prop}

\begin{Prop} Under the transformations, the action of $\cA$ on $\cH_{rep}$ is given by
\Eq{A=e^{2\pi bs},\tab B= e^{2\pi bp_s},}
and the action of $\hat{\cA}$ on $\hat{\cH}_{rep}$ is given by
\Eq{\hat{A}=e^{-2\pi bs},\tab \hat{B}=e^{2\pi bp_s}.}
\end{Prop}
\begin{proof} Let us demonstrate the use of transformation rule for say, the operator $B$. The rest is similar.
Recall that $B$ acts on $\cH$ as $e^{2\pi b s}e^{-2\pi bp_t}$. Hence under the transformation it becomes
\begin{eqnarray*}
e^{2\pi b s}e^{-2\pi bp_t}&\mapsto_{t\mapsto t+\frac{iQ}{2}}& e^{2\pi b s}e^{-2\pi bp_t}\\
&\mapsto_{\cF_s\inv}& e^{2\pi b p_s}e^{-2\pi bp_t}\\
&\mapsto_{t\mapsto t-s}& e^{2\pi b p_s}.
\end{eqnarray*}
\end{proof}

As a corollary, we have
\begin{Cor} As a representation of $\cA$ we have
\Eq{\cH\simeq \cH_{irr} \ox L^2(\R),}
where $\cH_{irr}$ is the canonical irreducible representation \eqref{Hirr} of $\cA$ on $L^2(\R)$.
\end{Cor}

Let us also note that 
\begin{Prop} The product of two elements in $\cA$ induces a twisted product on $\cH$ by
\Eq{F\cdot G = \cT(\cT\inv(F)\cdot \cT\inv(G)) = \int_\R F(s,\t)G(\t,t) d\t.}
\end{Prop}

We can do the same analysis to the dual space:
\begin{Prop} Under the transformations, the action of $\hat{\cA}$ on $\cH_{rep}$ is given by
\Eq{\hat{A}=e^{2\pi b(p_s+p_t)},\tab \hat{B}=G_b(\frac{Q}{2}+is-it)\circ e^{2\pi b(p_t+s)}\circ G_b(\frac{Q}{2}+is-it)\inv,}
and the action of $\cA$ on $\hat{\cH}_{rep}$ is given by
\Eq{A=e^{2\pi b(p_s+p_t)}, B=G_b(\frac{Q}{2}-is+it)\inv\circ e^{2\pi b(p_t-s)}\circ G_b(\frac{Q}{2}-is+it).}
\end{Prop}
\begin{Rem} Note that the action of $\cA$ on $\hat{\cH}$ is obtained by the $\xi$ transformation. Moreover, recall that $G_b(\frac{Q}{2}+is)$ is unitary, hence the action above is still positive.
\end{Rem}

Finally we will need to describe the action of the modular conjugation $J$ and $\hat{J}$ that are crucial in the construction of the multiplicative unitary for the quantum double.

\begin{Prop} The action of $J$ and $\hat{J}$ on $\cH_{rep}$ is given by
\Eq{J: F(s,t)\mapsto \over[F](t,s),}
\Eq{\hat{J}:F(s,t)\mapsto \over[F](-s,-t)e^{\pi is^2-\pi it^2}\nu^{i/8},}
where $\nu=e^{-2\pi Q^2}$.

Similarly, the action of $\hat{J}$ on $\hat{\cH}_{rep}$ is also given by
\Eq{\hat{J}: \hat{F}(s,t)\mapsto \overline{\hat{F}}(t,s),}
while the action of $J$ on $\hat{\cH}_{rep}$ is given by
\Eq{J:\hat{F}(s,t)\mapsto \overline{\hat{F}}(-s,-t)e^{\pi it^2-\pi is^2}\nu^{-i/8}.}
\end{Prop}

\begin{Prop} \label{AdJ} The action of $Ad J$ on $\cA$ is given by:
\Eq{J(x)J=\hat{R}_*(x^*),}
where $\hat{R}_*$ is an anti-homomorphism given by
\Eq{\hat{R}_*(A)=\hat{A}_*,\tab \hat{R}_*(B)=\hat{B}_*,}
where \Eq{\hat{A}_* =e^{2\pi bt},\tab \hat{B}_*=e^{-2\pi bp_t}.}
The action of $Ad \hat{J}$ on $\hat{\cA}$ is given by:
\Eq{\hat{J}(y)\hat{J}=R_*(y^*),}
where $R_*$ is an anti-homomorphism given by
\Eq{R_*(\hat{A})=A_*,\tab R_*(\hat{B})=B_*,}
where \Eq{A_* =e^{-2\pi bt}, \tab B_*=e^{-2\pi bp_t}.}
\end{Prop}

\begin{Cor}\label{Jcom} $J\circ \cA\circ J$ commutes with $\cA$  as operators on $\cH_{rep}$, and similarly $\hat{J}\circ \hat{\cA}\circ \hat{J}$ commutes with $\hat{\cA}$ as operators on  $\hat{\cH}_{rep}$, which verifies Theorem \ref{MvN}.
\end{Cor}
\subsection{The corepresentation space}\label{sec:Trans:corep}
Recall that $W$ satisfies the pentagon equation and the coaction axiom \eqref{DW}. In a more familiar form we see that $W'=W_{21}:=\S W \S \in M(\hat{\cA}\ox \cA)$ satisfies
\Eq{(1\ox \D)W' = W_{12}'W_{13}'.}
Hence if we treat $W'$ as an element in $M(\cB(\cH)\ox \cA)$ instead, we obtain the \emph{left regular corepresentation}
\begin{eqnarray} \Pi: \cH &\to& \cH \ox M(\cA)\nonumber\\
f&\mapsto& W'(f\ox 1)
\end{eqnarray}
which satisfies 
\Eq{ (1\ox \D)\circ \Pi = (\Pi \ox 1)\circ \Pi.}

This corepresentation picture is useful when we study the corepresentation of the quantum double. Therefore we introduce the transformation from $\cH$ to $\cH_{corep}$ so that $\hat{\cA}$ acts canonically. Similarly, by considering $\hat{W}$, we also describe the transformation so that $\cA$ acts canonically on $\hat{\cH}$.
\begin{Def} We define the transformation $\cT_{co}:\cH\to\cH_{corep}$ by
\begin{eqnarray} 
\cT_{co}: f(s,t)&\mapsto& f(-s,s-t+\frac{iQ}{2})G_b(\frac{Q}{2}-is+it)\inv,\\
\cT_{co}\inv: f(s,t)&\mapsto& f(-s,-t-s+\frac{iQ}{2})G_b(-it).
\end{eqnarray}

Similarly, we define the transformation $\hat{\cT}_{co}:\hat{\cH}\to\hat{\cH}_{corep}$ by
\begin{eqnarray}
\hat{\cT}_{co}:\hat{f}(s,t)&\mapsto& f(-s,s-t+\frac{iQ}{2})G_b(\frac{Q}{2}+is-it)e^{2\pi is(s-t+\frac{iQ}{2}),}\\
\hat{\cT}_{co}\inv:\hat{f}(s,t)&\mapsto&f(-s,-t-s+\frac{iQ}{2})e^{2\pi ist}G_b(Q+it)\inv.
\end{eqnarray}
\end{Def}

\begin{Prop} Under the transformation $\cT_{co}$, the action of $\hat{\cA}$ on $\cH_{corep}$ is given by
\Eq{\hat{A}=e^{-2\pi bs},\tab \hat{B}=e^{2\pi bp_s}.}
Under the transformation $\hat{\cT}_{co}$, the action of $\cA$ on $\hat{\cH}_{corep}$ is given by
\Eq{A=e^{2\pi bs},\tab B=e^{2\pi bp_s}.}
Furthermore, the $L^2$ measure on $\cH_{corep}$ and $\hat{\cH}_{corep}$ becomes the standard Lebesgue measure.
\end{Prop}
\begin{proof} It suffices to see that the transformations can be written as
\Eq{\cT_{co}=(t\mapsto t-s)\circ (-s,-t)\circ (t\mapsto t+\frac{iQ}{2}) \circ G_b(-it)\inv,}
\Eq{\hat{\cT}_{co}=(t\mapsto t-s)\circ (-s,-t)\circ (t\mapsto t+\frac{iQ}{2}) \circ G_b(Q+it)e^{-2\pi ist},}
followed by applying the transformation rules.
\end{proof}

This choice allows us to reproduce the representation from the pairing between $\cA$ and $\hat{\cA}$ given by the corepresentation associated to the multiplicative unitary.

\begin{Prop} The left regular corepresentation associated to $W'$ on $\cA$ is given by
\Eq{f(s,t)\mapsto \int_{\R+i0} \frac{f(s+\t,t)}{G_b(Q+i\t)}B^{ib\inv\t} A^{ib\inv s}d\t .}
\end{Prop}
\begin{Rem} Note that the left regular corepresentation is related to the left ``fundamental" representation
\Eq{\label{qplanecorep}\int f(s)x^{b\inv (-\frac{Q}{2}+is)}ds\mapsto \left(\int f(s)(xA+B)^{b\inv (-\frac{Q}{2}+is)}ds\right)\d^{-\frac{1}{2}},}
i.e.
\Eq{f(s)\mapsto \int f(s+\t)\frac{G_b(-i\t)G_b(\frac{Q}{2}-is)}{G_b(\frac{Q}{2}-i\t-is)}A^{ib\inv s-\frac{Q}{2b}}B^{ib\inv \t}A^{\frac{Q}{2b}}d\t}
$$=\int f(s+\t)\frac{G_b(\frac{Q}{2}+is+i\t)}{G_b(Q+i\t)G_b(\frac{Q}{2}+is)}B^{ib\inv \t}A^{ib\inv s}d\t$$
by multiplication by $G_b(\frac{Q}{2}+is)\inv$. Here $\d=A^{-\frac{Q}{b}}$ is the modular element.
\end{Rem}

\begin{Prop}\label{Bqaction} Under the pairing $\<1\ox \cA_q^*, W'(f\ox 1)\>$, the left regular representation of the quantum plane algebra $\cB_q=\cA_q^*$ is given by the canonical action
\Eq{X=e^{2\pi bs},\tab Y=e^{2\pi b(p_s+s)}.}
Similarly we can define the action of its modular double counterpart
\Eq{\til[X]=e^{2\pi b\inv s},\tab \til[Y]=e^{2\pi b\inv (p_s+s)},}
so that it extends to a representation of the modular double $\cB_{q\til[q]}$. Under a unitary transform by multiplication by $e^{\pi is^2}$ on $\cH_{corep}$, the action becomes the canonical action of $\cB_{q\til[q]}$.
\end{Prop}
\begin{proof}
It follows directly from Proposition \ref{lamda} since $\l$ is the identity map, hence the elements of $\hat{\cA}$ are sent to the corresponding actions. Now the formula follows from the definition $X=\hat{A}\inv, Y=q\hat{B}\hat{A}\inv$.
\end{proof}

Similarly the corepresentation associated to $\hat{W}'$ is given by
\Eq{f(s)\mapsto \int f(s+\t)G_b(-i\t)\hat{B}^{ib\inv \t}\hat{A}^{ib\inv s}d\t.}

For completeness, by composing $\cT_{co}$ with $\cT\inv$, we obtain the transformation $\cS:\cH_{rep}\to\cH_{corep}$
\begin{Prop}\label{Strans}The transformation $\cS:\cH_{rep}\to \cH_{corep}$ is given by
\begin{eqnarray}
\cS: f(s,t) &\mapsto&\int_\R f(\a-s, \a-t)G_b(\frac{Q}{2}-is+it)\inv e^{2\pi i(\a-s)s} d\a,\\
\cS\inv: f(s,t)&\mapsto& \int_\R f(\a-s,\a-t)G_b(\frac{Q}{2}+is-it) e^{-2\pi i(\a-s)s} d\a.
\end{eqnarray}

Similarly we define the corresponding transformation for $\cA$ using $\hat{\cS}:\hat{\cH}_{rep}\to\hat{\cH}_{corep}$ by
\begin{eqnarray}
\hat{\cS}: \hat{f}(s,t)&\mapsto& \int_\R f(\a-s,\a-t)G_b(\frac{Q}{2}+is-it)e^{-2\pi i(\a-s)s}d\a,\\
\hat{\cS}\inv: \hat{f}(s,t)&\mapsto& \int_\R f(\a-s,\a-t)G_b(\frac{Q}{2}-is+it)\inv e^{2\pi i(\a-s)s}d\a.
\end{eqnarray}
\end{Prop}


\section{The quantum double construction}\label{sec:double}
\subsection{Definitions}\label{sec:double:Def}
In this section we will describe the quantum double group construction given by \cite{PW} (see also \cite{Pu}) associated with the quantum plane, and show that the object we obtain is exactly the quantum ``semigroup" $GL_q^+(2,\R)$, also called the split quantum Minkowski spacetime, which is a generalization of the compact Minkowski spacetime introduced in \cite{F, FJ}.

\begin{Def} We define the split quantum Minkowski Spacetime $\cM_q^+(\R)$ as the Hopf *-algebra generated by positive self-adjoint operators $\z[ij], i,j\in\{1,2\}$ such that the following relations hold:

\begin{eqnarray*}
\label{com}
\left[\z[11],\z[12]\right]&=&0,\\
\left[\z[21],\z[22]\right]&=&0,\\
\left[\z[11],\z[22]\right]&=&[\z[12],\z[21]],\\
\z[11]\z[21]&=&q^{2}\z[21]\z[11],\\
\z[12]\z[22]&=&q^{2}\z[22]\z[12],\\
\z[12]\z[21]&=&q^{2}\z[21]\z[12],\\
\end{eqnarray*}
and the coproduct is given by
\begin{eqnarray}
\D(\z[11])&=&\z[11]\ox \z[11]+\z[12]\ox \z[21],\\
\D(\z[12])&=&\z[11]\ox \z[12]+\z[12]\ox \z[22],\\
\D(\z[21])&=&\z[21]\ox\z[11]+\z[22]\ox \z[21],\\
\D(\z[22])&=&\z[21]\ox\z[12]+\z[22]\ox \z[22].
\end{eqnarray}

It can also be realized as $GL_q^+(2,\R)$ in matrix form:
\Eq{\veca{\z[11]&\z[12]\\\z[21]&\z[22]}}
so that the coproduct is simply given by
\Eq{\D\left(\veca{\z[11]&\z[12]\\\z[21]&\z[22]}\right)=\veca{\z[11]&\z[12]\\\z[21]&\z[22]}\ox \veca{\z[11]&\z[12]\\\z[21]&\z[22]}.}

The quantum determinant $N$ is the positive self-adjoint operator defined by 
\Eq{N=\z[11]\z[22]-\z[12]\z[21]=\z[22]\z[11]-\z[21]\z[12],}
and we have
\Eq{N\z[11]=\z[11]N,\tab N\z[12]=q^{-2}\z[12]N,\tab N\z[21]=q^2\z[21]N,\tab N\z[22]=\z[22]N.}
\end{Def}

\begin{Prop} There is a projection map $\cP: GL_q^+(2,\R)\to SL_q^+(2,\R)$ given by
\Eq{\veca{\z[11]&\z[12]\\\z[21]&\z[22]}\mapsto \veca{a&b\\c&d}:=\veca{N^{-1/2}\z[11]&q^{-1/2}N^{-1/2}\z[12]\\q^{1/2}N^{-1/2}\z[21]&N^{-1/2}\z[22]},}
where $a,b,c,d$ satisfies the usual relations for $SL_q(2,\R)$:
$$ab=qba,\tab  ac=qca,\tab  ad=qda,\tab  bd=qdb,\tab  cd=qdc,$$
\Eq{bc=cb,\tab  ad-qbc = da-q\inv cb = 1.}
\end{Prop}

\begin{Prop}\label{gauss} There is a Gauss decomposition for $GL_q^+(2,\R)$ given by
\Eq{\veca{\z[11]&\z[12]\\\z[21]&\z[22]}=\veca{A&0\\B&1}\veca{1&\hat{B}\\0&\hat{A}}=\veca{A&A\hat{B}\\B&B\hat{B}+\hat{A}},}
where $A,B,\hat{A},\hat{B}$ are positive operators so that $\{A,B\}$ commutes with $\{\hat{A},\hat{B}\}$, with $AB=q^2BA$ and $\hat{A}\hat{B}=q^{-2}\hat{B}\hat{A}$. Furthermore we have $N=A\hat{A}$.
\end{Prop}

Now we will describe the quantum double group construction, and show that the result is precisely $GL_q^+(2,\R)$ together with the above Gauss decomposition.

\begin{Def} The quantum double group $\cD(\cA)$ is the Hopf algebra where as an algebra $\cD(\cA)\simeq\cA\ox \cA^{*op}=\cA\ox \hat{\cA}$ with the usual tensor product algebra structure, and with coproduct given by
\Eq{\D_m(x\ox \hat{x}) = (1\ox \s m \ox 1)(\D(x) \ox \hat{\D}(\hat{x})),}
where $\s$ is the permutation of the tensor product, and $m: M(\cA\ox \hat{\cA}) \to M(\cA\ox \hat{\cA})$ is called the matching, defined by
\Eq{m(x\ox\hat{x}) = W(x\ox \hat{x})W^*,}
with $W\in M(\cA\ox \hat{\cA})$ the multiplicative unitary defined in \eqref{W1}.
\end{Def}

Hence a general element in $\cD(\cA)$ can be written as 
\Eq{\iint\iint f(s_1,t_1)g(s_2,t_2)A^{ib\inv s_1}B^{ib\inv t_1}\hat{B}^{ib\inv t_2}\hat{A}^{ib\inv s_2} ds_1 dt_1 ds_2 dt_2}
or simply $f(s_1,t_1)g(s_2,t_2)$. For simplicity, we write $A:=A\ox1$, $\hat{A}:=1\ox\hat{A}$ and so on, and we will use this notation in the remaining of the paper.

\begin{Prop} \cite[Thm 4.1]{PW} $\D_m$ is coassociative, so it indeed defines a coproduct on $\cD(\cA)$. 
\end{Prop}
\begin{Prop} \cite[Thm 4.2]{PW}The Haar functional $h$ on $\cD(\cA)$ defined by 
\Eq{h = h_L\ox \hat{h}_R}
is both left and right invariant:
\Eq{(h\ox 1\ox 1)\D_m(x\ox \hat{x})=h(x\ox\hat{x})(1\ox 1),}
\Eq{(1\ox 1\ox h)\D_m(x\ox \hat{x})=h(x\ox\hat{x})(1\ox 1).}
In particular, the GNS representation of $\cD(\cA)$ is given by $\L_m:=\L\ox \hat{\L}_R$ on $\cH\ox\hat{\cH}$.
\end{Prop}
\begin{proof} Although the theorem in \cite{PW} applies only to compact quantum groups, the calculations using the graphical method there can be adapted in this setting without any changes, since it only depends on the invariances for $h_L, \hat{h}_R$ and the relations between the matching $m$ and the coproducts of $\cA$ and $\hat{\cA}$.
\end{proof}

\begin{Rem} The Haar functional on a general element $f(s_1,t_1)g(s_2,t_2)\in\cH\ox\hat{\cH}$ is thus given by 
\Eq{h(f\ox g) = f(0,iQ)g(-iQ,iQ)e^{-\pi iQ^2}.}
If we parametrize the element instead as
\Eq{f\ox g:=\iint\iint f(s_1,t_1)g(s_2,t_2)A^{ib\inv s_1}B^{ib\inv t_1}X^{ib\inv s_2}Y^{ib\inv t_2} ds_1 dt_1 ds_2 dt_2,}
where $X=\hat{A}\inv, Y = q\hat{B}\hat{A}\inv$, then it takes a more symmetric form
\Eq{h(f\ox g)=f(0,iQ)g(0,iQ),}
which means it only depends on the element $BY$. According to the Gauss decomposition, this is precisely
\Eq{\label{hyperbolic}B(q\hat{B}\hat{A}\inv)=qB\hat{B}AA\inv\hat{A}\inv=q\z[21]\z[12]N\inv,}
which corresponds under the projection to $SL_q^+(2,\R)$ the hyperbolic element $\ze:=bc$ that is crucial in the study of the $SU_q(2)$ and $SU_q(1,1)$ case in \cite{MMNNSU, MMNNU}.
\end{Rem}

\begin{Thm} By the Gauss decomposition, the Hopf algebra $GL_q^+(2,\R)$ can be naturally put into the $C^*$-algebraic setting, so that it is identified with the quantum double $\cD(\cA)$. Furthermore, the coproduct on $\cD(\cA)$ induces the same coproduct on the generators $z_{ij}$.
\end{Thm}
\begin{proof} By the Gauss decomposition, there is a one-to-one correspondence between the generators. Explicitly the inverse is given by
\begin{eqnarray*}
A&=&\z[11],\\
B&=&\z[21],\\
\hat{B}&=&\z[12]\z[11]\inv,\\
\hat{A}&=&N\z[11]\inv.
\end{eqnarray*}
We have to show that the coproduct is the same. The following calculations are very similar to the one given in \cite{Pu}.

\begin{Lem} We have the following commutation relations between $W$ and $\cD(\cA)$:
\begin{eqnarray*}
W(A\ox 1)W^* &=& A\ox 1+B\ox \hat{B},\\
W(A\ox\hat{A})W^* &=& A\ox \hat{A},\\
W(A\ox\hat{B})W^* &=& 1\ox \hat{B},\\
W(B\ox 1)W^* &=& B\ox \hat{A},\\
W^*(1\ox \hat{A})W &=& 1\ox\hat{A}+B\ox \hat{B}.\\
\end{eqnarray*}
\end{Lem}
\begin{proof} These follow from the summation properties \eqref{qsum1},\eqref{qsum2} for $g_b$, as well as the commutation relations \eqref{comW1},\eqref{comW2} for the exponentials.
\end{proof}

Now we proceed to the calculation of the coproduct on the generators:
\begin{eqnarray*}
\D(\z[11])&=&\D(A \ox 1)\\
&=&A \ox \s m(A\ox 1) \ox 1\\
&=&(A\ox 1 )\ox (A \ox 1)+(A\ox \hat{B})\ox( B \ox 1)\\
&=&\z[11]\ox \z[11]+\z[12]\ox \z[21],\\\\
\D(\z[12])&=&\D(A \ox \hat{B})\\
&=& A\ox \s m(A\ox 1)\ox \hat{B}+A \ox \s m(A\ox \hat{B}) \ox \hat{A}\\
&=&A\ox 1\ox A \ox \hat{B}+A\ox \hat{B}\ox B \ox \hat{B}+A\ox \hat{B}\ox 1\ox \hat{A}\\
&=&(A\ox 1)\ox (A \ox \hat{B})+(A\ox \hat{B})\ox (B \ox \hat{B}+1\ox \hat{A})\\
&=&\z[11]\ox\z[12]+\z[12]\ox\z[22],\\\\
\D(\z[21])&=&\D(B \ox 1)\\
&=& B\ox \s m(A\ox 1)\ox 1+1 \ox \s m(B\ox1) \ox 1\\
&=&B\ox 1\ox A \ox 1+B\ox \hat{B}\ox B\ox 1+1 \ox \hat{A}\ox B \ox 1\\
&=&(B\ox 1)\ox( A \ox 1)+(B\ox \hat{B}+1 \ox \hat{A})\ox (B \ox 1)\\
&=&\z[21]\ox\z[11]+\z[22]\ox\z[21],\\\\
\D(\z[22])&=&\D(B \ox \hat{B}+1\ox\hat{A})\\
&=&B\ox \s m(A\ox 1)\ox\hat{B}+1\ox \s m(B\ox 1) \ox \hat{B}+B\ox \s m(A\ox \hat{B})\ox \hat{A}+\\
&&1\ox \s m(B\ox \hat{B}+1\ox\hat{A})\ox\hat{A}.\\
\end{eqnarray*}
Since $W^*(1\ox\hat{A})W=1\ox \hat{A} + B\ox \hat{B}$, we have $$\s m(B\ox\hat{B}+1\ox \hat{A})=\hat{A}\ox 1.$$ Hence
\begin{eqnarray*}
\D(\z[22])&=&B\ox1\ox A\ox\hat{B}+B\ox\hat{B}\ox B\ox\hat{B}+1\ox\hat{A}\ox B\ox\hat{B}+\\
&&B\ox\hat{B}\ox1\ox\hat{A}+1\ox \hat{A}\ox 1\ox\hat{A}\\
&=&(B\ox 1)\ox (A\ox \hat{B}) + (B\ox\hat{B}+1\ox\hat{A})\ox(B\ox\hat{B}+1\ox\hat{A})\\
&=&\z[21]\ox\z[12]+\z[22]\ox\z[22].\\
\end{eqnarray*}

Finally, let us also derive the coproduct for the determinant $N$:
\begin{eqnarray*}
\D(N)&=&\D(A\ox \hat{A})\\
&=&A\ox \s m(A \ox \hat{A})\ox \hat{A}\\
&=&A\ox \hat{A}\ox A \ox \hat{A}\\
&=& N\ox N.
\end{eqnarray*}
\end{proof}
\subsection{The matrix coefficients and the fundamental corepresentation}\label{sec:double:matrix}
Recall that for $SL^+(2,\R)\sub SL(2,\R)$ the positive semigroup, the class of principal series representation can be expressed by considering the actions on homogeneous monomials (see e.g. \cite{Vi}):
\begin{eqnarray}
x^{l-i\mu}y^{l+i\mu} &\mapsto& (ax+cy)^{l-i\mu}(bx+dy)^{l+i\mu}\nonumber\\
&:=&\int_\R t_{\mu\nu}^l(g)x^{l-i\nu}y^{l+i\nu}d\nu,
\end{eqnarray}
where $\mu\in\R, l\in-\frac{1}{2}+i\R$, so that the representation is unitary.
In terms of coordinates, the representation acts on $L^2(\R)$ by
\Eq{g\cdot f(\mu)\mapsto \int_\R t_{\mu\nu}^l(g) f(\nu)d\nu.}

In \cite{F}, it is noted that for commuting variables $x,y$,
\Eq{[x\z[11]+y\z[21], x\z[12]+y\z[22]] = 0} 
by the commutation relations. Hence the following fundamental corepresentation for $\cD(\cA)=GL_q^+(2,\R)$ is well-defined:
\begin{Def}\label{fundT}
The fundamental corepresentation of $GL_q^+(2,\R)$ on $L^2(\R)$ is defined by
\Eq{T^\l: f(s)\mapsto \int_\R T_{s,\a}^\l(z) f(\a)d\a,}
where $T_{s,\a}^\l(z)$ is formally defined by
\Eq{(x\z[11]+y\z[21])^{b\inv(l-is)}(x\z[12]+y\z[22])^{b\inv(l+is)}N^{\frac{Q}{2b}}=:\left(\int_\R T_{s,\a}^\l(z) x^{b\inv(l-i\a)}y^{b\inv(l+i\a)}d\a\right)}
with $l:=-\frac{Q}{2}+i\l, \a,\l\in\R$ and $N=\z[11]\z[22]-\z[12]\z[21]=A\hat{A}$. 

More generally, we can introduce arbitrary character of the determinant, and will consider the corepresentation $T^{\l,t}$ defined by
\Eq{T_{s,\a}^{\l,t}(z):= N^{ib\inv(t-\l)/2}T_{s,\a}^\l(z) N^{ib\inv(t-\l)/2}.}
\end{Def}
\begin{Rem}The term $N^{\frac{Q}{2b}}$ is to make the corepresentation unitary. As seen from below, it is coming from the modular element $\d^{-\frac{1}{2}}$ in the fundamental corepresentation \eqref{qplanecorep} of the quantum plane. The factor chosen for $N$ in $T^{\l,t}$ is for later convenience when we obtain the fundamental representation.  
\end{Rem}
\begin{Rem}Using the Mellin transform, we can also write $T_{s,\a}^\l(z)$ as
\Eq{\label{mellinmatrix}T_{s,\a}^\l(z)=\frac{1}{2\pi b}\left(\int_0^\oo (x\z[11]+\z[21])^{b\inv(l-i\a)}(x\z[12]+\z[22])^{b\inv(l+i\a)} x^{b\inv(-l+is)} \frac{dx}{x}\right)N^{\frac{Q}{2b}}.}
This can be seen as the generalization of the matrix coefficients in the compact case, see e.g. \cite{FL}.
\end{Rem}
\begin{Prop} The matrix coefficient is given explicitly by
\begin{eqnarray}
T_{s,\a}^{\l,t}(z)&=&\int_{\R+i0}\veca{l+i\a\\i\t+i\a}_b\veca{l+i\t\\is+i\t}_be^{\pi i(t-\l)(s+\a+2\t)}e^{\pi Q(s+\t)}\cdot \nonumber\\
&&A^{ib\inv(t-s)}B^{ib\inv(s+\t)}\hat{B}^{ib\inv(\a+\t)}\hat{A}^{ib\inv(t-\t)}d\t\nonumber\\\\
&=&\int_{\R+i0} \frac{G_b(-i\t-i\a)G_b(\frac{Q}{2}+i\t-i\l)}{G_b(\frac{Q}{2}-i\a-i\l)}\frac{G_b(-i\t-is)G_b(\frac{Q}{2}+is-i\l)}{G_b(\frac{Q}{2}-i\t-i\l)}\cdot\nonumber\\
&&e^{\pi i(t-\l)(s+\a+2\t)}e^{\pi Q(s+\t)} A^{ib\inv(t-s)}B^{ib\inv(s+\t)}\hat{B}^{ib\inv(\a+\t)}\hat{A}^{ib\inv(t-\t)}d\t,\nonumber\\
\end{eqnarray}
where $l=-\frac{Q}{2}+i\l$ and $\veca{\a\\\b}_b=\frac{G_b(-\b)G_b(\b-\a)}{G_b(-\a)}$ is the $q$-binomial coefficient (cf. Lemma \ref{qbi}).

Hence under some changes of variables, the corepresentation $$f(s)\mapsto \int_\R T_{s,\a}^{\l,t}(z) f(\a)d\a$$ is given by
\begin{eqnarray}
f(s)&\mapsto& \int_{\R+i0}\int_{\R+i0}f(s+\a-\t) \frac{G_b(-i\a)}{G_b(Q+i\t)}\frac{G_b(\frac{Q}{2}-is+i\l+i\t)G_b(\frac{Q}{2}-is-i\l+i\t)}{G_b(\frac{Q}{2}-is+i\l)G_b(\frac{Q}{2}-is-i\l-i\a+i\t)}\nonumber\\
&&e^{\pi i\l(\t-\a)}e^{-2\pi i\t s}e^{\pi it(\t+\a)} A^{ib\inv(t-s)}B^{ib\inv \t}\hat{B}^{ib\inv\a}\hat{A}^{ib\inv(t+s-\t)}d\a d\t.\nonumber\\\label{fundcorep}
\end{eqnarray}
\end{Prop}
\begin{proof} We can make use of the Gauss decomposition (Proposition \ref{gauss}) and obtain
\Eq{\label{matrixdecomp}T_{s,\a}^{\l}(z)= \int t_{s,\t}^\l(\cA) \hat{t}_{\t,\a}^\l(\hat{\cA}) d\t,}
where the matrix coefficients correspond to the fundamental representation of the quantum plane \eqref{qplanecorep}:
\begin{eqnarray*}
t_{s,\a}^\l&:&\int f(s)x^{b\inv(l-is)}y^{b\inv(l+is)}ds\\
&\mapsto& \left(\int f(s)(xA+yB)^{b\inv(l-is)}y^{b\inv(l+is)}ds\right)A^{\frac{Q}{2b}}\\
&=&\int f(s)\veca{l-is\\i\a}_bA^{-\frac{Q}{2b}+ib\inv(\l-s-\a)}B^{ib\inv \a}A^{\frac{Q}{2b}} x^{b\inv(l-is-i\a)}y^{b\inv(l+is+i\a)} d\a ds\\
&=&\int f(s)\veca{l-is\\i\a-is}_be^{\pi Q(\a-s)}A^{ib\inv(\l-\a)}B^{ib\inv(\a-s)}x^{b\inv(l-i\a)}y^{b\inv(l+i\a)}d\a ds.\\
\end{eqnarray*}
Hence renaming $s\corr \a$, we see that the matrix coefficient is given by
\Eq{t_{s,\a}^\l(\cA)=\veca{l-i\a\\is-i\a}_be^{\pi Q(s-\a)}A^{ib\inv(\l-s)}B^{ib\inv(s-\a)}.}
A similar analysis using
$$\hat{t}_{s,\a}^\l:\int f(s)x^{b\inv(l-is)}y^{b\inv(l+is)}ds \mapsto\left(\int f(s)x^{b\inv(l-is)}(x\hat{B}+y\hat{A})^{b\inv(l+is)}ds\right)\hat{A}^{\frac{Q}{2b}}$$
shows that $\hat{t}_{s,\a}^\l(\hat{\cA})$ is given by
\Eq{\hat{t}_{s,\a}^\l(\hat{\cA})=\veca{-l+i\a\\-is+i\a}_b\hat{B}^{ib\inv(\a-s)}\hat{A}^{ib\inv(s+\l)}.}
Hence using \eqref{matrixdecomp} with the contour of $\t$ separating the poles of $G_b(\cdot \pm i\t)$, we obtain
\begin{eqnarray*}
&&T_{s,\a}^{\l}(z)\\
&=&\int_{\R-i0} \veca{l+i\a\\-i\t+i\a}_b\veca{l-i\t\\is-i\t}_be^{\pi Q(s-\t)} A^{ib\inv(\l-s)}B^{ib\inv(s-\t)}\hat{B}^{ib\inv(\a-\t)}\hat{A}^{ib\inv(\l+\t)}d\t\\
&=&\int_{\R+i0}\veca{l+i\a\\i\t+i\a}_b\veca{l+i\t\\is+i\t}_be^{\pi Q(s+\t)} A^{ib\inv(\l-s)}B^{ib\inv(s+\t)}\hat{B}^{ib\inv(\a+\t)}\hat{A}^{ib\inv(\l-\t)}d\t.
\end{eqnarray*}
Finally, for $T_{s,\a}^{\l,t}(z)$, by commuting $A$ and $\hat{A}$ of $N^{ib\inv(t-\l)/2} = (A\hat{A})^{ib\inv(t-\l)/2}$ to the corresponding sides, we pick up the factor $e^{\pi i(t-\l)(s+\a+2\t)}$, and obtain the desired formula.
\end{proof}
\begin{Cor}\label{matrixFb} The matrix coefficient $T_{s,\a}^\l(z)$ can also be expressed in closed form as
\begin{eqnarray}T_{s,\a}^\l(z)&=&A^{b\inv(l-is)}B^{b\inv(l+is)}\hat{B}^{b\inv(l+i\a)}\veca{2l\\l+is}_bF_b(-l-is,-l-i\a,-2l; -\ze\inv)N^{\frac{Q}{2b}},\nonumber\\\label{Tsa}
\end{eqnarray}
where $l=-\frac{Q}{2}+i\l$, $\ze=qB\hat{B}\hat{A}\inv$ is the hyperbolic element defined in \eqref{hyperbolic},
\Eq{F_b(\a,\b,\c; z):=\frac{G_b(\c)}{G_b(\a)G_b(\b)}\int_C (-z)^{ib\inv \t}e^{\pi i\t^2}\frac{G_b(\a+i\t)G_b(\b+i\t)G_b(-i\t)}{G_b(\c+i\t)}d\t} is the $b$-hypergeometric function (slightly modified from \cite{PT2}), which is the quantum analogue of the classical ${}_2F_1(a,b,c;z)$. Hence \eqref{Tsa} is exactly the quantum analogue of the classical formula for $SL^+(2,\R)$ given in \cite[VII. 4.1(4)]{Vi}, where the hyperbolic element is $\ze=\sinh^2\h$:
\Eq{T_{mn}^l=\frac{1}{2\pi i}A^{l-m}B^{l+m}\hat{B}^{l+n}\frac{\G(-l-m)\G(-l+m)}{\G(-2l)}{}_2F_1(-l-m,-l-n,-2l;-\frac{1}{\sinh^2\h}),}
where $A,B,\hat{B}$ are the corresponding variables in the Gauss decomposition for $SL_2^+(\R)$.
\end{Cor}
\subsection{The multiplicative unitary}\label{sec:double:multi}
Given two locally compact quantum groups $(M_1,\D_1)$ and $(M_2,\D_2)$, with a matching $m:M_1\ox M_2\to M_1\ox M_2$, the multiplicative unitary for the double crossed product is first constructed in \cite{BV}. The quantum double group construction is a special case given by $(M_1,\D_1)=(\cA,\D^{op})$ and $(M_2,\D_2)=(\hat{\cA},\hat{\D})$, while the matching $m$ is given by $m(x)=W(x)W^*$ as before \cite{K}. Let us restrict to the quantum plane case and describe the main ingredients needed.

\begin{Prop}
The multiplicative unitary operator is defined by
\Eq{\bW_m:=\bW=W_{13}Z_{34}^* \hat{W}_{24}Z_{34}\in\cB(\cH\ox\hat{\cH}\ox \cH\ox \hat{\cH}),}
where $Z\in\cB(\cH\ox \hat{\cH})$ is given by
\Eq{\label{defZ}Z=W(\hat{J_1}J_1\ox J_2\hat{J_2})W^*(\hat{J_1}J_1\ox \hat{J_2}J_2)=W(\hat{J}J\ox \hat{J}J)W^*(\hat{J}J\ox J\hat{J}),} and $\hat{W}=\S W^*\S=W_{21}^*$.

The coproduct is given by
\Eq{\D_m=(\iota\ox \s m\ox \iota)(\D_1^{op}\ox \D_2)=(\iota\ox \s m\ox \iota)(\D\ox \hat{\D}),}
where $m(z)=ZzZ^*$, and
\Eq{\D_m(z) = \bW^*(1\ox z)\bW}
for $z\in \cD(\cA)$.

The matching $m$ satisfies
\Eq{\label{mprop1}(\D^{op}\ox 1)m = m_{23}m_{13}(\D^{op}\ox 1),}
\Eq{\label{mprop2}(1\ox \hat{\D})m = m_{13}m_{12}(1\ox \hat{\D}).}
\end{Prop}

First of all we note the difference in the definition of the matching $m$. In fact they are the same.
\begin{Prop} We actually have
\Eq{m(x)=ZxZ^* =WxW^*}
for $x\in\cA\ox\hat{\cA}$. In particular from the pentagon equation of $W$ and $\hat{W}=W_{21}^*$, the matching satisfies the property \eqref{mprop1} and \eqref{mprop2}.
\end{Prop}
\begin{proof} Recall the conjugation properties from \eqref{R1}, \eqref{R2}:
\Eq{\label{conj1}(\hat{J}\ox J)W(\hat{J}\ox J)=W^*,}
\Eq{\label{conj2}(J\ox\hat{J})\hat{W}(J\ox\hat{J})=\hat{W}^*.}

Using also \eqref{nu}, we have
\begin{eqnarray*}
&&W(\hat{J}J\ox \hat{J}J)W^*(\hat{J}J\ox J\hat{J})\\
&=&\nu^{i/4}W(J\ox\hat{J})(\hat{J}\ox J)W^*(\hat{J}\ox J)(J\ox\hat{J})\\
&=&\nu^{i/4}W(J\ox\hat{J})W(J\ox\hat{J})
\end{eqnarray*} 

Hence using the definition \eqref{defZ} of $Z$, expanding $ZxZ^* =WxW^*$ we have

$$W(J\ox \hat{J})W(J\ox \hat{J}) x (J\ox \hat{J}) W^*(J\ox \hat{J}) W^* = WxW^*,$$
or
$$(J\ox \hat{J})W(J\ox \hat{J}) x = x(J\ox \hat{J}) W(J\ox \hat{J}),$$
that is, $(J\ox \hat{J}) W(J\ox \hat{J})$ commutes with every $x\in \cA\ox \hat{\cA}\sub \cB(\cH\ox \hat{\cH})$ as operators.
However, since $W\in \cA\ox\hat{\cA}$, from Corollary \ref{Jcom}, it is clear that
$$(J\ox \hat{J}) W(J\ox \hat{J})=g_b^*(q\hat{B}_*\hat{A}_*\inv \ox B_*)e^{\frac{i}{2\pi b^2}\log \hat{A}_*\ox \log A_*}$$ commutes with every $x\in \cA\ox\hat{\cA}$.
\end{proof}

For completeness, let us also reprove the coproduct formula.
\begin{proof}
Recall
$$\bW=W_{13}Z_{34}^* W_{42}^*Z_{34},$$
$$\D(x)=W^*(1\ox x)W,$$
$$\hat{\D}(y)=\hat{W}^*(1\ox y)\hat{W}=W_{21}(1\ox y)W_{21}^*.$$
Hence
\begin{eqnarray*}
\bW^*(1\ox z)\bW&=&Z_{34}^*W_{42}Z_{34}W_{13}^*(1\ox z)W_{13}Z_{34}^* W_{42}^*Z_{34}\\
&=&Z_{34}^*W_{42}Z_{34}(\D\ox \iota)(z)_{134}Z_{34}^* W_{42}^*Z_{34}\\
&=&Z_{34}^*W_{42}(\iota\ox m)((\D\ox \iota)(z))_{134}W_{42}^*Z_{34}\\
&=&Z_{34}^*((\iota\ox\iota\ox\hat{\D})(\iota\ox m)(\D\ox \iota)(z))_{1324}Z_{34}\\
&=&Z_{34}^*(m_{24}m_{23}(\iota\ox\iota\ox\hat{\D})(\D\ox \iota)(z))_{1324}Z_{34}\\
&=&Z_{34}^*(Z_{24}m_{23}(\iota\ox\iota\ox\hat{\D})(\D\ox \iota)(z)Z_{24}^*)_{1324}Z_{34}\\
&=&Z_{34}^*Z_{34}m_{32}(\iota\ox\iota\ox\hat{\D})(\D\ox \iota)(z)Z_{34}^*Z_{34}\\
&=&(\s m)_{23}(\iota\ox\iota\ox\hat{\D})(\D\ox \iota)(z)\\
&=&\D_m(z).
\end{eqnarray*}
\end{proof}

It is proved in \cite{BV} for general matching $m$ that $\bW_m:=\bW$ is a multiplicative unitary. Here we present a direct proof using the fact that $m$ is given by the multiplicative unitary $W$.

\begin{Thm} $\bW$ is a multiplicative unitary, i.e. it satisfies
$$\bW_{3456}\bW_{1234}=\bW_{1234}\bW_{1256}\bW_{3456}.$$
\end{Thm}
\begin{proof}
It suffices to check that
$$\bW_{3456}Z_{34}^*\hat{W}_{24}Z_{34}\bW_{3456}^*=Z_{34}^*Z_{56}^*\hat{W}_{24}\hat{W}_{26}Z_{34}Z_{56}$$
and
$$\bW_{3456}W_{13}\bW_{3456}^*=W_{35}W_{13}W_{35}^*=W_{13}W_{15}.$$

The second equation follows directly from the definition and the pentagon equation for $W$.

Let us write $Z=VV'$ where $V\in \cB(\cH\ox\hat{\cH})$ and $V'=(J\ox \hat{J}) V (J\ox\hat{J})$ is the copy of $W$ based on different spaces. Note that $V, V'$ commute. Also $V'$ commutes with entries in $\cA\ox\hat{\cA}\sub \cB(\cH\ox\hat{\cH})$.

From the pentagon equation for $W$, but with its legs sitting on different spaces, we have the relations (cf. Corollary \ref{DWW}):
\Eq{V_{23}W_{12}=W_{12}V_{13}V_{23}\in \cB(\cH\ox \cH\ox \hat{\cH}),}
\Eq{V_{13}\hat{W}_{23}=\hat{W}_{23}V_{13}V_{12}\in\cB(\cH\ox \hat{\cH}\ox \hat{\cH}).}

Using also \eqref{conj1} and \eqref{conj2}, we derive the relations
\begin{eqnarray*}
V'_{56}W_{35}^*&=&(J_5\ox \hat{J}_6) V_{56} (J_5\ox\hat{J}_6)(\hat{J}_3\ox J_5)W_{35}(\hat{J}_3\ox J_5)\\
&=&(\hat{J}_3\ox J_5\ox \hat{J}_6)V_{56}W_{35}(\hat{J}_3\ox J_5\ox \hat{J}_6)\\
&=&(\hat{J}_3\ox J_5\ox \hat{J}_6)W_{35}V_{36}V_{56}(\hat{J}_3\ox J_5\ox \hat{J}_6)\\
&=&W_{35}^* (\hat{J}_3\ox \hat{J}_6)V_{36} (\hat{J}_3\ox \hat{J}_6)(J_5\ox \hat{J}_6)V_{56} (J_5\ox \hat{J}_6)\\
&=&W_{35}^* (\hat{J}_3\ox \hat{J}_6)V_{36} (\hat{J}_3\ox \hat{J}_6) V_{56}',\\
\end{eqnarray*}

\begin{eqnarray*}
\hat{W}_{46}^*(\hat{J}_3\ox \hat{J}_6)V_{36}(\hat{J}_3\ox \hat{J}_6)&=&(J_4\ox \hat{J}_6)\hat{W}_{46}(J_4\ox \hat{J}_6)(\hat{J}_3\ox\hat{J}_6)V_{36}(\hat{J}_3\ox \hat{J}_6)\\
&=&(\hat{J}_3\ox J_4\ox \hat{J}_6)\hat{W}_{46}V_{36}(\hat{J}_3\ox J_4\ox \hat{J}_6)\\
&=&(\hat{J}_3\ox J_4\ox \hat{J}_6)V_{36}\hat{W}_{46}V_{34}^*(\hat{J}_3\ox J_4\ox \hat{J}_6)\\
&=&(\hat{J}_3\ox\hat{J}_6)V_{36}(\hat{J}_3\ox\hat{J}_6)(J_4\ox \hat{J}_6)\hat{W}_{46}(J_4\ox \hat{J}_6)(\hat{J}_3\ox J_4)V_{34}^*(\hat{J}_3\ox J_4)\\
&=&(\hat{J}_3\ox \hat{J}_6)V_{36}(\hat{J}_3\ox \hat{J}_6) \hat{W}_{46}^* V_{34}.
\end{eqnarray*}

Hence
\begin{eqnarray*}
&&W_{3456}Z_{34}^*\hat{W}_{24}Z_{34}W_{3456}^*\\
&=&W_{35}Z_{56}^*\hat{W}_{46}Z_{56}Z_{34}^*\hat{W}_{24}Z_{34}Z_{56}^*\hat{W}_{46}^*Z_{56}W_{35}^*\\
&=&(W_{35}V_{56}^*V_{56}^{*'}\hat{W}_{46}V_{34}^*V_{34}^{*'})\hat{W}_{24}(V_{34}V_{34}'\hat{W}_{46}^*V_{56}V_{56}'W_{35}^*).\\
\end{eqnarray*}
Now the right hand bracket is
\begin{eqnarray*}
&&V_{34}V_{34}'\hat{W}_{46}^*V_{56}V_{56}'W_{35}^*\\
&=&V_{34}\hat{W}_{46}^*V_{56}V_{56}'W_{35}^*V_{34}'\\
&=&V_{36}^*\hat{W}_{46}^*V_{36}V_{56}V_{56}'W_{35}^*V_{34}'\\
&=&V_{36}^*\hat{W}_{46}^*V_{36}V_{56}W_{35}^*(\hat{J}_3\ox\hat{J}_6)V_{36}(\hat{J}_3\ox\hat{J}_6)V_{56}'V_{34}'\\
&=&V_{36}^*\hat{W}_{46}^*W_{35}^*V_{56}(\hat{J}_3\ox\hat{J}_6)V_{36}(\hat{J}_3\ox\hat{J}_6)V_{56}'V_{34}'\\
&=&V_{36}^*W_{35}^*\hat{W}_{46}^*(\hat{J}_3\ox\hat{J}_6)V_{36}(\hat{J}_3\ox\hat{J}_6)Z_{56}V_{34}'\\
&=&V_{36}^*W_{35}^*(\hat{J}_3\ox\hat{J}_6)V_{36}(\hat{J}_3\ox\hat{J}_6)\hat{W}_{46}^*V_{34}Z_{56}V_{34}'\\
&=&V_{36}^*W_{35}^*(\hat{J}_3\ox\hat{J}_6)V_{36}(\hat{J}_3\ox\hat{J}_6)\hat{W}_{46}^*Z_{56}Z_{34}.\\
\end{eqnarray*}

Combining, we see that the left most 3 terms will pass through $\hat{W}_{24}$ and cancel with the left hand inverses.
Hence we obtain
\begin{eqnarray*}
&&W_{3456}Z_{34}^*\hat{W}_{24}Z_{34}W_{3456}^*\\
&=&Z_{34}^*Z_{56}^*\hat{W}_{46}\hat{W}_{24}\hat{W}_{46}^*Z_{56}Z_{34}\\
&=&Z_{34}^*Z_{56}^*\hat{W}_{24}\hat{W}_{26}Z_{56}Z_{34}\\
\end{eqnarray*}
as desired.
\end{proof}

Note that $\bW$ is actually a composition of 6 dilogarithm functions. We can simplify the expression to just 4 $g_b$'s by applying the pentagon equation.

\begin{Prop} We can rewrite $\bW$ as
\Eq{\label{Wmsimple}\bW=W_{13}V_{32}'' \hat{W}_{24}V_{32}^*,}
where 
\begin{eqnarray*}
V'' &:=& (J\ox J)V(J\ox J)\\
&=& e^{\frac{i}{2\pi b^2}\log\hat{A}_*\ox\log\hat{A}\inv}g_b^*(\hat{B}_*\ox \hat{B})\\
&=&e^{\frac{i}{2\pi b^2}\log\hat{A}_*\ox\log\hat{A}\inv}\int_{\R+i0}\hat{B}_*^{ib\inv \t}\ox \hat{B}^{ib\inv \t} e^{\pi i\t^2}G_b(-i\t)d\t.
\end{eqnarray*}
\end{Prop}
\begin{proof}
\begin{eqnarray*}
\bW&=&W_{13}Z_{34}^*\hat{W}_{24}Z_{34}\\
&=&W_{13} V_{34}^{*'}V_{34}^{*} \hat{W}_{24} V_{34} V_{34}'\\
&=&W_{13}V_{34}^{*'}\hat{W}_{24}V_{32}^* V_{34}'\\
&=&W_{13}V_{34}^{*'}\hat{W}_{24} V_{34}'V_{32}^*\\
&=&W_{13}(J_2\ox J_3\ox \hat{J}_4) V_{34}^*\hat{W}_{24}^* V_{34}(J_2\ox J_3\ox \hat{J}_4)V_{32}^*\\
&=&W_{13}(J_2 \ox J_3) V_{32} (J_2\ox J_3) \hat{W}_{24}V_{32}^*.\\
\end{eqnarray*}
The formula for $V''$ follows from the action of $Ad J$ given in Proposition \ref{RJ} on $\hat{\cH}$ and Proposition \ref{AdJ} on $\cH$.
\end{proof}

Since $\bW$ is a unitary operator on $\cH\ox \hat{\cH}\ox \cH\ox\hat{\cH}$, for completeness let us describe its action.

\begin{Prop} The action of $W$ on $\cH_{rep}\ox\cH_{rep}$ is given by
$$W=e^{\frac{i}{2\pi b^2}\log A\inv\ox \log \hat{A}}g_b(B\ox q\hat{B}\hat{A}\inv),$$

\Eq{W:F\ox G\mapsto \int_{\R+i0} F(s_1+\t,t_1)G(s_2-s_1-\t,t_2-s_1)\frac{e^{2\pi i\t(s_2-s_1-\t)}G_b(\frac{Q}{2}+is_2-it_2)}{G_b(\frac{Q}{2}+is_2-it_2-i\t)G_b(Q+i\t)}d\t,}
\Eq{W^*:F\ox G\mapsto \int_{\R+i0} F(s_1+\t,t_1)G(s_2+s_1,t_2+s_1+\t)\frac{e^{2\pi is_2\t}G_b(\frac{Q}{2}+is_2-it_2)G_b(-i\t)}{G_b(\frac{Q}{2}+is_2-it_2-i\t)}d\t.}
The action of $\hat{W}$ on $\hat{\cH}\ox\hat{\cH}_{rep}$ is given by
$$\hat{W}=e^{\frac{i}{2\pi b^2}\log \hat{A}\ox \log A}g_b^*(\hat{B}\ox q\inv BA\inv),$$
\begin{eqnarray}
\hat{W}_{24}:F\ox G&\mapsto& \int_{\R+i0} F(s_1+\t,t_1)G(s_2-s_1-\t,t_2-s_1)e^{2\pi i\t(\t+s_1-s_2)}\nonumber\\
&&\frac{G_b(\frac{Q}{2}-is_2+it_2+i\t)G_b(-i\t)}{G_b(\frac{Q}{2}-is_2+it_2)}d\t,
\end{eqnarray}
\Eq{\hat{W}_{24}^*:F\ox G\mapsto \int_{\R+i0} F(s_1+\t,t_1)G(s_2+s_1,t_2+s_1+\t)\frac{e^{-2\pi i\t s_2}G_b(\frac{Q}{2}-is_2+it_2+i\t)}{G_b(\frac{Q}{2}-is_2+it_2)G_b(Q+i\t)}d\t.}
The action of $V$ on $\cH_{rep}\ox\hat{\cH}_{rep}$ is given by
\Eq{\label{Z1}V:F\ox G\mapsto \int F(s_1+\t,t_1)G(s_2+\t,t_2)\frac{e^{2\pi is_2(s_1+\t)}}{G_b(Q+i\t)}d\t,}
\Eq{\label{Z*1}V^*:F\ox G\mapsto \int F(s_1+\t,t_1)G(s_2+\t,t_2)e^{-2\pi is_1(s_2+\t)}G_b(-i\t)d\t,}
and the action of $V''$ on $\cH_{rep}\ox\hat{\cH}_{rep}$ is given by
\Eq{\label{Z2}V'': F\ox G\mapsto \int_{\R+i0} F(s_1,t_1-\t)G(s_2+\t,t_2)e^{2\pi it_1s_2}e^{\pi i\t^2}G_b(-i\t) d\t,}
\Eq{\label{Z*2}V''^*: F\ox G\mapsto \int_{\R+i0} F(s_1,t_1-\t)G(s_2+\t,t_2)\frac{e^{2\pi i\t(s_2-t_1)}e^{-2\pi it_1s_2}e^{\pi i\t^2}}{G_b(Q+i\t)} d\t.}
Hence the action of $\bW$ and $\bW^*$ are the compositions of the corresponding operators.
\end{Prop}
\begin{proof} The actions follow directly from the closed form expression (Proposition \ref{W}) for $W$, using Corollary \ref{gbint} for the integral representations for $g_b$, as well as the action described in Section \ref{sec:Trans:rep}.
\end{proof}
\subsection{The left regular corepresentation}\label{sec:double:left}
With the construction of the multiplicative unitary $\bW$, we can talk about the corepresentation induced by it. By the unitary transformations given in Section \ref{sec:Trans:corep}, we can choose the Hilbert space to be $\cH_{corep}\ox\hat{\cH}_{corep}$ so that the action can be nicely described. 

Recall that the left regular corepresentation is given by 
\Eq{f\ox g \mapsto \bW'(f\ox g\ox 1\ox 1),}
where $f\ox g\in \cH_{corep}\ox\hat{\cH}_{corep}$ and $\bW' = \bW_{3412}\in \cB(\cH_{corep}\ox\hat{\cH}_{corep})\ox \cA\ox \hat{\cA}$. Here we realize $\cA \ox \hat{\cA}$ as operators in the 3rd and 4th components.

\begin{Thm} \label{leftco} The left regular corepresentation is given by
\begin{eqnarray*}&&f(s_1,t_1)g(s_2,t_2)\mapsto\\
&=&\int_{\R+i0} \int_\R \int_C f(s_1-\a+\t,t_1-\s)g(s_2+\s,t_2)\\
&&\frac{G_b(-i\s)G_b(-i\a+i\s)G_b(\frac{Q}{2}+is_1-it_1+i\t)G_b(\frac{Q}{2}+is_1-is_2-i\s+i\t)}{G_b(Q+i\t)G_b(\frac{Q}{2}+is_1-it_1-i\a+i\s+i\t)G_b(\frac{Q}{2}+is_1-is_2-i\a+i\t)}\\
&&e^{2\pi i\s(t_1-s_1-\t)}e^{2\pi is_1\t}A^{ib\inv s_1}B^{ib\inv \t}\ox\hat{B}^{ib\inv\a}\hat{A}^{ib\inv (t_1+s_2-s_1-\t)} d\s d\a d\t,\\
\end{eqnarray*}
where the contour $C$ of $\s$ goes above $\s=0$ and below $\s=\a$.
\end{Thm}

\begin{proof} Since we have used a permutation, now $\bW'$ reads
\Eq{\bW' = \bW_{3412} = W_{31}V''_{14}\hat{W}_{42}V_{14}^*.}

On the space $\cH_{corep}\ox \hat{\cH}_{corep}$, the coaction of the components of $\bW$ is given by:
\begin{eqnarray*}
W_{31}: f(s_1,t_1)g(s_2,t_2)&\mapsto& \int_{\R+i0} f(s_1+\t,t_1)g(s_2,t_2)\frac{e^{2\pi is_1\t}}{G_b(Q+i\t)}A^{ib\inv s}B^{ib\inv \t} d\t,\\
\hat{W}_{42}: f(s_1,t_1)g(s_2,t_2)&\mapsto& \int_{\R+i0} f(s_1,t_1)g(s_2+\s,t_2)G_b(-i\s)\hat{B}^{ib\inv\s}\hat{A}^{ib\inv s_2}d\s.\\
\end{eqnarray*}
The actions for $V$ and $V''$ are harder to describe since they involve formally $A^{ib\inv p}$. Therefore we use $\cS$ (cf. Proposition \ref{Strans}) to sent the first component back to $\cH_{rep}$ and obtain their actions.
\begin{eqnarray*}
&&\cS V_{14}^* \cS\inv:f\ox g \\
&\mapsto &\cS\cdot \left(\int_{\R+i0} B^{ib\inv \t}\ox \hat{B}^{ib\inv\t}\hat{A}^{-ib\inv \t}   G_b(-i\t)d\t \right) e^{\frac{i}{2\pi b^2} \log \hat{A}\ox\log A}\cdot\\
&& \int_\R f(\a-s_1,\a-t_1)g(s_2,t_2)G_b(\frac{Q}{2}-it_1+is_1) e^{-2\pi i(\a-s_1) s_1}d\a\\
&=&\cS\cdot \int_{\R+i0}\int_\R f(\a-s_1-\t,\a-t_1)g(s_2,t_2)G_b(\frac{Q}{2}-it_1+is_1+i\t)G_b(-i\t)  \cdot\\
&&e^{-2\pi i(\a-s_1-\t) (s_1+\t)} \hat{B}^{ib\inv\t}\hat{A}^{ib\inv s_1}d\a d\t\\
&=&\int_\R\int_{\R+i0}\int_\R f(\a-\b+s_1-\t,\a-\b+t_1)g(s_2,t_2)\frac{G_b(\frac{Q}{2}+it_1-is_1+i\t)G_b(-i\t)}{G_b(\frac{Q}{2}+it_1-is_1)}  \cdot\\
&&e^{-2\pi i(\a-\b+s_1-\t) (\b-s_1+\t)}e^{2\pi i(\b-s_1)s_1} \hat{B}^{ib\inv\t}\hat{A}^{ib\inv (\b-s_1)}d\a d\t d\b\\
&& \mbox{shifting $\a\mapsto \a+\b+\t, \b\mapsto \b+s_1$:}\\
&=&\int_\R\int_{\R+i0}\int_\R f(s_1+\a,t_1+\a+\t)g(s_2,t_2)\frac{G_b(\frac{Q}{2}+it_1-is_1+i\t)G_b(-i\t) }{G_b(\frac{Q}{2}+it_1-is_1)}\cdot\\
&&e^{-2\pi i(\a\b+\a\t+s_1\t)} \hat{B}^{ib\inv\t}\hat{A}^{ib\inv \b} d\a d\t d\b,\\
\end{eqnarray*}
and finally for $V''$:
\begin{eqnarray*}
&&\cS V_{14}'' \cS\inv:f\ox g \\
&\mapsto&\cS\cdot \left(e^{\frac{i}{2\pi b^2}\log\hat{A}_*\ox\log\hat{A}\inv}\int_{\R+i0}\hat{B}_*^{ib\inv \t}\ox \hat{B}^{ib\inv \t} e^{\pi i\t^2}G_b(-i\t)d\t\right)\cdot\\
&&\int_\R f(\a-s_1,\a-t_1)g(s_2,t_2)G_b(\frac{Q}{2}-it_1+is_1) e^{-2\pi i(\a-s_1) s_1}d\a\\
&=&\cS\cdot \int_{\R+i0}\int_\R f(\a-s_1\a-t_1+\t)g(s_2,t_2)G_b(\frac{Q}{2}-it_1+is_1+i\t) \\
&&e^{-2\pi i(\a-s_1) s_1}e^{\pi i\t^2}G_b(-i\t) \hat{A}^{-ib\inv t_1}\hat{B}^{ib\inv\t}d\a d\t\\
&=&\int_\R\int_{\R+i0}\int_\R f(\a-\b+s_1,\a-\b+t_1+\t)g(s_2,t_2)\frac{G_b(\frac{Q}{2}+it_1-is_1+i\t)G_b(-i\t)}{G_b(\frac{Q}{2}+it_1-is_1)} \\
&&e^{-2\pi i(\a-\b+s_1)(\b-s_1)}e^{\pi i\t^2}e^{2\pi i(\b-s_1)s_1}e^{2\pi i\t(t_1-\b)}\hat{B}^{ib\inv\t}\hat{A}^{ib\inv (t_1-\b)}d\a d\t d\b 
\end{eqnarray*}
\begin{eqnarray*}
&& \mbox{shifting $\a\mapsto \a+\b, \b\mapsto t_1-\b$:}\\
&=&\int_\R\int_{\R+i0}\int_\R f(s_1+\a,t_1+\a+\t)g(s_2,t_2)\frac{G_b(\frac{Q}{2}+it_1-is_1+i\t)G_b(-i\t)}{G_b(\frac{Q}{2}+it_1-is_1)} \\
&&e^{-2\pi i\a(t_1-s_1)}e^{2\pi i\b(\a+\t)}e^{\pi i\t^2}\hat{B}^{ib\inv\t}\hat{A}^{ib\inv \b} d\a d\t d\b.\\
\end{eqnarray*}
Combining, we have 
\begin{eqnarray*}
&&f\ox g\\
&\mapsto_{V^*}&\int_\R\int_{\R+i\e_\t}\int_\R d\a d\t d\b f(s_1+\a,t_1+\a+\t)g(s_2,t_2)\cdot\\
&&\frac{G_b(\frac{Q}{2}+it_1-is_1+i\t)G_b(-i\t) }{G_b(\frac{Q}{2}+it_1-is_1)}e^{-2\pi i(\a\b+\a\t+s_1\t)} \hat{B}^{ib\inv\t}\hat{A}^{ib\inv \b}\\
&\mapsto_{\hat{W}}&\int_{\R+i\e_\s}\int_\R\int_{\R+i\e_\t}\int_\R d\a d\t d\b d\s f(s_1+\a,t_1+\a+\t)g(s_2+\s,t_2)\cdot\\
&&\frac{G_b(\frac{Q}{2}+it_1-is_1+i\t)G_b(-i\t)G_b(-i\s) }{G_b(\frac{Q}{2}+it_1-is_1)}e^{-2\pi i(\a\b+\a\t+s_1\t)} e^{2\pi is_2\t}\hat{B}^{ib\inv(\s+\t)}\hat{A}^{ib\inv (\b+s_2)}.\\
\end{eqnarray*}
The integral for $\b$ and $\s$ are independent, hence we can interchange them. Furthermore from the decay properties for $f\ox g$ and the asymptotic properties for $G_b(-i\t)G_b(-i\s)$, the integral for $\a,\t$ and $\s$ is absolutely convergent, hence we can interchange the order so that $d\s$ goes to the inner most layer.

\begin{eqnarray*}
&\mapsto_{V''}&\int_\R\int_{\R+i\e_\t'}\int_\R\int_\R\int_{\R+i\e_\t}\int_\R \int_{\R+i\e_\s} d\s d\a d\t d\b  d\a' d\t' d\b'\\ &&f(s_1+\a+\a',t_1+\a+\a'+\t+\t')g(s_2+\s,t_2)\\
&&\frac{G_b(\frac{Q}{2}+it_1-is_1+i\t+i\t')G_b(-i\t)G_b(-i\s) }{G_b(\frac{Q}{2}+it_1-is_1+i\t')}\frac{G_b(\frac{Q}{2}+it_1-is_1+i\t')G_b(-i\t')}{G_b(\frac{Q}{2}+it_1-is_1)}\cdot\\
&&e^{-2\pi i(\a\b+\a\t+(s_1+\a')\t)} e^{2\pi is_2\t}e^{-2\pi i\a'(t_1-s_1)}e^{2\pi i\b'(\a'+\t')}e^{\pi i\t'^2}e^{2\pi i\b'(\s+\t)}\\
&&\hat{B}^{ib\inv(\s+\t+\t')}\hat{A}^{ib\inv (\b+\b'+s_2)}
\end{eqnarray*}
\begin{eqnarray*}
&&\mbox{shifting $\a\mapsto \a-\a', \b\mapsto \b-\b', \t\mapsto \t-\t'$:}\\
&=&\int_\R\int_{\R+i\e_\t'}\int_\R \int_\R\int_{\R+i\e_\t-i\e_\t'}\int_\R \int_{\R+i\e_\s}d\s d\a d\t d\b d\a' d\t' d\b'\\ &&f(s_1+\a,t_1+\a+\t)g(s_2+\s,t_2)\frac{G_b(\frac{Q}{2}+it_1-is_1+i\t)G_b(i\t'-i\t)G_b(-i\s)G_b(-i\t')}{G_b(\frac{Q}{2}+it_1-is_1)}\\
&&e^{2\pi i\b'(\a+\s+\t)}e^{2\pi i\a'(\b+s_1-t_1)}e^{2\pi i(\a+s_1-s_2)(\t'-\t)-2\pi i\a\b+\pi i\t'^2}\hat{B}^{ib\inv(\s+\t)}\hat{A}^{ib\inv (\b+s_2)}.\\
\end{eqnarray*}

The integral over $d\b d\a'$ are just Fourier transforms, hence we can integrate over $\a'$ to get $\b=t_1-s_1$:
\begin{eqnarray*}
&=&\int_\R\int_{\R+i\e_\t'}\int_{\R+i\e_\t-i\e_\t'}\int_\R \int_{\R+i\e_\s} d\s d\a d\t d\t' d\b'\\ &&f(s_1+\a,t_1+\a+\t)g(s_2+\s,t_2)\frac{G_b(\frac{Q}{2}+it_1-is_1+i\t)G_b(i\t'-i\t)G_b(-i\s)G_b(-i\t')}{G_b(\frac{Q}{2}+it_1-is_1)}\\
&&e^{2\pi i\b'(\a+\s+\t)}e^{2\pi i(\a+s_1-s_2)(\t'-\t)-2\pi i\a(t_1-s_1)+\pi i\t'^2}\hat{B}^{ib\inv(\s+\t)}\hat{A}^{ib\inv (t_1+s_2-s_1)}.\\
\end{eqnarray*}
Next from the decay properties for $G_b(i\t'-i\t)G_b(-i\t')$ this is integrable over $\t'$. By shifting $\t\mapsto \t-\s-\a$ we see also that the integrations over $\s,\t, \a$ are absolutely convergent, and furthermore integrations over $\a, \s$ do not depend on $\b'$. Hence we can interchange the order for $d\t d\t'$ and then $d\s d\a$ to get
$$\int_\R \int_{\R+i\e_\s}\int_\R\int_{\R+i\e_\t-i\e_\t'}\int_{\R+i\e_\t'}d\t' d\t d\b'  d\s d\a .$$
Using the reflection properties, we bring $G_b(-i\t')$ to the denominator. Notice that the contour goes above $\t'=0$ and below $\t'=\t$, hence we can integrate over $\t'$ using Lemma \ref{tau} to get

\begin{eqnarray*}
&=&\int_\R \int_{\R+i\e_\s}\int_\R\int_{\R+i\e_\t-i\e_\t'}d\t d\b'  d\s d\a f(s_1+\a,t_1+\a+\t)g(s_2+\s,t_2)\\ &&\frac{G_b(\frac{Q}{2}+it_1-is_1+i\t)G_b(-i\s)G_b(-i\t)G_b(\frac{Q}{2}-i\a+is_2-is_1)}{G_b(\frac{Q}{2}+it_1-is_1)G_b(\frac{Q}{2}-i\a-i\t+is_2-is_1)}\\
&&e^{2\pi i\b'(\a+\s+\t)}e^{-2\pi i(\a+s_1-s_2)\t-2\pi i\a(t_1-s_1)} \hat{B}^{ib\inv(\s+\t)}\hat{A}^{ib\inv (t_1+s_2-s_1)}\\
&=&\int_\R \int_{\R+i\e_\s}\int_\R\int_{\R+i\e_\t-i\e_\t'}d\t d\b'  d\s d\a f(s_1-\a,t_1-\a+\t)g(s_2+\s,t_2)\\ &&\frac{G_b(\frac{Q}{2}+is_1-it_1)G_b(-i\s)G_b(-i\t)G_b(\frac{Q}{2}+is_1-is_2-i\a+i\t)}{G_b(\frac{Q}{2}+is_1-it_1-i\t)G_b(\frac{Q}{2}+is_1-is_2-i\a)}\\
&&e^{2\pi i\b'(\s+\t-\a)}e^{2\pi i (t_1-s_1)(\t-\a)} \hat{B}^{ib\inv(\s+\t)}\hat{A}^{ib\inv (t_1+s_2-s_1)},\\
\end{eqnarray*}
where we used the reflection properties again, and flipping $\a\mapsto -\a$. Now we can integrate $\t$ and $\b'$ by Fourier transform to get $\t=\a-\s$. However, this should be interpreted as a function of $\s$ over $\R$, and then analytic continued to $\R+i\e_\s$, so that the contour for $\s$ will be pushed under the pole for $\s=\a$ in $G_b(-i\t)\mapsto G_b(i\s-i\a)$. We get
\begin{eqnarray*}
&=&\int_\R \int_C d\s d\a f(s_1-\a,t_1-\s)g(s_2+\s,t_2)\\ &&\frac{G_b(\frac{Q}{2}+is_1-it_1)G_b(-i\s)G_b(i\s-i\a)G_b(\frac{Q}{2}+is_1-is_2-i\s)}{G_b(\frac{Q}{2}+is_1-it_1-i\a+i\s)G_b(\frac{Q}{2}+is_1-is_2-i\a)}\\
&&e^{2\pi i (t_1-s_1)\s} \hat{B}^{ib\inv \a}\hat{A}^{ib\inv (t_1+s_2-s_1)}.\\
\end{eqnarray*}

Finally we apply $W_{31}$ to get
\begin{eqnarray*}
&\mapsto_W& \int_{\R+i\e_\t} \int_\R \int_C d\s d\a d\t f(s_1-\a+\t,t_1-\s)g(s_2+\s,t_2)\\
&&\frac{G_b(-i\s)G_b(-i\a+i\s)G_b(\frac{Q}{2}+is_1-it_1+i\t)G_b(\frac{Q}{2}+is_1-is_2-i\s+i\t)}{G_b(Q+i\t)G_b(\frac{Q}{2}+is_1-it_1-i\a+i\s+i\t)G_b(\frac{Q}{2}+is_1-is_2-i\a+i\t)}\\
&&e^{2\pi i\s(t_1-s_1-\t)}e^{2\pi is_1\t}A^{ib\inv s_1}B^{ib\inv \t}\ox\hat{B}^{ib\inv\a}\hat{A}^{ib\inv (t_1+s_2-s_1-\t).}
\end{eqnarray*}
\end{proof}

We can further simplify the expression by the following unitary transformations, which will also be used in the right picture:
\begin{eqnarray*}
s_2&\mapsto& s_2-t_1,\\
s_1&\mapsto& s_1+\frac{s_2}{2},\\
t_1&\mapsto& t_1+\frac{s_2}{2},\\
\end{eqnarray*}
or equivalently
\Eq{f\ox g\mapsto f(s_1+\frac{s_2}{2},\frac{s_2}{2}-t_1)g(t_1+\frac{s_2}{2},t_2),}
so that the expression becomes:
\begin{eqnarray}
&=& \int_{\R+i\e_\t} \int_\R \int_C d\s d\a d\t  f(s_1-\a+\t,t_1-\s)g(s_2,t_2)\nonumber\\
&&\frac{G_b(-i\s)G_b(-i\a+i\s)G_b(\frac{Q}{2}+is_1-it_1+i\t)G_b(\frac{Q}{2}+is_1+it_1-i\s+i\t)}{G_b(Q+i\t)G_b(\frac{Q}{2}+is_1-it_1-i\a+i\s+i\t)G_b(\frac{Q}{2}+is_1+it_1-i\a+i\t)}\nonumber\\
&&e^{2\pi i\s(t_1-s_1-\t)}e^{2\pi is_1\t}e^{\pi is_2\t}A^{ib\inv (s_1+\frac{s_2}{2})}B^{ib\inv \t}\ox\hat{B}^{ib\inv\a}\hat{A}^{ib\inv (\frac{s_2}{2}-s_1-\t)}.\nonumber\\\label{leftcorep}
\end{eqnarray}

\subsection{The right regular corepresentation}\label{sec:double:right}
In analogy to the left regular action, there is a notion of a right regular action. First we consider the quantum plane.
\begin{Prop} Given $W$ that defines the left regular corepresentation, 
\Eq{W_R:=(\hat{J}\ox\hat{J}) W_{21}^* (\hat{J}\ox\hat{J})\in \cB(\cH)\ox \cA}
is also a multiplicative unitary and satisfies
\Eq{(1\ox \D)W_R=W_{R,12}W_{R,13},}
i.e. it defines a corepresentation, called the right-regular corepresentation.
\end{Prop}
Using the realization of $Ad \hat{J}$, we see that $W_R$ is given by 
\Eq{W_R=e^{\frac{i}{2\pi b^2}\log A_*\ox \log A}g_b(B_*\ox B)=e^{\frac{i}{2\pi b^2}\log A_*\ox \log A}\int_{\R+i0} B_*^{ib\inv \t}\ox B^{ib\inv \t}\frac{e^{-\pi i\t^2}}{G_b(Q+i\t)}d\t,}
where on the space $\cH_{corep}$ we have 
$\hat{J}=\over[f](t,s)$ so that
\Eq{A_* = e^{-2\pi bt},\tab B_* = e^{-2\pi bp_t}.}
Therefore we see that the corepresentation is acting on the $t$-coordinate by:
\Eq{f(s,t)\mapsto \int_{\R+i0} f(s,t-\t)\frac{e^{-\pi i\t^2}}{G_b(Q+i\t)}A^{-ib\inv t}B^{ib\inv \t}d\t,}
in which under the pairing we have:
\Eq{X_r = e^{-2\pi b t}, \tab Y_r=e^{-2\pi b p_t},}
\Eq{\til[X]_r=e^{-2\pi b\inv t},\tab \til[Y]_r=e^{-2\pi b\inv p_t}.}
In other words, under the transform $t\mapsto -t$, we conclude together with Proposition \ref{Bqaction} that
\begin{Prop} We have the equivalence
\Eq{L^2(\cA)\simeq \cH_{irr}\ox \cH_{irr}}
as representation of $\cB_{q\til[q],L}\ox \cB_{q\til[q],R}\simeq \cB_{q\til[q]}\ox \cB_{q\til[q]}$, where $\cH_{irr}:=L^2(\R)$ is the canonical representation of $\cB_{q\til[q]}$ on $L^2(\R)$. Hence this is an analogue to the classical ``Peter-Weyl" theorem on ``functions on the quantum plane".
\end{Prop}

Now let us look at the case for $\cD(\cA)$. To define the right regular corepresentation, we need the corresponding multiplicative unitary $\bW_R$. It is known that $\bW_R$ is given by
\Eq{\bW_R=(\cU\ox 1)\bW_{3412}(\cU^*\ox 1),}
where 
\Eq{\cU=\hat{\bJ}_m\bJ_m = (J\hat{J}\ox \hat{J}J)Z\in \cA\ox \hat{\cA}.}
We simplify the expression and express $\bW_R$ in terms of 4 $g_b$'s as before.

\begin{Prop} $\bW_R$ is given by (cf. \eqref{Wmsimple})
\Eq{\bW_R=V_{32}W_{R,13}V_{32}^{''*}\hat{W}_{R,24},}
where
$V^{''*} = (\hat{J}\ox\hat{J})V^*(\hat{J}\ox\hat{J}).$
\end{Prop}
\begin{proof}
From (\cite{KV2} Prop 2.15), we have:
\begin{eqnarray*}
\hat{W}&=&W_{21}^*,\\
W_R&=&(\hat{J}\ox\hat{J})W_{21}^*(\hat{J}\ox\hat{J})\\
&=&(\hat{J}\ox\hat{J})(J\ox\hat{J})W_{21}(J\ox\hat{J})(\hat{J}\ox\hat{J})\\
&=&(\hat{J}J\ox 1)W_{21}(J\hat{J}\ox 1),\\
\hat{W}_R&=&(J\ox J)W(J\ox J)\\
&=&(J\ox J)(\hat{J}\ox J)W^*(\hat{J}\ox J)(J\ox J)\\
&=&(J\hat{J}\ox 1)W^*(\hat{J}J\ox 1).\\
\end{eqnarray*}
Hence we obtain (the indices indicate the legs in which the operators are acting):
\begin{eqnarray*}
\bW_R&=& ((J_1\hat{J}_1\hat{J}_2J_2)Z_{12})W_{31}Z_{12}^*\hat{W}_{42}Z_{12}(Z_{12}^* (\hat{J}_1J_1 J_2\hat{J}_2))\\
&=&(J_1\hat{J}_1\hat{J}_2J_2)Z_{12}W_{31}Z_{12}^*\hat{W}_{42}(\hat{J}_1J_1 J_2\hat{J}_2)\\
&=&(J_1\hat{J}_1\hat{J}_2J_2)V_{12}V_{12}'W_{31}V_{12}^{*'}V_{12}^*\hat{W}_{42}(\hat{J}_1J_1 J_2\hat{J}_2)\\
&=&(J_1\hat{J}_1\hat{J}_2J_2)(\hat{J}_2\hat{J}_3)V_{32}^*(\hat{J}_2\hat{J}_3)V_{12}W_{31}V_{12}^*\hat{W}_{42}(\hat{J}_1J_1 J_2\hat{J}_2)\\
&=&(J_1\hat{J}_1\hat{J}_2J_2)(\hat{J}_2\hat{J}_3)V_{32}^*(\hat{J}_2\hat{J}_3)W_{31}V_{32}\hat{W}_{42}(\hat{J}_1J_1 J_2\hat{J}_2)\\
&=&(J_1\hat{J}_1\hat{J}_2J_2)(\hat{J}_2\hat{J}_3)V_{32}^*(\hat{J}_2\hat{J}_3)W_{31}V_{32}(J_1 \hat{J}_1J_2\hat{J}_2)\hat{W}_{R,24}\\
&=&(J_1\hat{J}_1\hat{J}_2J_2)(\hat{J}_2\hat{J}_3)V_{32}^*(\hat{J}_2\hat{J}_3)(\hat{J}_1J_1 )W_{R,13}V_{32}(J_2\hat{J}_2)\hat{W}_{R,24}\\
&=&(\hat{J}_2J_2)(\hat{J}_2\hat{J}_3)V_{32}^*(\hat{J}_2\hat{J}_3)W_{R,13}V_{32}(J_2\hat{J}_2)\hat{W}_{R,24}\\
&=&(J_2\hat{J}_3)V_{32}^*(\hat{J}_2\hat{J}_3)W_{R,13}V_{32}(\hat{J}_2J_2)\hat{W}_{R,24}\\
&=&V_{32}(J_2\hat{J}_3)(\hat{J}_2\hat{J}_3)W_{R,13}V_{32}(\hat{J}_2J_2)\hat{W}_{R,24}\\
&=&V_{32}W_{R,13}(J_2\hat{J}_2)V_{32}(\hat{J}_2J_2)\hat{W}_{R,24}\\
&=&V_{32}W_{R,13}(J_2\hat{J}_2)(J_2\hat{J}_3)V_{32}^*(J_2\hat{J}_3)(\hat{J}_2J_2)\hat{W}_{R,24}\\
&=&V_{32}W_{R,13}(\hat{J}_2\hat{J}_3)V_{32}^*(\hat{J}_2\hat{J}_3)\hat{W}_{R,24}.\\
\end{eqnarray*}
\end{proof}

Now we can realize the action on $\cH_{corep}\ox\hat{\cH}_{corep}$ as in previous section:
\begin{Prop} The action of the components of $\bW_R$ is given by
\begin{eqnarray*}
\hat{W}_{R,24}&=&e^{-\frac{i}{2\pi b^2}\log \hat{A}_*\ox \log \hat{A}}g_b^*(\hat{B}_*\ox \hat{B})\\
f(s_1,t_1)g(s_2,t_2)&\mapsto& \int f(s_1,t_1)g(s_2,t_2-\t)e^{-2\pi i\t t_2+\pi i\t^2}G_b(-i\t)\hat{B}^{ib\inv \t} \hat{A}^{-ib\inv t_2}d\t,\\
W_{R,13}&=&e^{\frac{i}{2\pi b^2}\log A_*\ox \log A}g_b(B_*\ox B)\\
f(s_1,t_1)g(s_2,t_2)&\mapsto& \int f(s_1,t_1-\t)g(s_2,t_2)\frac{e^{-\pi i\t^2}}{G_b(Q+i\t)} A^{-ib\inv t_1}B^{ib\inv \t}d\t,\\
V_{32}&=&e^{\frac{i}{2\pi b^2}\log \hat{A}\inv\ox \log A\inv}g_b(q\hat{B}\hat{A}\inv\ox B)\\
f(s_1,t_1)g(s_2,t_2)&\mapsto&\int_\R\int_{\R+i0}\int_\R f(s_1,t_1)g(s_2+\a,t_2+\a+\t)\\
&&\frac{e^{2\pi i(\a\b+\a\t+\b\t+s_2\t)}G_b(\frac{Q}{2}+is_2-it_2)}{G_b(\frac{Q}{2}+is_2-it_2-i\t)G_b(Q+i\t)}A^{ib\inv \b}B^{ib\inv \t}d\a d\t d\b,\\
V_{32}^{''*}&=&e^{\frac{i}{2\pi b^2}\log A_*\ox \log A}g_b(B_*\ox B)\\
f(s_1,t_1)g(s_2,t_2)&\mapsto&\int_\R\int_{\R+i0}\int_\R g(s_2+\a,t_2+\a+\t)\\
&&\frac{e^{2\pi i\a(t_2-s_2-\b)}e^{-\pi i\t^2}G_b(\frac{Q}{2}+is_2-it_2)}{G_b(\frac{Q}{2}+is_2-it_2-i\t)G_b(Q+i\t)}A^{ib\inv \b}B^{ib\inv \t}d\a d\t d\b.\\
\end{eqnarray*}
\end{Prop}

\begin{Thm} The right regular corepresentation is given by
\begin{eqnarray*}
f(s_1,t_1)g(s_2,t_2)&\mapsto&\int_{\R+i\e_\t} \int_\R \int_C d\s d\a d\t f(s_1,t_1-\s)g(s_2+\s,t_2+\t-\a)\\
&&e^{\pi Q(-\a+\t)+2\pi i(-\s s_2+\a \t-\a t_2+\s t_2)}\\
&&\frac{G_b(-i\s)G_b(\frac{Q}{2}+is_2-it_2)G_b(i\s-i\t)G_b(\frac{Q}{2}-i\s+it_1-it_2)}{G_b(\frac{Q}{2}+is_2-it_2+i\s-i\t)G_b(Q+i\a)G_b(\frac{Q}{2}-i\t+it_1-it_2)}\\
&&A^{ib\inv(t_2-t_1-s_2)}B^{ib\inv \t}\hat{B}^{ib\inv \a} \hat{A}^{-ib\inv (t_2+\t)} d\a d\t d\s,
\end{eqnarray*}
where the contour $C$ of $\s$ goes above $\s=0$ and below $\s=\t$.
\end{Thm}

We have to follow the transformations given in the left regular picture. So we apply the transformations after Theorem \ref{leftco} and an extra one that does not affect the previous action:
\begin{eqnarray*}
s_2&\mapsto& s_2-t_1,\\
s_1&\mapsto& s_1+\frac{s_2}{2},\\
t_1&\mapsto& t_1+\frac{s_2}{2},\\
t_2&\mapsto& t_2+\frac{s_2}{2},\\
\end{eqnarray*}
and obtain
\begin{eqnarray}
f\ox g&\mapsto&\int_{\R+i\e_\t} \int_\R \int_C d\s d\a d\t f(s_1,t_1-\s)g(s_2,t_2+\t-\a)\nonumber\\
&&e^{\pi Q(-\a+\t)+2\pi i(\s t_1+\a \t-\a t_2+\s t_2)}e^{-\pi i \a s_2}\nonumber\\
&&\frac{G_b(-i\s)G_b(\frac{Q}{2}-it_1-it_2)G_b(i\s-i\t)G_b(\frac{Q}{2}-i\s+it_1-it_2)}{G_b(\frac{Q}{2}-it_1-it_2+i\s-i\t)G_b(Q+i\a)G_b(\frac{Q}{2}-i\t+it_1-it_2)}\nonumber\\
&&\label{rightcorep}A^{ib\inv(t_2-\frac{s_2}{2})}B^{ib\inv \t}\hat{B}^{ib\inv \a} \hat{A}^{ib\inv (-t_2-\t-\frac{s_2}{2})}d\a d\t d\s.
\end{eqnarray}


\section{Regular representation}\label{sec:regular}
To obtain a representation of $U_q(\sl(2,\R))$ on the space $L^2(GL_q^+(2,\R))$, we need to describe the non-degenerate pairing between them. Then we apply the pairing to the corepresentations constructed above and obtain the desired representations for $U_q(\gl(2,\R))$.


\subsection{Pairing with $U_q(\gl(2,\R))$}\label{sec:regular:pairing}
Let us first recall the following definition of $U_q(\gl(2,\R))$ serving as a dual space for $\cM_q=GL_q^+(2,\R)$ that is observed by I. Frenkel \cite{F}:
\begin{Def} As a Hopf *-algebra, $U_q(\gl(2,\R))$ is generated by positive self-adjoint operators $E,F,K,K_0$ such that
\begin{eqnarray*}
KE&=&q EK,\\
KF&=&q\inv FK,\\
EF-FE&=&\frac{K^2-K^{-2}}{q-q\inv},\\
\D(K)&=&K\ox K,\\
\D(K_0)&=&K_0\ox K_0,\\
\D(E)&=&K_0\inv K\inv\ox E+E\ox K_0 K,\\
\D(F)&=&K_0K\inv \ox F+F\ox K_0\inv K,
\end{eqnarray*}
and $K_0$ commutes strongly with $E,F,K$.
\end{Def}
In particular, we see that by setting $K_0=1$, we obtain the usual definition of $U_q(\sl(2,\R)$.
Now we can similarly bring $U_q(\gl(2,\R))$ into $C^*$-algebraic level by introducing continuous parameters. It is the algebra
\Eq{\{\int_\R\int_\R\int_{\R+i0}\int_{\R+i0} f(r_0,r,s,t)K_0^{ib\inv r_0}K^{ib\inv r}E^{ib\inv s}F^{ib\inv t}dsdtdr_0dr\},}
where $f(r_0,r,s,t)$ has simple poles at $s,t=-inb-imb\inv, n,m\in\Z_{\geq0}$, and has similar decays as $\cA$ along the real direction in the variables $r_0,r,s,t$. This is to ensure that the coproduct lies in the multiplier algebra as argued before. The $C^*$-norm can be introduced once we obtain its representation as certain unitary operators on $L^2(\R\x \R)$ under the pairing. By abuse of notation we still denote it by $U_q(\gl(2,\R))$, and a similar version with $K_0$ omitted by $U_q(\sl(2,\R))$.

Let us ignore the $*$-structure and consider only the polynomial algebra. Then there exists a non-degenerate pairing defined on the generators:
\begin{Prop}
The pairing between $E,F,K,K_0$ and the $z$-variable $\z[ij],i=1,2$ is given by
\begin{eqnarray*}
\<E,\z[21]\>=c,&\tab& \<F,\z[12]\>=c\inv,\\
\<K,\z[11]\>=q^{-1/2},&\tab& \<K,\z[22]\>=q^{1/2},\\
\<K_0,\z[11]\>=q^{-1/2},&\tab& \<K_0,\z[22]\>=q^{-1/2},\\
\<E,N\>=0,&\tab&\<F,N\>=0,\\
\<K,N\>=1,&\tab& \<K_0,N\>=q\inv,
\end{eqnarray*}
and zero otherwise, where $c\in\C$ is any constant. The pairing is then extended to any monomial by the coproduct and induction.
\end{Prop}

\begin{Cor}
The pairing between $E,F,K,K_0$ and the $\cD(\cA)$ generators is given by
\begin{eqnarray*}
\<K,A\>=q^{-\frac{1}{2}},&\tab& \<K,\hat{A}\>=q^{\frac{1}{2}},\\
\<K_0,A\>=q^{-\frac{1}{2}},&\tab &\<K_0,\hat{A}\>=q^{-\frac{1}{2}},\\
\<E,B\>=c,&\tab& \<F,\hat{B}\>=c\inv,
\end{eqnarray*}
and zero otherwise, where $c\in\C$ is any constant. The pairing is then extended to any monomial by the coproduct and induction.
\end{Cor}

By induction, we obtain
\begin{Prop} The pairing between any monomials is given by
\Eq{\<K^l E^mF^n, \z[11]^{L} \z[21]^{m'}\z[12]^{n'}N^{L'}\>=c^{m^2-n^2}q^{l(m-n-L)/2+mL+(n-m)L'-2nm}[n]_q![m]_q! \d_{nn'}\d_{mm'},}
or
\Eq{\<K^l E^mF^n, A^{L} B^{m'}\hat{B}^{n'}\hat{A}^{L'}\>=c^{m^2-n^2}q^{l(m+L'-L)/2+mL+nL'-nm}[n]_q![m]_q! \d_{nn'}\d_{mm'}.}
\end{Prop}

\begin{Rem} Restricting to $U_q(\sl(2,\R))$, the pairing is essentially the same as the one given in \cite{MMNNSU}, with the role of $E$ and $F$ interchanged. Hence this provides a more elegant formula for the pairing of a general monomial.
\end{Rem}

Hence as in Section \ref{sec:Dual:pairing}, we introduce $G_b$ for the expression in the general pairing with elements represented by integrations. It turns out that we need to replace $[n]_q!$ by 
\Eq{[n]_q!\mapsto\frac{G_b(Q+i\t)}{(1-q^2)^{ib\inv \t}}} so that the commutation relation $[E,F]=\frac{K^2-K^{-2}}{q-q\inv}$ is satisfied. Here we take $1-q^2$ as a complex number with $-\frac{\pi}{2}<\arg(1-q^2)<\frac{\pi}{2}$.

Therefore the general pairing is given by
\begin{eqnarray}
&&\< f(r_0,r,s,t),g(s_1,t_1)h(s_2,t_2)\>\nonumber\\
&=&\iint f(r_0,r,t_1,t_2)g(s_1,t_1)h(s_2,t_2)e^{-\pi i(r(t_1+s_2-s_1)/2-r_0(s_2+s_1)/2+t_1s_1+t_2s_2-t_1t_2)}\nonumber\\
&&c^{b^{-2}(t_2^2-t_1^2)}\frac{G_b(Q+it_1)G_b(Q+it_2)}{(1-q^2)^{ib\inv(t_1+t_2)}}ds_1ds_2dt_1dt_2dr_0 dr,
\end{eqnarray}
where we denoted by $$f(r_0,r,s,t):=\iiint f(r,s,t)K_0^{ib\inv r_0}K^{ib\inv r}E^{ib\inv s}F^{ib\inv t}dsdtdr_0dr,$$ and $$g(s_1,t_1)h(s_2,t_2):=\iint\iint g(s_1,t_1)h(s_2,t_2)A^{ib\inv s_1}B^{ib\inv t_1}\hat{B}^{ib\inv t_2}\hat{A}^{ib\inv s_2}ds_1ds_2dt_1dt_2.$$
Finally, in order to obtain a positive representation in later sections, it turns out we need to choose the pairing constant to be $c=-iq^{1/2}=-q^{\frac{Q}{2b}}$.
\subsection{The principal series representation}\label{sec:regular:principal}
By taking the pairing above with the fundamental corepresentation defined in Section \ref{sec:double:matrix}, we obtain:
\begin{Thm}\label{fundrep} The fundamental corepresentation $T^{\l,t}$ corresponds to the representation $\cP_{\l,t}$, called the \emph{principal series representation}, given on $L^2(\R,ds)$ by:
\begin{eqnarray*}
K_0 &=& e^{\pi bt},\\
K&=&e^{-\pi bs},\\
E &=&[\frac{Q}{2b}-\frac{i}{b}(s-\l)]_q e^{-2\pi bp_{s}},\\
F &=&[\frac{Q}{2b}+\frac{i}{b}(s+\l)]_q e^{2\pi bp_{s}},\\
\end{eqnarray*}
where $[n]_q=\frac{q^n-q^{-n}}{q-q\inv}$. The operators are all positive essentially self-adjoint.
\end{Thm}
\begin{proof}
Recall that we choose $c=-iq^{1/2}$ in the pairing.  Recall also that we have to multiply the expressions by $\frac{G_b(Q+i\a)G_b(Q+i\t)}{(1-q^2)^{ib\inv (\a+\t)}}$. From the expression \eqref{fundcorep} we obtain for $K_0$ and $K$ that $\t=\a=0$. The factor becomes $\dis\lim_{\a\to0}G_b(-i\a)G_b(Q+i\a)=1$ and the rest becomes
\begin{eqnarray*} 
K_0:f(s)&\mapsto&f(s)\cdot (q^{-1/2})^{ib\inv(t-s)}(q^{-1/2})^{ib\inv(t+s)}\\
&=&e^{\pi bt}f(s),\\
K:f(s)&\mapsto&f(s)\cdot (q^{-1/2})^{ib\inv(t-s)}(q^{1/2})^{ib\inv(t+s)}\\
&=&e^{-\pi bs}f(s).\\
\end{eqnarray*}
For $E$, we have $\t=-ib,\a=0$ so that the factor $\dis\lim_{\a\to0}G_b(-i\a)G_b(Q+i\a)=1$ as well, and we obtain
\begin{eqnarray*}
E:f(s)&\mapsto&\frac{-iq^{1/2}f(s+ib)}{1-q^2}\frac{G_b(\frac{Q}{2}-is+i\l+b)}{G_b(\frac{Q}{2}-is+i\l)}e^{\pi \l b-2\pi bs+\pi bt}\cdot (q^{ib\inv(t-s)})\\
&=&\frac{iq^{-1/2}}{q-q\inv} (1-e^{2\pi i b(\frac{Q}{2}-is+i\l)})e^{\pi b(\l-s)} f(s+ib)\\
&=&\frac{i}{q-q\inv}(q^{1/2}e^{\pi b(s-\l)}+q^{-1/2}e^{-\pi b(s-\l)}) f(s+ib)\\
&=&[\frac{Q}{2b}-\frac{i}{b}(s-\l)]_q f(s+ib).
\end{eqnarray*}
Finally, for $F$ we have $\a=-ib,\t=0$, the factor $\dis\lim_{\a\to-ib}G_b(-i\a)G_b(Q+i\a)=-q^2$, and we obtain
\begin{eqnarray*}
F:f(s)&\mapsto&\frac{iq^{-1/2}(-q^2)f(s-ib)}{1-q^2}\frac{G_b(\frac{Q}{2}-is-i\l)}{G_b(\frac{Q}{2}-is-i\l-b)}e^{-\pi \l b+\pi bt}\cdot (q^{ib\inv(t+s)})\\
&=&\frac{iq^{1/2}}{q-q\inv} (1-e^{2\pi i b(\frac{Q}{2}-is-i\l-b)})e^{-\pi b(\l+s)} f(s-ib)\\
&=&\frac{i}{q-q\inv}(q^{-1/2}e^{\pi b(s+\l)}+q^{1/2}e^{-\pi b(s+\l)}) f(s-ib)\\
&=&[\frac{Q}{2b}+\frac{i}{b}(s+\l)]_q f(s-ib).
\end{eqnarray*}
Since $\frac{i}{q-q\inv}=\frac{1}{2\sin(\pi b^2)}>0$, we immediately see that the operators are positive. We note that the expression for $E$ can be rewritten as
\Eq{E=\frac{i}{q-q\inv} (e^{\pi b(s-\l-p_s)}+e^{\pi b(\l-s-p_s)}),}
where the summand $q^2$-commutes. Hence using Proposition \ref{x+p}, it is unitary equivalent to $e^{\pi b(s-\l-p_s)}$ which in turn is equivalent to $e^{-\pi b(p_s+\l)}$ by multiplication by $e^{-\pi is^2}$. Hence it is essentially self-adjoint. Similar analysis applies to $F$, hence all operators are positive essentially self-adjoint.
\end{proof}

\begin{Rem} We note that by restricting to $U_q(\sl(2,\R))$, this is precisely (the Fourier transform of) the continuous series representation $\cP_{s}$ obtained in \cite{BT} given by
\begin{eqnarray*}
K\d_k&=&e^{-\pi bk}\d_k,\\
E_s\d_k&=&[\frac{Q}{2b}-\frac{i}{b}(k-s)]_q \d_{k+ib},\\
F_s\d_k&=&[\frac{Q}{2b}+\frac{i}{b}(k+s)]_q\d_{k-ib}\\
\end{eqnarray*}
under the correspondences $s\corr k, \l\corr s$.
\end{Rem}
\subsection{Casimir operator}\label{sec:regular:Casimir}
Next we will study the action induced from the pairing with the multiplicative unitary $\bW$. It turns out that the technical difficulty comes from the analysis of the following positive operator:
\Eq{\bC=e^{2\pi bx}+e^{-2\pi bx}+e^{-2\pi bp},}
which is studied in detail in \cite{Ka4} in a different context. This expression comes from the Casimir operator defined by 
\Eq{C=FE+\frac{qK^2+q\inv K^{-2}-2}{(q-q\inv)^2}.}
More precisely, under the pairing given in the next section, the Casimir operator is given by
\Eq{C=\left(\frac{i}{q-q\inv}\right)^2 (e^{2\pi bx}+e^{-2\pi bx}+e^{-2\pi bp}+2).}
\begin{Prop} For $\l \in \R^+$, 
\Eq{\Phi_\l(s_2):=S_b(-ix+i\l)S_b(-ix-i\l)}
is an eigenfunction for the operator $\bC$ with eigenvalue $e^{2\pi b\l}+e^{-2\pi b\l}$.
\end{Prop}
\begin{proof} We need to solve
\Eq{\bC\Phi_\l(x) = (e^{2\pi b\l}+e^{-2\pi b\l}) \Phi_\l(x).}

We consider $$\bC- (e^{2\pi b\l}+e^{-2\pi b\l}) = \bC-2\cos(2\pi ib\l)=2\cos(2\pi ib x)-2\cos(2\pi ib\l)+e^{2\pi b p_{s_2}}.$$

From the angle sum formula, we have
\begin{eqnarray*}
&&2\cos(2\pi ib x)-2\cos(2\pi ib\l)\\
&=&-4\sin(\pi ib(x-\l))\sin(\pi ib(x+\l)),\\
\end{eqnarray*}
and from the functional equation \eqref{funceq} $$S_b(x+b)=2\sin(\pi bx)S_b(x),$$ we immediately see that
\begin{eqnarray*}
&&e^{2\pi bp_{x}}\cdot S_b(-ix+i\l)S_b(-ix-i\l)\\
&=&S_b(-ix+i\l+b)S_b(-ix-i\l+b)\\
&=&4\sin(\pi ib(-x+\l))\sin(\pi ib(-x-\l))S_b(-ix+i\l)S_b(-ix-i\l)\\
&=&4\sin(\pi ib(x-\l))\sin(\pi ib(x+\l))S_b(-ix+i\l)S_b(-ix-i\l).
\end{eqnarray*}
Hence $\Phi_\l(x) =S_b(-ix+i\l)S_b(-ix-i\l)$ satisfies the above eigenvalue equation.
\end{proof}

Therefore we introduce the following transformation
\begin{Def} We define the integral transformation
$$\Phi\inv: L^2(\R^+, d\mu(\l))\to L^2(\R)$$
\Eq{F(\l)\mapsto f(x):=\lim_{\e\to 0}\int_0^\oo \Phi_\l(x+i\e)e^{-2\pi x\e} F(\l) d\mu(\l)}
where $d\mu(\l)=|S_b(Q+2i\l)|^2d\l$.
\end{Def}

\begin{Thm}\label{Casimir} $\Phi\inv$ is a unitary transformation that intertwines $e^{2\pi b\l}+e^{-2\pi b\l}$ and $\bC$.
\end{Thm}
\begin{proof} We compute the inner product for $f,g\in \cW^+$ where $\cW^+$ is the dense subspace
$$\cW^+:=\{f(x)\in C_\oo(\R^+): F(y):=f(e^y)\in \cW\}.$$
Note that
$$\over[\Phi_\l(x)]=\frac{1}{S_b(Q-ix-i\l)S_b(Q-ix+i\l)}.$$ We have
\begin{eqnarray*}
&&\<\Phi\inv f, \Phi\inv g\>\\
&=&\lim_{\e\to 0}\int_\R \int_0^\oo \int_0^\oo \frac{S_b(\e-ix-i\l)S_b(\e-ix+i\l)e^{-4\pi x\e}}{S_b(Q-\e-ix-i\b)S_b(Q-\e-ix+i\b)}f(\l)\over[g(\b)]d\mu(\l) d\mu(\b) dx\\
&=&\lim_{\e\to 0}\int_\R \int_0^\oo \int_0^\oo \frac{G_b(\e+ix-i\l)G_b(\e+ix+i\l)e^{\pi i\l^2-\pi i\b^2}e^{-2\pi x(Q-4\e)}}{G_b(Q-\e+ix-i\b)G_b(Q-\e+ix+i\b)}f(\l)\over[g(\b)]d\mu(\l) d\mu(\b) dx\\
\end{eqnarray*}
which means the contour for $x$ separates the poles of the numerator and denominator.

Now by Corollary \ref{asymp} the integrand has asymptotics in $x$
$$\left|\frac{G_b(\e+ix-i\l)G_b(\e+ix+i\l)e^{\pi i\l^2-\pi i\b^2}e^{-2\pi x(Q-4\e)}}{G_b(Q-\e+ix-i\b)G_b(Q-\e+ix+i\b)}\right|=\left\{\begin{array}{cc}e^{-2\pi x(Q-4\e)}&x\to\oo\\e^{2\pi Qx}&x\to-\oo\end{array}.\right.$$
Hence the integral is absolutely convergent in $x,\l,\b$ and we can interchange the order of integration to bring integration over $x$ inside.

Shifting $x\mapsto x+\b$, and using the 4-5 relation (Lemma \ref{45}), we obtain
\begin{eqnarray*}
&=&\lim_{\e\to 0}\int_0^\oo \int_0^\oo \frac{G_b(Q-4\e)G_b(2\e+i\b-i\l)G_b(2\e+i\b+i\l)}{G_b(Q-2\e+i\b-i\l)G_b(Q-2\e+i\b+i\l)}\cdot\\
&&e^{\pi i\l^2-\pi i\b^2}e^{-2\pi Q\b+8\pi\b\e} f(\l)\over[g(\b)]|S_b(Q+2i\l)|^2|S_b(Q+2i\b)|^2d\l d\b.
\end{eqnarray*}
Note that the integration for $\b$ separates the poles of the numerator and denominator as well, and also since $\l,\b>0$, the case $\b=-\l$ can be ignored, hence we obtain $\d(\l-\b)$ using Lemma \ref{delta}:
\begin{eqnarray*}
&=&\int_0^{\oo}\frac{G_b(2i\l)e^{-2\pi Q\l}}{G_b(Q+2i\l)} f(\l)\over[g(\l)]|S_b(Q+2i\l)|^4d\l\\
&=&\int_0^{\oo} f(\l)\over[g(\l)]|S_b(2i\l)|^2 |S_b(Q+2i\l)|^4d\l\\
&=&\int_0^{\oo} f(\l)\over[g(\l)]|S_b(Q+2i\l)|^2d\l.\\
\end{eqnarray*}
This shows that the map is an isometry. The converse showing the eigenfunctions are also complete is more difficult, and can be found in \cite{Ka4}. The technique involves writing the measure $$|S_b(Q+2i\l)|^2=4\sinh(\pi b\l)\sinh(\pi b^{-1}\l)$$ as linear combinations of 4 exponentials, and calculate the integral directly using hints from hyperfunctions.

From the proof we therefore deduce the inverse map of the unitary transformation to be the complex conjugate
$$\Phi: L^2(\R)\to L^2(\R^+, d\mu(\l))$$
$$f(x)\mapsto F(\l):= \int_{\R-i0} \frac{f(x)}{S_b(Q-ix-i\l)S_b(Q-ix+i\l)} dx,$$
where the contour goes below the poles.

Since $\Phi_\l(x)$ are the eigenfunctions satisfying $\bC\Phi_\l(x) = (e^{2\pi b\l}+e^{-2\pi b\l})\Phi_\l(x)$, formally the integral transformation with $\Phi_\l(x)$ as kernel will intertwine the action.
\end{proof}

\subsection{Left regular representation}\label{sec:regular:left}
Now we apply the pairing to the left corepresentation \eqref{leftcorep}.

\begin{Thm} \label{ThmLeft}The representation of $U_q(\gl(2,\R))$ is equivalent to $$\int_\R\int_{\R^+}^\o+\cP_{\l,s_2/2}|S_b(Q+2i\l)|^2d\l ds_2\ox L^2(\R, dt_2),$$ where the representation space $\cP_{\l,s_2/2}=L^2(\R,ds_1)$ is defined in Theorem \ref{fundrep}.
\end{Thm}
\begin{proof} Recall that we have to multiply the expression by $\frac{G_b(Q+i\a)G_b(Q+i\t)}{(1-q^2)^{ib\inv (\a+\t)}}$ under the pairing. Also since the resulting expression is analytic in $\a, \t$ on the lower half plane, by taking $\a, \t$ to $0, -ib$ we mean the analytic continuation from $\R$ to the corresponding point, while respecting the contour of $\s$.

The pairing of $K_0$ and $K$ is then given by $\t=0, \a=0$.
The factor $$\frac{G_b(-i\s)G_b(-i\a+i\s)}{G_b(Q+i\t)} G_b(Q+i\a)G_b(Q+i\t)$$ becomes $\d(\s)$ which means we also set $\s=0$. Hence we obtain
\Eq{K_0\cdot f(s_1,t_1)g(s_2,t_2)= f(s_1,t_1)g(s_2,t_2)e^{-\frac{\pi b}{2}s_2},}
\Eq{K\cdot f(s_1,t_1)g(s_2,t_2)= f(s_1,t_1)g(s_2,t_2)e^{\pi bs_1}.}

The pairing of $E$ is given by $\t=-ib$ and $\a=0$, so that again $\s=0$ and is given by
\begin{eqnarray*}
E\cdot f(s_1,t_1)g(s_2,t_2)&=&-iq^{1/2} f(s_1-ib,t_1)g(s_2,t_2)\frac{e^{2\pi bs_1+\pi bs_2}e^{-\pi bs_1-\frac{\pi b}{2}s_2}}{1-q^2}\\
&=&-iq^{1/2}f(s_1-ib,t_1)g(s_2,t_2)\frac{e^{\pi bs_1+\frac{\pi b}{2}s_2}}{1-q^2},
\end{eqnarray*}
or simply
\Eq{E=\frac{i}{q-q\inv} e^{\frac{\pi b}{2}s_2} e^{\pi bs_1+2\pi bp_{s_1}}.}

Finally the pairing of $F$ is given by $\t=0, \a=-ib$. By Corollary \ref{delta} we have
\begin{eqnarray*}
&&\frac{G_b(-i\s)G_b(-i\a+i\s)}{G_b(Q+i\t)} G_b(Q+i\a)G_b(Q+i\t)\\
&=&e^{-\pi i\s^2-\pi Q\s}\frac{G_b(i\s-b)G_b(Q+b)}{G_b(Q+i\s)}\\
&=&-q^2\d(\s)-q^2\d(\s+ib),
\end{eqnarray*}
so that we obtain
\begin{eqnarray*}
F\cdot f(s_1,t_1)g(s_2,t_2)&=&iq^{-1/2}\frac{-q^2}{1-q^2}e^{\pi bs_1-\frac{\pi b}{2}s_2}\left(e^{2\pi b(t_1-s_1)}f(s_1+ib,t_1+ib)g(s_2,t_2) \right.\\
&&\left.+ f(s_1+ib,t_1)g(s_2,t_2)\frac{G_b(\frac{Q}{2}+is_1-it_1)G_b(\frac{Q}{2}+is_1+it_1)}{G_b(\frac{Q}{2}+is_1-it_1-b)G_b(\frac{Q}{2}+is_1+it_1-b)}\right).
\end{eqnarray*}
Using $\frac{G_b(ix)}{G_b(ix-b)} = (1-e^{2\pi ib(ix-b)})$, we expand the expression:
\begin{eqnarray*}
F&=&\frac{iq^{1/2}e^{-\frac{\pi b}{2}s_2}}{q-q\inv}\left(e^{-\pi b(s_1-2t_1)}e^{-2\pi b(p_{t_1}+p_{s_1})}+\right.\\
&&\left.(1+q\inv e^{-2\pi b(s_1-t_1)})(1+q\inv e^{2\pi b(-s_1-t_1)})e^{2\pi b s_1}e^{-2\pi bp_{s_1}}\right)\\
&=&\frac{ie^{-\frac{\pi b}{2}s_2}}{q-q\inv}\left(e^{-\pi b(2p_{s_1}-s_1)}+e^{-\pi b(2p_{s_1}+3s_1)}+(e^{2\pi bt_1}+e^{-2\pi bt_1}+e^{-2\pi b(p_{t_1}-t_1)})e^{-\pi b(2p_{s_1}+s_1)}\right),\\
\end{eqnarray*}
so that under multiplication by $e^{-\pi it_1^2}$ which changes $p_{t_1}\mapsto p_{t_1}+t_1$, we get

\Eq{F=\frac{ie^{-\frac{\pi b}{2}s_2}}{q-q\inv}\left(e^{-\pi b(2p_{s_1}-s_1)}+e^{-\pi b(2p_{s_1}+3s_1)}+\bC e^{-\pi b(2p_{s_1}+s_1)}\right),}
where $\bC$ depending only on $t_1$ is the operator defined in the previous section.
Therefore under the unitary transformation $\Phi$ defined previously, the representation is equivalent to
\begin{eqnarray*}
K_0&=&e^{\frac{\pi b}{2}s_2},\\
K&=&e^{\pi bs_1},\\
E&=&\frac{i}{q-q\inv} e^{\frac{\pi b}{2}s_2} e^{\pi bs_1+2\pi bp_{s_1}},\\
F&=&\frac{ie^{-\frac{\pi b}{2}s_2}}{q-q\inv}\left(e^{-\pi b(2p_{s_1}-s_1)}+e^{-\pi b(2p_{s_1}+3s_1)}+(e^{2\pi b\l}+e^{-2\pi b\l})e^{-\pi b(2p_{s_1}+s_1)}\right),
\end{eqnarray*}
where the measure $d\mu(\l)$ is now $|S_b(Q+2i\l)|^2d\l$.

Motivated from the expression of the corepresentation, we multiply the space by $G_b(\frac{Q}{2}+is_1-i\l)\inv$, and in effect changes the action of $e^{2\pi bp_{s_1}}$ to $(1+qe^{-2\pi bs_1+2\pi b\l})e^{2\pi bp_{s_1}}$. The action of $E$ and $F$ then become:

$$E=\frac{i}{q-q\inv}e^{\frac{\pi b}{2}s_2}(e^{\pi bs_1+2\pi bp_{s_1}}+e^{-\pi bs_1+2\pi b\l+2\pi bp_{s_1}}),$$
$$F=\frac{i}{q-q\inv}e^{-\frac{\pi b}{2}s_2}(e^{\pi bs_1-2\pi bp_{s_1}}+e^{-\pi bs_1-2\pi b\l-2\pi bp_{s_1}}).$$
Now we make the following changes of variables (unitary transformations):
\begin{eqnarray*}
\circ s_1\mapsto -s_1 &:& p_{s_1}\mapsto -p_{s_1},\\
\circ e^{\pi i s_1\l} &:& p_{s_1}\mapsto p_{s_1}-\frac{\l}{2},\tab p_{\l}-\frac{s_1}{2},\\
\circ e^{-\frac{\pi i}{2}s_1s_2}&:& p_{s_1}\mapsto p_{s_1}-\frac{s_2}{4},\tab p_{s_2}\mapsto p_{s_2}-\frac{s_1}{4}.\\
\end{eqnarray*}
In other words, the transformation:
$$f(s_1,t_1)g(\l,t_2)\mapsto f(-s_1,t_1)g(\l,t_2)e^{\frac{\pi i}{2}s_1(2\l-s_2)}.$$

Then the representation (denoted by $\cP_{\l,t_1}$) becomes:
\begin{eqnarray*}
K_0 &=& e^{\frac{\pi b}{2}s_2},\\
K&=& e^{-\pi bs_1},\\
E&=&\frac{i}{q-q\inv}(e^{\pi b(-2p_{s_1}-s_1+\l)}+e^{\pi b(-2p_{s_1}+s_1-\l)})\\
&=&\frac{i}{q-q\inv}(q^{1/2}e^{\pi b(s_1-\l)}+q^{-1/2}e^{-\pi b(s_1-\l)})e^{-2\pi bp_{s_1}}\\
&=&[\frac{Q}{2b}-\frac{i}{b}(s_1-\l)]_q e^{-2\pi bp_{s_1}},\\
F&=&\frac{i}{q-q\inv}(e^{\pi b(2p_{s_1}-s_1-\l)}+e^{\pi b(2p_{s_1}+s_1+\l)})\\
&=&\frac{i}{q-q\inv}(q^{-1/2}e^{\pi b(s_1+\l)}+q^{1/2}e^{-\pi b(s_1+\l)})e^{2\pi bp_{s_1}}\\
&=&[\frac{Q}{2b}+\frac{i}{b}(s_1+\l)]_q e^{2\pi bp_{s_1}}.\\
\end{eqnarray*}
\end{proof}
\begin{Rem} We note that there is a freedom of choice in multiplication by $G_b$ in the above action. In fact this precisely says that
\begin{eqnarray}
\cP_{\l,t}&\simeq& \cP_{-\l,t}\\
f(x)&\mapsto& \frac{G_b(\frac{Q}{2}+i\l-ix)}{G_b(\frac{Q}{2}-i\l-ix)}f(x)
\end{eqnarray}
as observed in \cite{PT2}.
\end{Rem}

Since the pairing is non-degenerate, we obtain the following important realization of the multiplicative unitary $\bW$:
\begin{Cor}\label{Wdecompose} As a corepresentation on $L^2(\R)$, $\bW$ is equivalent to
\Eq{\bW=\int_\R\int_{\R_+}^{\o+} T^{\l,t/2} d\mu(\l)dt,}
where $d\mu(\l)=|S_b(Q+2i\l)|^2 d\l$.
\end{Cor}

This is a generalization of Podle\'{s}-Woronowicz's definition \cite{PW} of a multiplicative unitary for a compact group
\Eq{u = {\sum_{\a\in \hat{G}}}^{\o+} u^\a,}
where $\hat{G}$ is the set of all (finite dimensional) irreducible unitary representation of $G$, and 
$$u^\a:V_\a\mapsto V_\a\ox G$$
is the irreducible corepresentation on some finite dimensional vector space $V_\a$. In particular, this defines the multiplicative unitary in a coordinate free canonical way.

Furthermore, by applying the transformations that was done in the proof of Theorem \ref{ThmLeft} to the expression \eqref{leftcorep}, and comparing the result with the expression of $T^{\l,t/2}$ given in \eqref{fundcorep}, we recover the ``6-9 relation" of $G_b$ that is first observed by Volkov \cite{Vo}. See Section \ref{sec:meaning} for further details.

\subsection{Right regular representation}\label{sec:regular:right}
Similarly by applying the pairing to the right regular corepresentation we obtain
\begin{Thm}\label{ThmRight} The right regular corepresentation is equivalent to $$ L^2(\R, ds_1)\ox\int_\R\int_{\R^+}^\o+ \cP_{\l,-s_2/2}|S_b(Q+2i\l)|^2d\l ds_2.$$
\end{Thm}
\begin{proof} We apply the pairing to \eqref{rightcorep} and get
{
\allowdisplaybreaks
\begin{eqnarray*}
K_0&=& e^{-\pi b\frac{s_2}{2}},\\
K&=&e^{\pi bt_2},\\
E&=&\frac{(-iq^{1/2})(-q^2)e^{-\pi iQb}}{1-q^2} e^{\frac{\pi b}{2}s_2}e^{-\pi bt_2}\left( e^{\pi b(2t_1+2t_2)}e^{2\pi b(-p_{t_1}+p_{t_2})}+\right.\\
&&\left. \frac{G_b(\frac{Q}{2}-it_1-it_2)G_b(\frac{Q}{2}+it_1-it_2)}{G_b(\frac{Q}{2}-it_1-it_2-b)G_b(\frac{Q}{2}+it_1-it_2-b)}e^{2\pi bp_{t_2}}\right)\\
&=&\frac{i}{q-q\inv}e^{\frac{\pi b}{2}s_2}e^{-\pi bt_2}\left(e^{\pi b(-2p_{t_1}+2t_1+2p_{t_2}+t_2)}+\right.\\
&&\left.q^{1/2}(1+q\inv e^{2\pi b(t_1+t_2)})(1+q\inv e^{2\pi b(-t_1+t_2)})e^{2\pi b p_{t_2}}\right)\\
&=&\frac{i}{q-q\inv}e^{\frac{\pi b}{2}s_2}\left(e^{\pi b(2p_{t_2}-t_2)}+e^{\pi b(2p_{t_2}+3t_2)}+(e^{2\pi b t_1}+e^{-2\pi bt_1}+e^{\pi b(-2p_{t_1}+2t_1)})e^{\pi b(2p_{t_2}+t_2)}\right),\\\\
F&=&\frac{iq^{-1/2}e^{\pi iQb}}{1-q^2}e^{\frac{\pi b}{2}s_2+\pi bt_2}e^{-2\pi bt_2-\pi bs_2}e^{-2\pi bp_{t_2}}\\
&=&\frac{i}{q-q\inv}e^{-\frac{\pi b}{2}s_2}e^{-\pi b(t_2+2p_{t_2})}.\\
\end{eqnarray*}
}

Again by shifting $p_{t_1}\mapsto p_{t_1}+t_1$ as in the left representation, the action of $E$ becomes
$$E=\frac{i}{q-q\inv}e^{\frac{\pi b}{2}s_2}\left(e^{\pi b(2p_{t_2}-t_2)}+e^{\pi b(2p_{t_2}+3t_2)}+\bC e^{\pi b(2p_{t_2}-t_2)}\right).$$

Hence under the unitary transformation $\Phi$, $E$ becomes
$$E=\frac{i}{q-q\inv}e^{\frac{\pi b}{2}s_2}\left(e^{\pi b(2p_{t_2}-t_2)}+e^{\pi b(2p_{t_2}+3t_2)}+(e^{2\pi b\l}+e^{-2\pi b\l}) e^{\pi b(2p_{t_2}-t_2)}\right).$$

Finally, the remaining transformation used in the left regular representation doesn't change these actions.
In exactly the same way if we apply:
\begin{eqnarray*}
multiplication&by& G_b(\frac{Q}{2}+i\l-it_2)\inv,\\
\circ t_2\mapsto -t_2&:& t_2\mapsto -t_2, p_{t_2}\mapsto -p_{t_2},\\
\circ e^{-\pi it_2 \l}&:& p_{\l}\mapsto p_{\l}-\frac{t_2}{2},p_{t_2}\mapsto p_{t_2}-\frac{\l}{2},\\
\circ e^{-\frac{\pi i}{2} s_2t_2}&:&p_{s_2}\mapsto p_{s_2}+\frac{t_2}{4},p_{t_2}\mapsto p_{t_2}+\frac{s_2}{4},\\
\end{eqnarray*}
that does not change the left representation either, it becomes:
\begin{eqnarray*}
K_0&=& e^{-\pi b\frac{s_2}{2}}\\
K&=&e^{-\pi bt_2}\\
E&=&\frac{i}{q-q\inv}(q^{\frac{1}{2}}e^{\pi b(t_2-\l)}+q^{-\frac{1}{2}}e^{-\pi b(t_2-\l)})e^{-2\pi bp_{t_2}}\\
&=&[\frac{Q}{2b}-\frac{i}{b}(t_2-\l)]_q e^{-2\pi bp_{t_2}}\\
F&=&\frac{i}{q-q\inv}(q^{-\frac{1}{2}}e^{\pi b(t_2+\l)}+q^{\frac{1}{2}}e^{-\pi b(t_2+\l)})e^{2\pi bp_{t_2}}\\
&=&[\frac{Q}{2b}+\frac{i}{b}(t_2+\l)]_q e^{2\pi bp_{t_2}}\\
\end{eqnarray*}
which is precisely the representation denoted by $\cP_{\l, -s_2/2}$.
\end{proof}

Combining Theorem \ref{ThmLeft} and \ref{ThmRight}, and note that all the transformations used are unitary, we conclude the main theorem:
\begin{Thm} \label{main} As a representation of $U_q(\gl(2,\R))_L \ox U_q(\gl(2,\R))_R$, 
\Eq{L^2(GL_q^+(2,\R))\simeq \int_\R^\o+\int_{\R^+}^\o+ \cP_{\l,s}\ox \cP_{\l,-s} |S_b(Q+2i\l)|^2 d\l ds}
and the equivalence is unitary.
\end{Thm}
In particular, we proved Ponsot-Teschner's claim in \cite{PT1} for $U_q(\sl(2,\R))$ by setting $s_2=0$ which corresponds to the determinant element $N=A\hat{A}$, as well as the action of the central element $K_0$.

\subsection{The modular double}\label{sec:regular:MD}
It is known that the representation space $\cP_{\l,t}$ is actually a nontrivial representation for the modular double $U_q(\gl(2,\R))\ox U_{\til[q]}(\gl(2,\R))$. In this section we describe its action using the same multiplicative unitary $\bW_m$.

The modular double counterpart $U_{\til[q]}(\gl(2,\R))$ is defined in the same way with $q$ replaced by $\til[q]$ throughout.  In particular we can define the pairing between $\til[E],\til[F],\til[K],\til[K_0]$ with $\til[A],\til[B],\hat{\til[A]},\hat{\til[B]}$ as before. Recall that $\til[A]=A^{1/b^2}$ and similarly for the other quantum plane variables. Hence by the $b\longleftrightarrow b\inv$ duality of $\bW_m$, we can pair it with $U_{\til[q]}(\gl(2,\R))$, and obtain:

\begin{Prop} The left regular action of $U_{\til[q]}(\gl(2,\R))$ on $L^2(GL_q^+(2,\R))$ is given by replacing $b$ with $b\inv$, i.e.
\begin{eqnarray*}
\til[K_0] &=& e^{\frac{\pi b\inv}{2}s_2},\\
\til[K]&=& e^{-\pi b\inv s_1},\\
\til[E]&=&\frac{i}{\til[q]-\til[q]\inv}(\til[q]^{1/2}e^{\pi b\inv (s_1-\l)}+\til[q]^{-1/2}e^{-\pi b\inv(s_1-\l)})e^{-2\pi b\inv p_{s_1}},\\
\til[F]&=&\frac{i}{\til[q]-\til[q]\inv}(\til[q]^{-1/2}e^{\pi b\inv (s_1+\l)}+\til[q]^{1/2}e^{-\pi b\inv(s_1+\l)})e^{2\pi b\inv p_{s_1}}.\\
\end{eqnarray*}
We note that $\til[E]$ and $\til[F]$ are positive operators only when $\frac{1}{2n+1}<b^2<\frac{1}{2n}$ for $n\in \N$, and they are not necessarily essentially self-adjoint on the natural domain $\cW$ \cite{Ru}. To obtain essential self-adjointness, we have to apply Proposition \ref{x+p} and define the operators on the respective transformed domain $g_b\cdot \cW$.
\end{Prop}
Ignoring the factors, we note from Theorem \ref{ThmLeft} that the summand of $E$ and $F$ are $q$-commuting, hence we immediately obtain using Corollary \ref{u+v} (see also \cite[Cor 1]{BT}).
\begin{Prop}
As positive self-adjoint operators we have
\begin{eqnarray*}
K_0^{1/b^2}&=&\til[K_0],\\
K^{1/b^2}&=&\til[K],\\
2\sin(\pi b^{-2})\til[E] &=& (2\sin(\pi b^{2})E)^{1/b^2},\\
2\sin(\pi b^{-2})\til[F] &=& (2\sin(\pi b^{2})F)^{1/b^2}.
\end{eqnarray*}
\end{Prop}

Similar analysis works as well for the right regular representation, hence Theorem \ref{main} is actually an equivalence as a representation of the modular double $$U_{q\til[q]}(\gl(2,\R))_L\ox U_{q\til[q]}(\gl(2,\R))_R,$$ 
where $U_{q\til[q]}(\gl(2,\R)):=U_q(\gl(2,\R))\ox U_{\til[q]}(\gl(2,\R))$ denotes the modular double.
\section{Representation theoretic meaning for certain integral transforms}\label{sec:meaning}
In this section, we state without proof certain integral transformations of $G_b$ that arise in the calculations of certain representation relations.
\begin{Prop} The pentagon equation $W_{23}W_{12}=W_{12}W_{13}W_{23}$ is equivalent to the 4-5 relation (Lemma \ref{45})
\Eq{\frac{G_b(\a)G_b(\b)G_b(\c)}{G_b(\a+\c)G_b(\b+\c)}=\int_C d\t e^{-2\pi \c\t}\frac{G_b(\a+i\t)G_b(\b+i\t)}{G_b(\a+\b+\c+i\t)G_b(Q+i\t)}.}
By a change of variables and using the reflection formula, it can be rewritten as
\Eq{\frac{G_b(\a+\c)G_b(\b+\c)G_b(\a)G_b(\b)}{G_b(\a+\b+\c)}=\int_C d\t e^{-2\pi i(\b+i\t)(\a+i\t)}G_b(\a+i\t)G_b(\b+i\t)G_b(\c-i\t)G_b(-i\t).}
Scaling all the variables by $b$ and taking the limit $b\to 0$ by applying Theorem \ref{classicallimit}, we recover precisely Barnes' first lemma:
\Eq{\frac{\G(a+c)\G(b+c)\G(a)\G(b)}{\G(a+b+c)}=\frac{1}{2\pi }\int_C\G(a+i\t)\G(b+i\t)\G(c-i\t)\G(-i\t) d\t} in the special case for $d=0$.
\end{Prop}

Next, as noted in Section \ref{sec:regular:left}, by comparing $\bW_m=\iint^{\o+} T^{\l,t/2} dt d\mu(\l)$ as corepresentation, we obtain the following relation that is first observed by Volkov \cite{Vo}:
\begin{Prop}\label{69}
The 6-9 relation for $G_b(x)$ can be written as
\begin{eqnarray}
&&\int_C e^{-2\pi \t (\d-i\t)} \frac{G_b(\a+i\t)G_b(\b+i\t)G_b(\c+i\t)G_b(\d-i\t)G_b(-i\t)}{G_b(\a+\b+\c+\d+i\t)} d\t\nonumber\\
&=&\frac{G_b(\a)G_b(\b)G_b(\c)G_b(\a+\d)G_b(\b+\d)G_b(\c+\d)}{G_b(\a+\b+\d)G_b(\a+\c+\d)G_b(\b+\c+\d)},
\end{eqnarray}where the contour goes along $\R$ and separates the increasing and decreasing sequence of poles. By the asymptotic properties of $G_b$, the integral converges for any choice of parameters.

Again by scaling all the variables by $b$ and applying Theorem \ref{classicallimit}, we recover Barnes' second lemma:
\begin{eqnarray}
&&\int_C \frac{\G(a+i\t)\G(b+i\t)\G(c+i\t)\G(d-i\t)\G(-i\t)}{\G(a+b+c+d+i\t)}d\t\nonumber\\
&=&\frac{\G(a)\G(b)\G(c)\G(a+d)\G(b+d)\G(c+d)}{\G(a+b+d)\G(a+c+d)\G(b+c+d)},
\end{eqnarray}
which in turn is a generalization of Pfaff-Saalsch\"{u}tz's sum
\Eq{\sum_{j\geq 0} \frac{(n+m+l+k-j)!}{(m-j)!(l-j)!(k-j)!(n+j)!j!}=\frac{(n+m+l)!(n+m+k)!(n+l+k)!}{m!l!k!(n+m)!(n+l)!(n+k)!}.}
\end{Prop}

Finally, in \cite{WZ} an alternative description of the multiplicative unitary lying instead in $\cA\ox \cA$ is defined to be (slightly modified to fit our definition):
\Eq{V=g_b(B\inv\ox q\inv BA\inv)e^{\frac{i}{2\pi b^2}\log (qAB\inv)\ox \log A\inv},}
so that
\Eq{V(x\ox 1)V^* = \D(x).}
On comparing
\Eq{W^*(1\ox x)W = \D(x) = V(x\ox 1)V^*}
as operators on $\cH\ox\cH$, we obtain a new relation involving $G_b(x)$, which is essentially the same as the relation in \cite[Theorem 5.6.7]{VdB} in the case $n=1$:
\begin{Prop}\label{32} The 3-2 relation is given by
\Eq{\int_C G_b(\a+i\t)G_b(\b-i\t)G_b(\c-i\t)e^{-2\pi i(\b-i\t)(\c-i\t)}d\t = G_b(\a+\c)G_b(\a+\b),}
where the contour goes along $\R$ and separates the poles for $i\t$ and $-i\t$. By the asymptotic properties for $G_b$, the integral converges for $Re(\a-\b-\c)<\frac{Q}{2}$.
\end{Prop}


\addcontentsline{toc}{section}{References}
\begin{thebibliography}{99}


  \bibitem{BV}
  S. Baaj and S. Vaes,
  \textit{Double crossed products of locally compact quantum groups},
  J. Inst. Math. Jussieu \textbf{4}, 135-173, (2005)
\bibitem{BT}
  A.G. Bytsko, K. Teschner,
  \textit{R-operator, co-product and Haar-measure for the modular double of $U_q(\sl(2,\R))$},
  Comm. Math. Phys. \textbf{240}, 171-196, (2003)
  \bibitem{CF}
  L. Chekhov and V.V. Fock,
  \textit{A quantum Teichm\"{u}ller space},
  Theor. Math. Phys. \textbf{120}, 511-528, (1999)
  \bibitem{DMMZ}
  E. E. Demidov, Yu. I. Manin, E. E. Mukhin and D. V. Zhdanovich,
  \textit{Non-standard quantum deformations of $GL(n)$ and constant solutions of the Yang-Baxter equation},
  Prog. Theor. Phys. Supp No.102, 203-218, (1990) 
\bibitem{Fa}
  L.D. Faddeev,
  \textit{Discrete Heisenberg-Weyl group and modular group},
  Lett. Math. Phys. \textbf{34}, 249-254, (1995)  
\bibitem{FV}
L.D. Faddeev, A. Yu. Volkov,
\textit{Abelian current algebra and the Virasoro algebra on the lattice},
Phys.Lett. B315, \textbf{3,4}, 311-318, (1993)
\bibitem{FG}
  V.V. Fock, A.B. Goncharov,
  \textit{The quantum dilogarithm and representations of the quantum cluster varieties},
  Inventions Math. \textbf{175}, 223-286, (2009)
\bibitem{F}
  I. Frenkel
  \textit{Lectures on quaternionic analysis and representation theory}, Yale University (2006), unpublished
  \bibitem{FJ}
  I. Frenkel, M. Jardim,
  \textit{Quantum instantons with classical moduli spaces},
  Comm. Math. Phys., \textbf{237}(3), 471-505, (2003)
  \bibitem{FK}
  I. Frenkel, H. Kim,
  \textit{Quantum Teichm\"{u}ller space from quantum plane},
  arXiv:1006.3895v1, (2010)
\bibitem{FL}
  I. Frenkel, M. Libine, 
 \textit{Quaternionic analysis, representation theory and physics}, 
 Advances in Math \textbf{218}, 1806-1877, (2008)
  \bibitem{Gon}
  A. B. Goncharov,
\textit{Pentagon relation for the quantum dilogarithm and quantized $\cM_{0,5}^{cyc}$},
Progress in Mathematics, \textbf{256}, 413-426, (2007)
\bibitem{GKK}
W. Groenevelt, E. Koelink, J. Kustermans,
\textit{The dual quantum group for the quantum group analog of the normalizer of $SU(1,1)$ in $SL(2,\C)$},
Int. Math. Res. Not. \textbf{7}, 1167-1314, (2010)
  \bibitem{Ip}
  I. Ip,
  \textit{The classical limit of representation theory of the quantum plane},
   arXiv:1012.4145v1 [math.RT], (2010)
  \bibitem{Ip2}
  I. Ip,
  \textit{The graphs of quantum dilogarithm},
   arXiv:1108.5376v1 [math.QA], (2011)
  \bibitem{K}
  Byung-Jay Kahng,
  \textit{Twisting of the quantum double and the Weyl algebra},
  arXiv:0809.0098v1 [math.OA], (2008)
  \bibitem{Ka1}
  R.M. Kashaev
  \textit{The Heisenberg double and the pentagon relation},
  St. Petersburg Math. J. \textbf{8}, 585-592, (1997)  
\bibitem{Ka2}
  R.M. Kashaev,
  \textit{The hyperbolic volume of knots from the quantum dilogarithm},
  Lett. Math. Phys. \textbf{39}, 269-275, (1997)
\bibitem{Ka3}
  R.M. Kashaev,
  \textit{Quantization of Teichm\"{u}ller spaces and the quantum dilogarithm},
  Lett. Math. Phys. \textbf{43}, 105-115, (1998)
\bibitem{Ka4}
R.M. Kashaev,
\textit{The quantum dilogarithm and Dehn twist in quantum Teichm\"{u}ller theory},
Integrable Structures of Exactly Solvable Two-Dimensional Models of Quantum Field Theory (Kiev, Ukraine, September 25-30, 2000), NATO Sci. Ser. II Math. Phys. Chem., vol. 35, Kluwer, Dordrecht, 211-221 (2001)
\bibitem{KaN}
  R.M. Kashaev, T. Nakanishi,
 \textit{Classical and quantum dilogarithm identities},
 arXiv:1104.4630v2 [math.QA], (2011)  
\bibitem{KLS}
  S. Kharchev, D. Lebedev, M. Semenov-Tian-Shanksy,
  \textit{Unitary representations of $U_q(\sl(2,\R))$, the modular double, and the multiparticle $q$-deformed Toda chain},
  Comm. Math. Phys. \textbf{225}(3), 573-609, (2002)
\bibitem{KK}
E. Koelink, J. Kustermans,
\textit{A locally compact quantum group analogue of the normalizer of $SU(1,1)$ in $SL(2,\C)$},
Comm. Math. Phys. \textbf{233} no. 2, 231-296, (2003)
  \bibitem{KV1}
  J. Kustermans and S. Vaes,
  \textit{Locally compact quantum groups},
  Ann. Sci. Ecole Norm. Sup. (4) \textbf{33}, 837-934, (2000)
  \bibitem{KV2}
  J. Kustermans and S. Vaes,
  \textit{Locally compact quantum groups in the von Neumann algebraic setting},
  Math. Scand. \textbf{92}, 68-92, (2003)
  \bibitem{Ma}
  S. Majid
  \textit{Foundations of quantum group theory},
  Cambridge University Press, Cambridge, (1995)
   \bibitem{MMNNSU}
  T. Masuda, K. Mimachi, Y. Nakagami, M. Noumi, Y. Saburi, K. Ueno,
  \textit{Unitary representations of the quantum group $SU_q(1,1)$: I, II},
  Letters in Math. Phys. 19, 187-204, (1990)  
 \bibitem{MMNNU}
  T. Masuda, K. Mimachi, Y. Nakagami, M. Noumi, K. Ueno,
  \textit{Representations of the quantum group $SU_q(2)$ and the little $q$-Jacobi polynomials},
  Journal of Functional Analysis 99, 357-386, (1991)
\bibitem{PK}
  R.B. Paris, D. Kaminski,
  \textit{Asymptotics and Mellin-Barnes integrals},
  Encyclopedia of Mathematics and its Applications,
  \textbf{85}, Cambridge University Press, Cambridge,  (2001)
\bibitem{PW}
  P. Podle\'{s}, S.L. Woronowicz,
  \textit{Quantum deformation of Lorentz group},
  Commun Math. Phys. 130, 381-431 
  (1990)
\bibitem{PT1}
  B. Ponsot, J. Teschner,
  \textit{Liouville bootstrap via harmonic analysis on a noncompact quantum group},
  arXiv: hep-th/9911110, (1999)
  \bibitem{PT2}
  B. Ponsot, J. Teschner,
  \textit{Clebsch-Gordan and Racah-Wigner coefficients for a continuous series of representations of $\cU_q(\mathfrak{sl}(2,\R))$},
  Comm. Math. Phys \textbf{224}, 613-655, (2001)

\bibitem{Puk}
  L. Puk\'{a}nszky
  \textit{On the Kronecker products of irreducible representations of the $2\x 2$ real unimodular group. I},
  Trans. Amer. Math. Soc., Vol. 100, No. 1, 116-152, (1961)
\bibitem{Pu}
  W. Pusz,
  \textit{Quantum $GL(2,\C)$ group as double group over $'az+b'$ quantum group},
  Reports on Mathematical Physics, Vol. 49, 113-122,  (2002)
\bibitem{PuW}
  W. Pusz, S.L. Woronowicz,
  \textit{A new quantum deformation of '$ax+b$' group},
  Commun Math. Phys. 259, 325-362, (2005)
\bibitem{Ru}
  S.N.M. Ruijsenaars,
  \textit{A unitary joint eigenfunction transform for the $A\D O$'s $\exp(ia_\pm d/dz)+\exp(2\pi z/a_\mp)$},
  J. Nonlinear Math. Phys. \textbf{12} Suppl. 2, 253-294, (2005)
\bibitem{SWZ}
A. Schirrmacher, J. Wess and B. Zumino,
\textit{The two-parameter deformation of $GL(2)$, its differential calculus, and Lie algebra},
Z. Phys. C - Particles and Fields \textbf{49}, 317-324, (1991)
\bibitem{Sch}
K. Schm\"{u}dgen,
\textit{Operator representations of $\R_q^2$},
Publ. RIMS Kyoto Univ. \textbf{29}, 1030-1061, (1993)
\bibitem{Ta}
 M. Takesaki,
\textit{Tomita's theory of modular Hilbert algebras and its applications},
Lecture Notes in Math. 128, Springer-Verlag, Berlin, (1970)
\bibitem{T}
T. Timmermann,
\textit{An invitation to quantum groups and duality},
EMS Textbooks in Mathematics, (2008)
\bibitem{VD}
A. van Daele,
\textit{Multiplier Hopf algebras},
Trans. Amer. Math. Soc. \textbf{342} 917-932, (1994) 
\bibitem{VdB}
F. J. van de Bult, 
\textit{Hyperbolic hypergeometric functions}, 
Ph. D. thesis, University of Amsterdam, (2007)
 \bibitem{Vi}
  N. Vilenkin,
  \textit{Special functions and the theory of group representations},
  Monographs. \textbf{22}, Amer. Math. Soc., Providence, RI,  (1968)
\bibitem{Vo}
  A. Yu. Volkov,
  \textit{Noncommutative hypergeometry},
  Comm. Math. Phys. \textbf{258}(2), 257-273, (2005)
\bibitem{W}
  S.L. Woronowicz,
  \textit{From multiplicative unitaries to quantum groups},
  Int. J. Math. \textbf{7}(1), 127-149, (1996)
\bibitem{WZ}
  S.L. Woronowicz, S. Zakrzewski,
  \textit{Quantum '$ax+b$' group},
  Rev. Math. Phys. \textbf{14}(7\& 8), 797-828, (2002)
\end{thebibliography}
\end{document}